\documentclass[Unicode,OJMO,english,manuscript,NoThmDefs]{cedram}

\usepackage[all]{xypic}

\usepackage{csquotes}  % automatically chooses the right quotation marks
\usepackage{acro}      % convenient management of acronyms
\usepackage{nicefrac}  % used for inline fractions or fractions in matrices
\usepackage[capitalise]{cleveref}  % care-free references

\usepackage{mathtools} % \mathclap and \shortintertext

\usepackage{tabularx}
\usepackage{booktabs}

\theoremstyle{plain}
\newtheorem{assumption}{Assumption}
\crefname{assumption}{Assumption}{Assumptions}

\newtheorem{mytheorem}{Theorem}
\crefname{mytheorem}{Theorem}{Theorems}

\newtheorem{mylemma}{Lemma}
\crefname{mylemma}{Lemma}{Lemmata}

\newtheorem{mycorollary}{Corollary}
\crefname{mycorollary}{Corollary}{Corollaries}

\theoremstyle{remark}
\newtheorem*{myremark*}{Remark}
\newtheorem{myremark}{Remark}
\crefname{myremark}{Remark}{Remarks}

\theoremstyle{definition}
\newtheorem{mydefinition}{Definition}
\crefname{mydefinition}{Definition}{Definitions}

\crefname{step}{step}{steps}

\newcommand{\CC}{\mathtt{c}} % generic constant

\newcommand*{\tangcone}[2]{\mathcal T_{#1}\left(#2\right)}

\newcommand{\abs}[1]{\left| #1 \right|}
\newcommand*{\norm}[1]{\left\|#1\right\|} 

\newcommand*{\ve}[1]{\pmb{#1}} % vector shorthand (makes letter bold)
\newcommand*{\mat}[1]{\ve{#1}} % matrix shorthand (same as \ve)

\newcommand*{\eye}[1]{\mat{I}_{#1,#1}} % \eye{n} Identity matrix of dimension n

\newcommand{\grad}{\ve{\nabla}\!}
\newcommand{\jac}{\mat{D}\!}
\newcommand{\hess}{\mat{H}\!}

\newcommand*{\tran}[1]{#1^{\mkern-1.5mu\mathsf{T}}}

\newcommand*{\kk}{k}  % letter used as iteration index

\newcommand*{\iterUpCmd}[1]{\ensuremath{^{(#1)}}}

\newcommand*{\iterDownCmd}[1]{\ensuremath{_{(#1)}}}
\newcommand*{\iterDown}{\iterDownCmd{\kk}}
\newcommand*{\iterSubCmd}[1]{\ensuremath{_{#1}}}
\newcommand*{\iterSub}{\iterSubCmd{\kk}}

\newcommand*{\itCmd}[1]{\iterUpCmd{#1}}
\newcommand*{\itDeltaCmd}[1]{\iterDownCmd{#1}}
\newcommand*{\itthetaCmd}[1]{\iterSubCmd{#1}}

\newcommand*{\KK}{\itCmd{\kk}}  % shorthand to typeset iteration index

\newcommand*{\normTr}[1]{\norm{#1}_{\mathrm{tr}}}
\newcommand*{\normTrIt}[1]{\norm{#1}_{\mathrm{tr,\kk}}}
\newcommand*{\normTwo}[1]{\norm{#1}_2}
\newcommand*{\normInf}[1]{\norm{#1}_\infty}
\newcommand*{\normIt}[1]{\norm{#1}\iterSub}

\newcommand*{\feas}{\mathcal X} % feasible set
\newcommand*{\vx}{\ve{x}}       % vector x
\newcommand*{\vZ}{\ve{0}}       % vector zero

\newcommand*{\vd}{\ve{d}}
\newcommand*{\vs}{\ve{s}}

\newcommand*{\vn}{\ve{n}}
\newcommand*{\vdk}{\vd\KK}
\newcommand*{\vnk}{\vn\KK}
\newcommand*{\vsk}{\vs\KK}

\newcommand*{\vxNk}{\vxk_{\mathrm{n}}}
\newcommand*{\vxN}{\vx_{\mathrm{n}}}
\newcommand*{\vxSk}{\vxk_{\mathrm{s}}}
\newcommand*{\vxS}{\vx_{\mathrm{s}}}

\newcommand{\vxk}{\vx\KK}

\newcommand*{\fL}{f}
\newcommand*{\gL}{g}
\newcommand*{\hL}{h}

\newcommand*{\surrogateSign}{\hat}

\newcommand*{\vf}{\ve{\fL}}
\newcommand*{\vh}{\ve{\hL}}%{\ve{c}_{\mathrm{e}}}
\newcommand*{\vg}{\ve{\gL}}%{\ve{c}_{\mathrm{e}}}

\newcommand{\df}{\grad{\fL}}
\newcommand{\dg}{\grad{\gL}}
\renewcommand*{\dh}{\grad{\hL}}

\newcommand*{\DF}{\mat{\MakeUppercase{\fL}}}
\renewcommand*{\DH}{\mat{\MakeUppercase{\hL}}}
\newcommand*{\DG}{\mat{\MakeUppercase{\gL}}}

\newcommand*{\DFk}{\DF\iterSub}

\newcommand*{\varDim}{n}
\newcommand*{\fDim}{K}
\newcommand*{\hDim}{M}
\newcommand*{\gDim}{P}

\newcommand*{\mfL}{\surrogateSign{\fL}}
\newcommand*{\mhL}{\surrogateSign{\hL}}
\newcommand*{\mgL}{\surrogateSign{\gL}}

\newcommand*{\mf}{\ve{\surrogateSign{\fL}}}
\newcommand*{\mg}{\ve{\surrogateSign{\gL}}}
\newcommand*{\mh}{\ve{\surrogateSign{\hL}}}

\newcommand*{\mfk}{\mf\KK}
\newcommand*{\mgk}{\mg\KK}
\newcommand*{\mhk}{\mh\KK}

\newcommand*{\mgkSub}{\mg\iterSub}

\newcommand{\dmf}{\grad{\mf}}

\newcommand{\dmh}{\grad{\mh}}

\newcommand*{\DmF}{\mat{\surrogateSign{\MakeUppercase{\fL}}}\iterSub}
\newcommand*{\DmH}{\mat{\surrogateSign{\MakeUppercase{\hL}}}\iterSub}
\newcommand*{\DmG}{\mat{\surrogateSign{\MakeUppercase{\gL}}}\iterSub}

\newcommand*{\scal}[1]{\Phi\bigl[ #1 \bigr]}
\newcommand*{\scalf}[1]{\scal{\vf}\left(#1\right)}
\newcommand*{\scalShort}{\Phi} 
\newcommand*{\scalk}{\scalShort\KK}

\newcommand{\stepsize}{\sigma}
\newcommand{\stepsizek}{\stepsize\iterDown}

\newcommand*{\linFeasLetter}{L}
\newcommand*{\linFeas}{\mathcal{\linFeasLetter}}
\newcommand*{\linFeasApprox}{\mathcal{\linFeasLetter}}
\newcommand*{\linFeasApproxk}{\mathcal{\linFeasLetter}\iterSub}

\newcommand*{\omegak}{\surrogateSign{{\omega}}\KK}

\newcommand*{\chik}{{\chi}\KK}
\newcommand*{\chikPre}{\bar{{\chi}}\KK}

\newcommand*{\Deltak}{\Delta\itDeltaCmd{\kk}}
\newcommand*{\DeltaMax}{\Delta_{\text{max}}}
\newcommand*{\DeltaMin}{\Delta_{\min}}

\newcommand{\DeltaPre}{\bar{\Delta}}
\newcommand{\DeltaPrek}{\DeltaPre\itDeltaCmd{\kk}} 

\newcommand*{\ballLetter}{B}
\newcommand*{\ball}[2]{\ballLetter\left(#1;#2\right)}
\newcommand{\trk}{\ballLetter\KK} 

\newcommand*{\thetak}{{\theta}\itthetaCmd\kk}

\newcommand*{\covering}{\mathcal C(\mathcal X)}
\newcommand{\filter}{\mathcal F}
\newcommand*{\filterConst}{\gamma_{{\theta}}}

\newcommand*{\CCDelta}{\CC_{\Delta}}
\newcommand*{\CCmu}{\CC_{\mu}}

\newcommand*{\CCdeltan}{\delta_{\mathrm{n}}}
\newcommand*{\CCubn}{\CC_{\mathrm{ubn}}}

\newcommand*{\CCfl}{\mathtt{e}}
\newcommand*{\CCflDiff}{\dot{\CCfl}}
\newcommand*{\LipschitzBoundModel}[1]{L^{#1}}
\newcommand*{\LipschitzBoundModels}{L^{\mfL}}

\newcommand*{\CCflF}{\CCfl_{\fL}}
\newcommand*{\CCflG}{\CCfl_{\gL}}
\newcommand*{\CCflH}{\CCfl_{\hL}}
\newcommand*{\CCflDiffF}{\CCflDiff_{\fL}}
\newcommand*{\CCflDiffG}{\CCflDiff_{\gL}}
\newcommand*{\CCflDiffH}{\CCflDiff_{\hL}}

\newcommand*{\CCnorm}{\mathsf{c}}
\newcommand*{\CCcrit}{\mathtt{w}}

\newcommand*{\armijoConstant}{\mathsf{a}}
\newcommand*{\backtrackConstant}{\mathsf{b}}
\newcommand*{\CCsd}{\mathtt{c}_{\mathrm{sd}}}%{\kappa_{sd}^{\mfL}\,}
\newcommand*{\CCsdT}{\tilde{\mathtt{c}}_{\mathrm{sd}}}

\newcommand*{\CChessF}{\mathtt{H}}%_{\mfL}}

\newcommand*{\chikl}{{\chi}\itCmd{\kk_{\ell}}}
\newcommand*{\thetakl}{{\theta}\itthetaCmd{\kk_{\ell}}}
\newcommand*{\vnkl}{\vn\itCmd{\kk_{\ell}}}

\newcommand*{\vxkl}{\vx\itCmd{\kk_{\ell}}}
\newcommand*{\Deltakl}{\Delta\itDeltaCmd{\kk_{\ell}}}
\newcommand*{\vxNkl}{\vxN\itCmd{\kk_{\ell}}}
\newcommand*{\vxSkl}{\vxS\itCmd{\kk_{\ell}}}
\newcommand*{\scalkl}{\scalShort\itCmd{\kk_{\ell}}}
\newcommand*{\mfkl}{\mf\itCmd{\kk_{\ell}}}

\newcommand*{\linFeasApproxkl}{\mathcal{\linFeasLetter}\iterSubCmd{\kk_{\ell}}}

\newcommand*{\DeltaInit}{\Delta\itDeltaCmd{0}}
\newcommand*{\gammass}{\gamma_{0}}
\newcommand*{\gammas}{\gamma_{1}}
\newcommand*{\gammag}{\gamma_{2}}
\newcommand*{\CCaccept}{\nu_1}
\newcommand*{\CCsuccess}{\nu_0}
\newcommand*{\CCcritTestchi}{\varepsilon_{{\chi}}}
\newcommand*{\CCcritTesttheta}{\varepsilon_{{\theta}}}

\newcommand*{\CCMu}{\mathtt{M}}
\newcommand*{\CCBeta}{\mathtt{B}}
\newcommand*{\CCalpha}{\alpha}

\newcommand*{\CCcritDenom}{\mathtt{W}}

\newcommand{\vxNext}{\vx^{(\kk+1)}}

\newcommand{\DeltaPreNext}{\DeltaPre_{(\kk+1)}}
\newcommand{\None}{\mathtt{undef}}

\newcommand*{\acceptIndices}{\mathcal{A}}
\newcommand*{\restoreIndices}{\mathcal{R}}
\newcommand*{\filterIndices}{\mathcal{Z}}

\newcommand*{\CCubj}{\CC_{\mathrm{ubj}}}

\newcommand*{\CCaccN}{\mathtt{e}_{\text{err}}}
\newcommand*{\CCaccS}{\CCaccN}

\newcommand*{\CCthetaub}{\CC_{\mathrm{ub}{\theta}}}

\newcommand*{\CCepschi}{\epsilon}

\DeclareAcronym{mop}{
  short=MOP,
  long=multi-objective optimization problem,
}

\DeclareAcronym{moo}{
  short=MOO,
  long=multi-objective optimization,
}

\DeclareAcronym{rbf}{
  short=RBF,
  long=radial basis function,
}

\DeclareAcronym{wlog}{
  short=w.l.o.g.,
  long=without loss of generality,
}

\DeclareAcronym{lhs}{
  short=LHS,
  long=left-hand side,
}

\DeclareAcronym{rhs}{
  short=RHS,
  long=right-hand side,
}

\DeclareAcronym{mfcq}{
  short=MFCQ,
  long={{M}angasarian-{F}romovitz constraint qualification}
}
\title[Filter Trust-Region Algorithm for MOPs with Constraints]{Multi-Objective Trust-Region Filter Method for Nonlinear Constraints using Inexact Gradients}

\author{\firstname{Manuel} \lastname{Berkemeier}}
\address{Department of Computer Science\\
Paderborn University\\
Germany}
\email{manuelbb@math.upb.de}
\author{\firstname{Sebastian} \lastname{Peitz}}
\address{Department of Computer Science\\
  Paderborn University\\
  Germany}
\email{sebastian.peitz@upb.de}

\begin{abstract} 
  In this article, we build on previous work to present an optimization algorithm for non-linearly
  constrained multi-objective optimization problems.
  The algorithm combines a surrogate-assisted derivative-free trust-region approach with the
  filter method known from single-objective optimization.
  Instead of the true objective and constraint functions,
  so-called \emph{fully linear} models are employed, and we show how
  to deal with the gradient inexactness in the composite step setting,
  adapted from single-objective optimization as well.
  Under standard assumptions, we prove convergence of a subset of iterates to a quasi-stationary point and,
  if constraint qualifications hold, then the limit point is also a KKT-point of the multi-objective problem.
\end{abstract}

\begin{document}

% Use the \maketitle command after the abstract
\maketitle

\section{Introduction}
%\subsection{Setting and Literature Review}

\Acfp{mop} arise naturally in many areas of mathematics, engineering, in the natural sciences or
in economics.
The goal of \acf{moo} is to find acceptable trade-offs between the competing objectives of a
\ac{mop}.
Generally, there are multiple solutions constituting the so-called Pareto Set in variable space
and the Pareto Front in objective space.
In this article, we consider the problem
\begin{equation}
    \min_{\vx\in \mathbb{R}^\varDim}
      \begin{bmatrix}
        \fL_1( \vx)
        \\
        \vdots
        \\
        \fL_\fDim( \vx )
      \end{bmatrix}
      =
      \min_{\vx \in \mathbb{R}^\varDim} \vf(\vx)
    \quad\text{s.t.}\quad
    \vh(\vx) = \tran{[\hL_1 (\vx), \ldots, \hL_\hDim(\vx)]} = \vZ,
    \;
    \vg(\vx) = \tran{[\gL_1 (\vx), \ldots, \gL_\gDim(\vx)]} \le \vZ
    ,
  \tag{MOP}
  \label{eqn:mop_nonlin}
\end{equation}
where all functions are twice continuously differentiable.
A global Pareto-optimal point $\vx^*$ is feasible and non-dominated, i.e.,
there is no other feasible $\vx \ne \vx^*$ with $\vf(\vx) \le \vf(\vx^*)$ and
$\fL_{\ell}(\vx) < \fL_{\ell}(\vx^*)$ for some ${\ell}\in\{1,\ldots ,\fDim\}$.
Throughout this article we will refer to the feasible set as
$\feas :=
  \left\{
    \vx \in \mathbb{R}^\varDim :
    \vh(\vx) = \vZ, \,
    \vg(\vx) \le \vZ
  \right\}$.

There is a multitude of methods available to approximate a single solution or the entire
Pareto Set/Front (or a superset thereof)
and the choice of method is heavily dependent on the structure of the problem at hand and
the demands of the person seeking a solution.
The references%
~\cite{ehrgott_multicriteria_2005,%
  jahn_vector_2011,%
  miettinen_nonlinear_2013,%
  eichfelder_twenty_2021%
} all provide an extensive overview of the topic.
Amongst others, there are scalarization approaches~\cite{eichfelder_adaptive_2008}
and adaptions of single-objective descent
methods to the multi-objective case -- for different problem classes as defined by the
properties of their objectives and constraint functions, see%
~\cite{%
  fukuda_survey_2014,%
  fliege_svaiter,%
  grana_drummond_steepest_2005,%
  lucambio_perez_nonlinear_2018,%
  gebken_constraints,%
  wilppu_new_2014,%
  Gebken_Peitz_2021%
}
for examples of well-known scalar techniques adjusted for constrained and
unconstrained \acp{mop}
in the smooth and non-smooth case.
For global approximations there are, e.g., evolutionary algorithms
(see~\cite{coello2013evolutionary} for an overview and~\cite{deb2002fast} for a prominent example)
and structure-exploiting methods%
~\cite{%
  hillermeier_nonlinear_2001,%
  pareto_tracer,%
  gebken_hierarchical_2019%
}.
Of course, there is also research in combining global and local
techniques~\cite{%
  schutze_pareto_2020,%
  peitz_gradient-based_2018%
}.
Whilst there are natural applications of \ac{moo} for machine
learning tasks~\cite{sener_koltun,gan_training},
machine learning techniques can conversely be used to assist the search
for optimal points~\cite{pfront_training}.
In case of expensive objective or constraints, surrogate models can be employed~\cite{%
  peitz_survey_2018,%
  chugh_survey_2019,%
  deb_surrogate_2020%
}.

In our setting we assume (some) objective and
constraint functions to be computationally
expensive and without exact derivative-information available.
This motivates the use of \emph{derivative-free} optimization methods, which also have been
adapted to the multi-objective case.
Most prominently, there are direct search algorithms~\cite{%
  audet_mesh_2010%
}
and surrogate-assisted trust-region algorithms~\cite{%
  qu_goh_liang,%
  thomann_paper,%
  villacorta%
}.
Based on those trust-region algorithms, we have in a previous article~\cite{ours}
presented a trust-region algorithm for
problems with a feasible set
that is convex and compact (or $\feas = \mathbb{R}^\varDim$).
The algorithm uses fully linear models
(e.g., Lagrange polynomials or \acp{rbf} as in~\cite{wild_orbit_2008}) to approximate the
objective functions and passes the exact constraints to an inner solver.
We build upon this work to accommodate general non-linear constraints
by also modelling them and solving inexact sub-problems.

To this end, we transfer the techniques from \cite{fletcher_global_2002}
to the multi-objective case with inexact derivatives.
Inexact derivative have already been handled in single-objective optimization 
in a similar manner:
For example, the authors of~\cite{ferreira_global_2017} provide
strong convergence results for inexact objectives and objective derivatives with similar model
accuracy requirements to our case.
We will also discuss the single-objective algorithm presented in~\cite{eason_trust_2016}, 
as it is also based on fully linear models and employs a special kind of criticality 
check.
Our work, however, is more along the lines of~\cite{walther2016}, 
where the filter trust-region algorithm and 
the composite step framework are likewise modified
to handle inexact derivatives of both the objective and constraint functions.
But in~\cite{walther2016} the derivatives (approximated via automatic differentiation)
can eventually become exact, in contrast to our setting, which is why
--without constraint qualifications-- we can only prove convergence to a \emph{quasi-stationary}
point, similar to~\cite{echebest_inexact_2017}.
The algorithms in \cite{echebest_inexact_2017} and \cite{garreis_inexact_2019}
also use fully linear surrogate models satisfying the same error bounds 
as assumed in this article within single-objective filter algorithms.
In~\cite{echebest_inexact_2017}, convergence to KKT-points is proven under constraint qualifications.
We will see, that our result also requires similar constraint qualifications.
Such qualifications are not needed in \cite{garreis_inexact_2019},
but an additional assumption (eq. (3.10)) on the model accuracy
is made to show convergence to critical points.
In their case, the additional assumption is justified by domain-specific surrogates.

In next section, we will state the relevant optimality conditions for the multi-objective 
case. Afterwards, the main building blocks of our algorithm are described 
in~\cref{sec:trm_ideas}. The algorithm itself is given in~\cref{section:algorithm} and
convergence is shown in~\cref{sec:quasi_convergence,section:kkt_convergence}.
Finally, two numerical examples are discussed in~\cref{sec:examples} with a brief
discussion in~\cref{sec:conclusion}.

\subsection{Optimality Conditions}

We assume the reader to be familiar with the concept of Pareto-optimality
in the context of \ac{moo} (else see, for example,~\cite{miettinen_nonlinear_2013}).
In this subsection we introduce necessary conditions for a point
to be locally Pareto-optimal.
To state the criteria, we require the following assumption to
hold, so that all functions are sufficiently smooth:

\begin{assumption}%
  \label{ass:lipschitz_gradients_exact}
  The objective functions $\fL_{\ell},\, {\ell}=1,\ldots ,\fDim$, and
  the constraint functions $\gL_{\ell},\, 1, \ldots , \gDim$, and
  $\hL_{\ell},\,{\ell}=1,\ldots ,\hDim$,
  are twice continuously differentiable in an open domain
  containing $\feas$ and have Lipschitz continuous gradients on $\feas$.
\end{assumption}

There then is a formulation of Fermat's Theorem for optimization with multiple objectives:
\begin{mytheorem}[see {\cite{attouch}}]%
	\label{thm:fermats_theorem}
	Suppose \cref{ass:lipschitz_gradients_exact} holds.
	If $\vx^* \in \feas$ is locally Pareto-optimal for~\eqref{eqn:mop_nonlin},
	then there is no $\vd \in \tangcone{\feas}{\vx^*}$ such that
	for all ${\ell}\in \{1,\ldots , \fDim\}$ it holds that
	\(
    \langle \df_{\ell}(\vx^*), \vd \rangle < 0,
  \)
  where $\tangcone{\feas}{\vx^*}$ is the tangent cone of $\feas$ at $\vx^*$.
\end{mytheorem}

We call a point $\vx^* \in \feas$
satisfying the criterion in \cref{thm:fermats_theorem} \emph{Pareto-critical}.
\Cref{thm:fermats_theorem} can be used to motivate the following problem to compute
a descent direction and check for criticality:
\begin{equation}
	\min_{\vd \in \tangcone{\feas}{\vx}, \normTwo{\vd} \leq 1}
	\max_{{\ell} = 1, \ldots , \fDim}
	\langle
  \df_{\ell} (\vx),
  \vd
	\rangle.
  \label{eqn:DemandIdenticalOval}
\end{equation}
For a critical point $\vx\in\feas$ the optimal value is zero, else the minimizer is
a multi-descent direction (cf.~\cite[Th. 1.9.]{attouch},~\cite{fliege_svaiter}).
The choice of norm in the above problem is not really important, and so
we could use a linear norm and assume that at $\vx \in \feas$ certain constraint qualifications hold,
ensuring that the tangent cone equals the set of linearized directions
\[
L(\vx) =
\left\{
  \vd \in \mathbb{R}^\varDim :
  \tran{\vd} \dh_{\ell}(\vx) = 0, \, {\ell} = 1, \ldots , \hDim, \;
  \tran{\vd} \dg_{\ell}(\vx) \le 0, \, {\ell} \in A(\vx)
\right\},
\]
with $A(\vx) := \{{\ell} \in \{1,\ldots ,\gDim\}: g_{\ell}(\vx) = 0 \}$,
to obtain a linear problem, that also indicates criticality for~\eqref{eqn:mop_nonlin}.
Under constraint qualifications one can also derive KKT conditions which
then provide an equivalent definition of Pareto-criticality~\cite{hillermeier}.
In our algorithm, however, we do not use the set $L$.
%have slightly different direction
%sets to work with.
Instead, at the iterate $\vxk$ (which does not have to be feasible)
we use an approximation of the linearized feasible set,
\[
  \linFeas(\vxk) :=
  \left\{
    \vxk + \vd \in \mathbb{R}^\varDim:
    \vh(\vxk) + \DH(\vxk)\vd = \vZ, \;
    \vg(\vxk) + \DG(\vxk)\vd \le \vZ
  \right\},
\]
where $\DH(\vxk)$ and $\DG(\vxk)$ denote the \emph{full} Jacobian matrices of
the constraints $\vh$ and $\vg$, respectively.
The set $\linFeas(\vxk) - \vxk$ is not necessarily a cone,
but for feasible $\vxk$, it is a subset of $L(\vxk)$, and
intuitively it should not matter which of the sets is used near critical points.
Indeed, we have the following theorem:

\begin{mytheorem}%
  \label{thm:subproblems_implies_kkt}
  Suppose that $\vx$ is feasible and that a suitable constraint qualification holds.
  % and that the \ac{mfcq} hold at $\vx$,
  % i.e.,
  % the rows of $\DH(\vx)$ are linearly independent,
  % and there is a direction $\vd \in \mathbb{R}^\varDim$
  % such that $\DH(\vx)\vd = \vZ$ and
  % $\tran{\vd} \dg_{\ell}(\vx) < 0$ for all ${\ell} \in A(\vx)$.
  If the linear optimization problem
  \begin{equation}
    -
    \min_{
      \substack{
        \vd \in \linFeas(\vx) - \vx,\\
        \normInf{\vd} \leq 1
      }
    }
    \max_{{\ell} = 1, \ldots , \fDim}
    \langle
		\df_{\ell} (\vx),
		\vd
    \rangle
    \label{eqn:linear_opt_problem}
  \end{equation}
  has zero as its optimal value,
  then $\vx$ is also a KKT-point
  of~\eqref{eqn:mop_nonlin}
  (as defined in~\cite{pareto_tracer,hillermeier}).
\end{mytheorem}

\begin{proof}
	Dropping the argument $\vx$ for notational convenience,
  and denoting by $\DF$ the objective Jacobian,
	the linear problem \eqref{eqn:linear_opt_problem} is equivalent to
	\begin{equation}
			\max_{\ve d, {\beta}^-}
        \begin{bmatrix}
          \tran{\ve 0_\varDim} & 1
        \end{bmatrix}
        \begin{bmatrix}
          \ve d \\ {\beta}^-
        \end{bmatrix}
        \quad\text{s.t.}\quad
        \ve d \in \mathbb{R}^\varDim,
        \;
        {\beta}^- \in \mathbb R,
        \;
        \begin{bmatrix}
          -\eye{\varDim} & \ve 0_\varDim \\
          \eye{\varDim} & \ve 0_\varDim \\
          \DF & \ve 1_\fDim \\
          \DH & \ve 0_\hDim \\
          \DG & \ve 0_\gDim
          % \DmG & -\DmG & \ve 0_{\gDim} \\
        \end{bmatrix}%
			\begin{bmatrix}
					\ve d \\ {\beta}^-
			\end{bmatrix}
			\begin{matrix}
			\le \\ \le \\ \le \\ = \\ \le
			\end{matrix}
			\begin{bmatrix}
					\ve 1_\varDim\\
					\ve 1_\varDim\\
					\ve 0_\fDim \\
					\ve 0_\hDim\\
					-\vg \\
					%-\mgkSub
			\end{bmatrix}%
      .
		\label{eqn:legitimizersSannop}
    \tag{P}
	\end{equation}
	Consider also the dual problem:
	\begin{equation}
    \begin{aligned}
		&\min_{\ve y^1, \ldots , \ve y^5 }
			\begin{bmatrix}
				\tran{\ve 1_{\varDim}} & \tran{\ve 1}_{\varDim} &  \tran{\ve 0}_{\fDim} & \tran{\ve 0}_{\hDim} & - \tran{\vg}
			\end{bmatrix}
			\begin{bmatrix}
				\ve y^1 \\ \vdots \\ \ve y^5
			\end{bmatrix}
			\quad
      \text{s.t.}
      \quad
			\ve y^1 \ge \vZ_{\varDim},
			\ve y^2 \ge \vZ_{\varDim},
			\ve y^3 \ge \vZ_{\fDim},
			\ve y^4 \in \mathbb{R}^\hDim,
			\ve y^5 \ge \vZ_{\gDim}
      \quad
      \text{and}
      \\
			&\begin{bmatrix}
				-\eye{\varDim} & \eye{\varDim} & \tran{\DF} & \tran{\DH} & \tran{\DG}\\
				\tran{\ve 0}_{\varDim} & \tran{\ve 0}_{\varDim} & \tran{\ve 1}_{\fDim} & \tran{\ve 0}_{\hDim} & \tran{\ve 0}_{\gDim}
			\end{bmatrix}%
			\begin{bmatrix}
				\ve y^1 \\ \vdots \\ \ve y^5
			\end{bmatrix}
		  =
			\begin{bmatrix}
				\ve 0_{\varDim} \\ 1
			\end{bmatrix}
        .
    \end{aligned}
    \tag{D}
		\label{eqn:buntingGeologically}
	\end{equation}
	If $0$ is the optimal value of \eqref{eqn:legitimizersSannop},
  then ${\beta}^-=0$.
  By strong duality, the dual problem is feasible with optimal value 0,
	implying $\ve{y}^1 = \vZ$, $\ve{y}^2 = \vZ$ and $-\tran{\ve g} \ve y^5= 0$.
	The KKT equations immediately follow from the remaining constraints
  and from the complementary slackness property of dual solution pairs.
\end{proof}

% -----------------------sec01 ---------------------------------
\section{Trust-Region Concepts and Surrogates}\label{sec:trm_ideas}
As mentioned in the previous section,
we assume (at least some) of the objectives or constraints to be computationally expensive.
To approximate a Pareto-critical point whilst avoiding expensive function evaluations,
a trust-region approach is used:
The true functions $\vf,\vh$ and $\vg$ are modeled by $\mfk,\mhk$ and $\mgk$ respectively.
These models are constructed to be
sufficiently accurate within iteration-dependent trust-regions
\begin{equation*}
  \trk =
  B(\vxk; \Deltak)
  :=
  \left\{
    \vx \in \mathbb{R}^\varDim : \normTrIt{ \vx - \vxk } \le \Deltak
  \right\}
  .
\end{equation*}
Using the surrogate models, a step $\vsk$ and a step-size $\stepsize\KK$ are
determined -- in such a way as to reduce the constraint violation or
achieve an objective value reduction at the
trial point $\vxk + \stepsize\KK \vsk$,
which is tested as a candidate for the next iterate.
The step $\stepsize\KK\vsk$ is computed so that the trial point is contained
in both the trust-region and the approximate linearized feasible set at $\vxk$:
\begin{equation}
  \linFeasApproxk =
  \linFeasApproxk(\vxk) =
  \left\{
      \vxk + \vs \in \mathbb{R}^\varDim:
      %\begin{aligned}
              \mhk(\vxk) + \DmH(\vxk) \cdot \vs = \vZ, \;
              \mgk(\vxk) + \DmG(\vxk)\cdot \vs \le \vZ
      %\end{aligned}
  \right\}
  ,
  \label{eqn:orgiastThripses}
\end{equation}
where $\DmH$ and $\DmG$ are now the full model Jacobians.
This introduces uncertainty and iterates might no longer be feasible.
Hence, we treat the constraints as \emph{relaxable}, and we might have to evaluate $\vf$
(and $\vg$ and $\vh$) outside of $\feas$, which motivates the next assumption.
It states that all functions have to be available in all possible trust-regions.
Note, that the algorithm actually can use two trust-region sizes,
the preliminary size $\DeltaPrek$ at the beginning of an iteration 
and $\Deltak \leq \DeltaPrek$ after the \emph{Criticality Routine},
which is described in \cref{section:crit_routine}.

\begin{assumption}%
  \label{ass:covering}
  There is a constant $0<\DeltaMax<\infty$ such that
  for every $\kk \in \mathbb{N}_0$ the trust region sizes conform to
  $0 < \Deltak \le \DeltaPrek < \DeltaMax$.
  All functions (the true functions and their models)
  are defined in these regions, i.e.,
  on the set
  $\covering = \bigcup_{\kk\in \mathbb{N}_0}
    % B\KK.
    \ball{\vxk}{\DeltaPrek}
  $
\end{assumption}

% In the convexly constrained case, convergence
% %of a sub-sequence $\{\kk_{\ell}\}$
% to a critical point is proven via
% \[
% \lim_{\kk \to \infty} {\omega}\left(  \vxk; \vf  \right)  = 0,
% \]
% where ${\omega}$ is the optimal value for the same problem \eqref{eqn:sub_problem},
% except that the true objective functions $\vf$ (and their derivatives) are used,
% which makes it a criticality measure for \eqref{eqn:mop_nonlin}
% (cf. \cite{ours,conn_trust_region_methods}).

\subsection{Fully Linear Models}

In this subsection, we first want to explain what it means for the model functions
to be sufficiently accurate.
Although the true derivative information is not used in any computation,
we will need it for the convergence
analysis and hence require the following generalization of~\cref{ass:lipschitz_gradients_exact}:
\begin{assumption}%
    \label{ass:lipschitz_gradients_true}
    The objective functions $\fL_{\ell},\, {\ell}=1,\ldots ,\fDim$,
    the constraint functions $\gL_{\ell},\, 1, \ldots , \gDim$,
    and $\hL_{\ell},\,{\ell}=1,\ldots ,\hDim$,
    are twice continuously differentiable in an open domain containing $\covering$
    and have Lipschitz continuous gradients on $\covering$.
\end{assumption}

With \cref{ass:lipschitz_gradients_true} we can define models satisfying
error bounds w.r.t. to their construction radius:

\begin{mydefinition}[Fully Linear Models]%
  \label{def:fully_linear_models}
  Let $\DeltaMax > 0$ be a given constant and let $\fL\colon \mathbb{R}^n\to \mathbb{R}$ be a
  scalar-valued function
  satisfying \cref{ass:lipschitz_gradients_true}.
  A set of model functions $\mathcal M = \{ \mfL \colon \mathbb{R}^n \to \mathbb{R} \}$ is called a
  fully linear class of models for $f$ if the following hold:
  \begin{enumerate}
      \item\label{def:fully_linear_models1}
        There are positive constants $\CCfl,\CCflDiff$ and $\LipschitzBoundModels$ such
        that for any $\Delta\in (0,\DeltaMax]$
        and for any $\vx \in \covering$ there is a model function $\mfL \in \mathcal M$,
        with Lipschitz continuous gradient and corresponding Lipschitz constant
        bounded by $\LipschitzBoundModels$, such that
        the error between the gradients and the error between the values satisfy
        \begin{equation*}
          \normTwo{\df(\ve {\xi}) - \dmf(\ve {\xi})}
          \le \CCflDiff \Delta
          \quad\text{and}\quad
          \left|
            \fL(\ve {\xi}) - \mfL(\ve {\xi})
          \right|
          \le \CCfl \Delta^2, \qquad\forall\ve {\xi}
          \in \ball{\vx}{\Delta} \cap \covering.
        \end{equation*}
      \item\label{def:fully_linear_models2}
        For this class $\mathcal M$ there exists a
        ``model-improvement'' algorithm that, in a finite, uniformly bounded number of steps can either establish that a given model $\mfL\in \mathcal M$ is fully linear on $\ball{\vx}{\Delta}$
        or find a model $\mfL \in \mathcal M$ that is fully linear on $\ball{\vx}{\Delta}$.
  \end{enumerate}
\end{mydefinition}

\begin{mydefinition}%
  \label{def:fully_linear_vector_models}
  Let $\DeltaMax>0$ be a given constant and let $\vf = \tran{[\fL_1, \ldots , \fL_\fDim]}$
  be a vector of functions satisfying the requirements of
  \cref{def:fully_linear_models} with classes $\mathcal M_{\ell}$ and
  constants $(\CCfl_{\ell}, \CCflDiff_{\ell}, \LipschitzBoundModels_{\ell}), {\ell}=1,\ldots ,\fDim$.
  Then
  \[
    \mathcal M =
    \{ \mf = \tran{[\mfL_1, \ldots , \mfL_\fDim]} :
    \mfL_1 \in \mathcal M_1, \
    \ldots , \
    \mfL_\fDim \in \mathcal M_\fDim \}
  \]
  is a class of fully linear, vector-valued model functions with constants
  $\max_{\ell} \CCfl_{\ell} > 0$,
  $\max_{\ell} \CCflDiff_{\ell}>0$
  and
  $\max_{\ell}\LipschitzBoundModels_{\ell}>0$.
  A vector of functions $\mf \in \mathcal M$ is deemed fully linear, if all components $\mfL_{\ell}$ are fully linear.
  The improvement algorithms for $\mathcal M_{\ell}$ are applied component-wise.
\end{mydefinition}
A function that is fully linear for $\Deltak > 0$ is automatically
fully linear for any smaller radius $0<\Delta<\Deltak$.
When the trust-region radius is bounded above,
then the constants $\CCfl>0$ and $\CCflDiff > 0$ can be chosen large
enough such that a fully linear model stays fully linear in enlarged trust-regions:

\begin{mylemma}[Lemma 10.25 in~\cite{conn_introduction_to_df_optim}]%
  \label{thm:fully_linear_enlarged_trust_region}
  For $\vx \in \covering$ and $\Delta \in (0,\DeltaMax]$ consider a function $\fL$
  and a fully linear model $\mfL$ of $\fL$ with constants
  $\CCfl, \CCflDiff, \LipschitzBoundModels > 0$.
  Let $L^{\fL} > 0$ be a Lipschitz constant of $\df$.
  Assume \ac{wlog} that
  \begin{equation}
    \LipschitzBoundModels + L^{\fL} \le \CCfl
    \text{ and }
    \nicefrac{\CCflDiff}{2} \le \CCfl.
    \label{eqn:large_fully_linear_constants}
  \end{equation}
  Then $\mfL$ is fully linear on $\ball{\vx}{\tilde{\Delta}}$ for any
  $\tilde{\Delta} \in [\Delta, \DeltaMax]$ with the same constants.
\end{mylemma}

\begin{assumption}%
  \label{ass:fully_linear_models}
  For any $\kk \in \mathbb{N}_0$ the models $\mfk$ are fully linear on
  $\ball{\vxk}{\Deltak}$ as in \cref{def:fully_linear_vector_models}
  with constants $\CCflF, \CCflDiffF$ and $\LipschitzBoundModel{\vf}$
  that are chosen large enough such that~\eqref{eqn:large_fully_linear_constants} is fulfilled globally.
  The same holds for the models $\mgk$ of $\vg$ with constants $\CCflG, \CCflDiffG, \LipschitzBoundModel{\vg}$
  and the models $\mhk$ of $\vh$ with constants $\CCflH, \CCflDiffH$ and $\LipschitzBoundModel{\vh}$.
  We also assume that all models are \textbf{interpolating} at $\vxk$.
\end{assumption}

\subsection{Composite Step Approach and Sub-Problems}
Just like in previous articles on multi-objective trust-region
algorithms~\cite{thomann_paper,ours,villacorta,qu_goh_liang,ryu},
we use the maximum-scalarization
\[
  \scalShort(\vx) = \scal{\vf}(\vx) := \max_{{\ell}=1,\ldots ,\fDim} \fL_{\ell}(\vx)
  \quad\text{and}\quad
  \scalk(\vx) = \scal{\mfk}(\vx) := \max_{{\ell}=1,\ldots ,\fDim} \mfL_\ell(\vx)
\]
to determine objective reduction.
The idea then is to find a step $\vsk$ that approximately solves
\begin{equation*}
  %\begin{aligned}
  %&
  \min_{\vs \in \mathbb{R}^\varDim}
  \scal{\mfk}(\vxk + \vsk)
%  &
    \quad\text{s.t.}\quad
    \vs \in \left( \linFeasApproxk - \vxk \right),
%        \\
%            &\mhk( \vxk ) + \DmH( \vxk ) \vs = \vZ,
%        \\
%            &\mgk( \vxk ) + \DmG( \vxk ) \vs \le \vZ,
%        \\
%        &
    \,
    \normTrIt{\vs} \le \Deltak
    ,
%    \end{aligned}
\end{equation*}
with the inexact linearized feasible set defined in~\eqref{eqn:orgiastThripses}.
Without constraints, inexact line-search can be used.
In our case, we use composite-step approach and
$\vsk$ is split into a normal component $\vnk$
towards feasibility and a descent direction $\vdk$
(see~\cite{fletcher_global_2002,walther2016} for details).
Then, the normal step can be computed with
\begin{equation}
%    \begin{aligned}
        %&\min \normTrIt{ \vn } &\text{s.t.}
%        &
  \min \normTwo{ \vn}^2
  \quad\text{s.t.}\quad
  \vn \in \left( \linFeasApproxk - \vxk \right).
%        \\
%        &\mhk(\vxk) + \DmH(\vxk) \cdot \vn = \vZ,
%        \\
%        &\mgk(\vxk) + \DmG(\vxk)\cdot \vn \le \vZ,
        %&\normTrIt{\vn} \le \DeltaPrek,
%    \end{aligned}
    \label{eqn:itr_normal}
    \tag{ITRN\ensuremath{\KK}}
\end{equation}
%where $\mhk$ and $\mgk$ are the surrogate models and $\DmH, \DmG$ their jacobians.
If a normal step $\vnk$ has been found
-- and if $\normTrIt{\vn} \le \DeltaPrek$ --
the descent direction $\vdk$ can be taken as the minimizer of
% \begin{equation}
%   \begin{aligned}
%     {\omega}\left(  \vxNk; \mfk, \linFeasApproxk, \normIt{{\bullet}} \right)
%     =-
%     &
%       \min_{\ve d \in \mathbb{R}^\varDim, {\beta}\in \mathbb{R}} {\beta} &\text{s.t.}
%     \\
%     &\DmF(\vxNk) \cdot \vd \le {\beta},
%     \\
%     &\vd \in \left(\linFeasApproxk - \vxNk\right)
%     \\
%     &\normIt{\vd} \le 1,
%   \end{aligned}
%   \tag{ITRT\ensuremath{\KK}}
%   \label{eqn:itr_tangential}
% \end{equation}
\begin{equation}
  {\omega}\left(  \vxNk; \mfk, \linFeasApproxk, \normIt{{\bullet}} \right)
  =-
  \min_{
    \beta \in \mathbb{R},
    \vd \in \mathbb{R}^\varDim
  }
  \beta
  \quad\text{s.t.}\quad
    \DmF(\vxNk) \cdot \vd \le {\beta}, \
    \vd \in \left(\linFeasApproxk - \vxNk\right), \
    \normIt{\vd} \le 1,
  \tag{ITRT\ensuremath{\KK}}
  \label{eqn:itr_tangential}
\end{equation}
where $\linFeasApproxk$ is the approximate linearized feasible set
according to~\eqref{eqn:orgiastThripses} and
$\vxNk = \vxk + \vnk$.\\
We actually only compute the descent direction if there is enough movement possible after
performing the normal step, i.e., the normal step must not be too large.
We call such a step \emph{compatible}, and it is defined with respect to the
preliminary radius
$\DeltaPrek$
%(which defines the region for which the descent step is sought):
%(before the so called Criticality routine)
:
\begin{mydefinition}%
    \label{def:compatible_normal_step}
    Let $\CCDelta \in (0,1], \CCmu > 0$ and $\mu\in (0,1)$
    be constants.
    The minimizer $\vnk$ of~\eqref{eqn:itr_normal} is called compatible if
    \begin{equation}
        \normTrIt{\vnk}
        \le
        \CCDelta \DeltaPrek
            \min \left\{ 1, \CCmu \DeltaPrek^{\mu} \right\}
            .
        \label{eqn:compatible_normal_step}
    \end{equation}
\end{mydefinition}
Furthermore, we assume a normal step to exist if the true constraint violation is not too large, as
measured by the infeasibility function
\[
  {\theta}\left(  \vx  \right)
  =:
  \max\left\{
      \max_{{\ell}=1, \ldots , \hDim} | \hL_{{\ell}}(\ve x ) |, \max_{{\ell}=1,\ldots ,\gDim} \gL_{\ell}(\vx)
  \right\}.
\]
%In that case, the normal step is also bounded:
\begin{assumption}[Existence and Boundedness of Normal Step]%
  \label{ass:normal_step_exists_bounded}
  If $\thetak = {\theta}\left( \vxk \right) \le \CCdeltan$, for a constant $\CCdeltan > 0$,
  then $\vnk$ exists and there is a constant $\CCubn > 0$ such that
  \begin{equation}
      \normTrIt{\vnk} \le \CCubn \thetak.
      \label{eqn:normal_step_exists_bounded}
  \end{equation}
\end{assumption}
\Cref{ass:normal_step_exists_bounded} is a standard assumption and the reasoning
behind it can be found in~\cite{fletcher_global_2002}.
Finally, the step-size $\stepsize\KK$ is determined in such a way that
$\normTrIt{\vxk + \vnk + \stepsize\KK \vdk} \le \Deltak$
and a \emph{sufficient decrease condition} (for the objective surrogates) is satisfied,
which is described and justified in \cref{section:sufficient_decrease}.
The term $\stepsize\KK \vdk$ is also called \emph{tangential step} and
altogether the step $\vsk$ results in 
the \emph{trial-point} $\vxSk = \vxk + \vsk =  \vxNk + \stepsizek \vdk$.

\subsection{Filter Mechanism}
As mentioned already, iterates can become infeasible.
We employ a so-called \emph{Filter} to drive them back towards 
the feasible set.
Thus, the trial point $\vxSk = \vxk + \vsk$
is tested not only
for actual decrease of the original functions 
but also against previous
iterates stored in the filter $\filter$.
%The handling of $\vxSk$ is determined using a so called filter $\filter$.
$\filter$ is a set of tuples $({\theta}_j, \Phi_j)$ describing a forbidden area
in image space.
In fact, the tuples in $\filter$ (w.r.t. $\vf$) are currently non-dominated for
the bi-objective optimization problem of minimizing both $\theta(\vx)$ and $\scalf{\vx}$.
The trial point $\vxSk$ is only acceptable for $\filter$ if its value
tuple is also (sufficiently) non-dominated.
An acceptable trial point is further tested and if it sufficiently reduces the \emph{true}
objectives it is kept as the next iterate.
The current iterate might be added to the filter if the predicted objective decrease
is small compared to the constraint violation, or if the latter is too large.
More formally, we use the following definition:

\begin{mydefinition}[Multi-Objective Filter]%
  \label{def:filter}
  A filter $\filter$ with respect to some function 
  $\vf\colon \mathbb{R}^\varDim \to \mathbb{R}^\fDim$
  is a discrete set of tuples $\{({\theta}_j, \Phi_j)\}\subset \mathbb{R}^2$,
  and a point $\vx$ is \emph{acceptable} for the filter iff
  \begin{equation*}
      {\theta}\left( \vx  \right) \le (1-\filterConst) {\theta}_j
      \quad\text{or}\quad
      \scalf{\vx} \le \Phi_j - \filterConst  {\theta}_j
      \quad \forall ({\theta}_j, \Phi_j)\in \filter.
  \end{equation*}
  When a point $\vx$ is \emph{added} to the filter the tuple
  $({\theta}\left( \vx \right), \scalf{\vx})$ is added to the set $\filter$, but all
  tuples $({\theta}_j, \Phi_j)$ with
  \[
    {\theta}_j \ge {\theta}\left( \vx \right)
    \text{  and  }
    \Phi_j - \filterConst {\theta}_j
    \ge
    \scalf{\vx} - \filterConst {\theta}\left( \vx \right)
  \]
  are removed from the set.
\end{mydefinition}

As can be seen from~\cref{def:filter} a filter strengthens non-dominance testing
by employing a positive offset $\filterConst {\theta}_j$.
Note, that we could also use a stricter,
$(\fDim + 1)$-dimensional filter instead of the $2$-dimensional filter,
similar to~\cite{gould2004multidimensional}, by using $\vf({\bullet})$ instead of $\scalf{{\bullet}}$.

% Similar to the single objective case, a point $\vxOpt$ is critical for \eqref{eqn:mop_nonlin} iff
% it is feasible, a constraint qualification holds, and if ${\omega}_2 \left( \vxOpt; \vf, \linFeas( \vxOpt ) \right) = 0$,
% where $\linFeas$ is the exact linearized feasible set and the 2-norm is used.
% \todo{Explain "criticality measure", provide optimality conditions, projection context etc.}
% We show that our algorithm produces iterates with $\vxNk \to \bar{\vx}$ and ${\theta}(\bar{\vx}) = 0$ as well as
% ${\omega}_2 \left(\vxNk; \mfk, \linFeasApproxk \right) \to 0$.
% In accordance with \cite{echebest_inexact_2017} we could call such limit points \emph{quasi-stationary}.
% From this, it is easy to see that it must also hold that
% ${\omega}_2\left(\vxNk; \vf, \linFeasApproxk \right) \to 0$,
% where the true \emph{objective} derivatives $\DF$ are used instead of the surrogates $\DmF$.
% Finally, it remains to show, that it must then hold (under the MF constraint qualifications) that
% ${\omega}_2\left( \bar{\vx}; \vf, \linFeas \right) = 0$ as well, which means that $\bar{\vx}$ is truly critical and a KKT-point.
% So, in contrast to \cite[Lemma 3.2]{fletcher_global_2002} we cannot directly use the
% inexact measures to show criticality, but have to take further considerations in
% \cref{section:kkt_convergence}.

\subsection{Sufficient Decrease}%
\label{section:sufficient_decrease}

In this sub-section we want to explain what is meant by ``sufficient decrease''.
In short, we have to relate the model predicted objective reduction to the
criticality value.
The criticality is the optimal value of \eqref{eqn:itr_tangential},
and it is thus dependent on the norm.
As indicated above, we allow for the usage of iteration
dependent problem norms $\normIt{{\bullet}}$ and trust-region norms $\normTrIt{{\bullet}}$.
This way, the algorithm can be implemented by using norms that are suited best
for the problem geometry or the inner solver(s).
We only require that these norms be uniformly equivalent to the Euclidean norm:

\begin{assumption}%
  \label{ass:equivalent_norms}
  There is a constant $\CCnorm \ge 1$
  such that for all $\vx\in \covering$ and all $\kk\in \mathbb{N}_0$ and $\norm{{\bullet}}_* = \normTrIt{{\bullet}}$ or $\norm{{\bullet}}_* = \normIt{{\bullet}}$ it holds that
  \begin{equation}
      \frac{1}{\CCnorm}
      \norm{\vx}_*
      \le
      \normTwo{\vx}
      \le
      \CCnorm \norm{\vx}_*
      .
      \label{eqn:equivalent_norms}
  \end{equation}
\end{assumption}

\begin{myremark}%
  \label{remark:equivalent_norms}
  Any two norms that are uniformly equivalent to $\normTwo{{\bullet}}$ with constant $\CCnorm$ are pairwise equivalent with constants $\nicefrac{1}{\CCnorm^2}$ and $\CCnorm$.
\end{myremark}

\begin{mylemma}%
  \label{thm:critical_values_equivalent}
  Suppose \cref{%
  ass:covering,%
  ass:lipschitz_gradients_true,%
  ass:equivalent_norms,%
  ass:fully_linear_models%
  } hold.
  Denote by
  $
  \omegak := {\omega}( \vxNk; \mfk, \linFeasApproxk, \normIt{{\bullet}} )
  $
  the optimal value of~\eqref{eqn:itr_tangential} and 
  by $\omegak_2$ the optimal value if the 2-norm is used,
  $\omegak_2 =  {\omega}( \vxNk; \mfk, \linFeasApproxk, \normTwo{{\bullet}} )$.
  There is a constant $\CCcrit \ge 1$ such that for any 
  $\kk \in \mathbb{N}_0$
  for which the normal step exists,
  it holds that 
  \begin{equation}
      \frac{1}{\CCcrit} ~
      \omegak
      \le
      \omegak_2
      \le
      \CCcrit ~ \omegak
      .
      \label{eqn:critical_values_equivalent}
  \end{equation}
\end{mylemma}
The proof of \cref{thm:critical_values_equivalent} can be found in the appendix.
The lemma shows that we can relate the iteration dependent
inexact critical values $\omegak$ and $\omegak_2$.
Furthermore, we show in the appendix that it is then sensible to assume
a sufficient decrease condition for the objective
surrogate functions as per~\cref{ass:sufficient_decrease},
as long as an additional (standard) assumption holds:
\begin{assumption}%
    \label{ass:hess_norm_bounded}
    The norm of all model Hessians of the objective function surrogates is
    uniformly bounded above,
    i.e., there is a positive constant
    $\CChessF > 0$ such that for all $\kk \in \mathbb{N}_0$
    \[
      \norm{ \hess\mfL\KK_{\ell} (\ve {\xi}) }_2 \le \CChessF
      \text{ for all ${\ell}=1,\ldots ,\fDim$, %
        and all $\ve {\xi}\in \linFeasApprox\KK(\vxk) \cap \trk$.%
      }
    \]
\end{assumption}
We know how to find suitable model functions that
satisfy~\cref{ass:hess_norm_bounded},
including Lagrange interpolation polynomials or \acp{rbf},
so that finally we can
state the sufficient decrease assumption:
\begin{assumption}[Sufficient Decrease]%
  \label{ass:sufficient_decrease}
  Suppose
  \cref{%
  ass:covering,%
  ass:lipschitz_gradients_true,%
  ass:fully_linear_models,%
  %ass:normal_step_exists_bounded,%
  %ass:equivalent_norms,%
  %ass:hess_norm_bounded,%
  }
  hold and that
  $\Deltak \in (0,\DeltaMax)$.
  Let $\scalk := \scal{\mfk}$.
  If $\vn\KK$ is compatible, and $\vdk$ is a minimizer of~\eqref{eqn:itr_tangential} at $\vxk$ for
  $\normIt{{\bullet}}$ and $\omegak_2$ is the optimal value for~\eqref{eqn:itr_tangential} if the $2$-norm is used,
  then there is a step-length $\stepsizek \ge 0$ such that
  $\vxk + \vnk + \stepsizek \vdk \in \linFeasApproxk \cap \trk$ and
  \begin{equation}
      \scalk \left( \vxNk \right)
      -
      \scalk \left( \vxNk + \stepsizek \vdk \right)
      \ge
      \CCsd \omegak_2
      \min \left\{
      \frac{\omegak_2}{\CCcritDenom},
      \Deltak,
      1
      \right\},
      \label{eqn:sufficient_decrease}
  \end{equation}
  for constants $\CCsd \in (0,1)$ and $\CCcritDenom \ge 1$.
\end{assumption}

Throughout the rest of this article we use a slightly modified measure
for notational convenience:
%to have a the sufficient decrease condition in the usual form
%(similar to unconstrained optimization):
\begin{mycorollary}[Modified Criticality Measure]%
  \label{thm:ScaryPitchOld}
    For any $\kk\in \mathbb{N}_0$ and $\omegak$ as in~\cref{thm:critical_values_equivalent},
    define the criticality measure $\chik := \min \left\{ 1, \omegak \right\}$
    and denote by $\chik_2$ the corresponding value, if the 2-norm is used
    % for the computation of $\omegak_2$
    instead of $\normIt{{\bullet}}$.
    Then $\lim_{k\to \infty} \omegak_2 = 0$ if and only if
    $\lim_{k\to \infty} \chik_2 = 0$,
    and if \cref{ass:sufficient_decrease} holds,
    then it also follows (with $\CCcritDenom \ge 1$)
    that
    \begin{equation*}
        \scalk \left( \vxNk \right)
        -
        \scalk\left( \vxNk + \stepsizek \vdk \right)
        \ge
        \CCsd \chik_2
        \min \left\{
        \frac{\chik_2}{\CCcritDenom},
        \Deltak
        \right\}.
        \label{eqn:sufficient_decrease_critMin}
  \end{equation*}
\end{mycorollary}
%
%-----------------------sec02 ---------------------------------
\section{Discussion of the Algorithm}%
\label{section:algorithm}

The behavior of the algorithm is determined by several additional algorithmic parameters:
\begin{center}
    \begin{tabularx}{\linewidth}{cX}%
    \toprule
         Parameter(s) & Description \\
         \midrule
         $0 < \DeltaInit \le \DeltaMax < \infty$ & initial and maximum trust-region radius
         \\
         $0<\gammass \le \gammas < 1 \le \gammag$ & shrinking and growing parameters for the trust-regions
         \\
         $0 < \CCaccept \le \CCsuccess < 1$ & acceptance thresholds for trial point test
         \\
         $0 < \CCcritTestchi < 1$, $0 \le \CCcritTesttheta \le \CCdeltan$ & thresholds for criticality test
         \\
         ${\kappa}_{\theta} \in (0,1), \psi > \frac{1}{1+{\mu}}$ &
         threshold parameters in~\eqref{eqn:algo_model_dec}
         \\
         $0 < \CCBeta < \CCMu$ and $\CCalpha \in (0,1)$ &
         Criticality Routine threshold factors and backtracking constant
         \\
         $\CCsd \in (0,1)$ & \emph{s}ufficient \emph{decrease} constant in \cref{ass:sufficient_decrease}
         \\
         $\CCDelta \in (0,1), \CCmu > 0, \mu\in (0,1)$ & constants defining
                                                         compatibility in \cref{def:compatible_normal_step}
         \\
         $\CCdeltan > 0, \, \CCubn > 0$ & existence of normal step in \cref{ass:normal_step_exists_bounded}
    \\
    \bottomrule
    \end{tabularx}
\end{center}

Below, the \emph{Criticality Test} is used instead of the original stopping criterion
``${\chi}(\vxk) = 0$'' because the surrogate models are inexact.
If an iterate is nearly feasible and nearly critical for the surrogate problem, then
the criticality routine is entered and the trust-region radius is reduced
to make the models more precise.
At a truly critical point the routine loops infinitely, but if a point is not
critical for the true functions,
we exit and continue with regular iterations.
For further details of the Criticality Test and the Criticality Routine we refer to~\cite{conn_introduction_to_df_optim}.

\subsection*{Algorithm}
\begin{enumerate}
    \setcounter{enumi}{-1}
    \item\label[step]{algo:init} \textbf{Initialization:}
    Let $\kk \leftarrow 0$, $\filter \leftarrow \emptyset$ and $\vnk = \None$.
    Evaluate $\vf( \vxk), \vg(\vxk), \vh(\vxk)$ and compute $\thetak = {\theta}\left( \vx \right)$.
    Build surrogate models $\mfk,\mgk,\mhk$ that are
    fully linear in $\trk$ with radius $\DeltaPrek := \Delta\itDeltaCmd{0}$.
    \item\label[step]{algo:compat_test} \textbf{Compatibility Test:}
    \begin{itemize}
        \item If $\vnk \ne \None$ and $\vnk$ is compatible w.r.t. $\DeltaPrek$,
        go to step~\ref{algo:descent}.
        \item If $\vnk = \None$, try to compute $\vnk$.
        If $\vnk$ exists and is compatible, go to step \ref{algo:descent}.
    \end{itemize}
    \item\label[step]{algo:restoration} \textbf{Restoration}
    Add $\vxk$ to the filter, set
    $\Deltak \leftarrow \DeltaPrek, \chik \leftarrow \chikPre$,
    and attempt to find a \emph{restoration step}
    $\ve r\KK$ and $\DeltaPreNext>0$ for which
    ~\eqref{eqn:itr_normal} is compatible at
    $(\vxk + \ve r\KK, \DeltaPre\itDeltaCmd{\kk+1})$ and for which
    $\vxk + \ve r\KK$ is acceptable for $\filter$.
    Set $\vxNext\leftarrow \vxk + \ve r\KK$, keep $\DeltaPreNext$ and go to step~\ref{algo:updates}.
    \item\label[step]{algo:descent} \textbf{Descent Step:}
    Compute a descent direction $\vdk$ and $\chikPre$ with~\eqref{eqn:itr_tangential}.
    \item\label[step]{algo:crit_test} \textbf{Criticality Test:}
    If $\thetak < \CCcritTesttheta$
    \textbf{and}
    $\bigl(${}$\chikPre < \CCcritTestchi$ and
    $\DeltaPrek > \CCMu \chikPre${}$\bigr)$,
    then enter the \textbf{Criticality Routine} to get $\Deltak$ and $\chik$ and
    update $\vdk$ and $\vnk$.
    Else, set $\Deltak \leftarrow \DeltaPrek$ and
    $\chik \leftarrow \chikPre$.
    \item\label[step]{algo:acceptance_test} \textbf{Acceptance Test:}
    Compute a step-size $\stepsizek > 0$ such that \cref{ass:sufficient_decrease} is fulfilled.
    Set $\vxSk = \vxk + \vnk + \stepsizek \vdk$.
    Compute $\vf(\vxSk)$ and ${\theta}\left( \vxSk \right)$.
    \begin{itemize}
        \item If $\vxSk$ is not acceptable for the augmented filter
        $\filter \cup \{ (\thetak, \scalf{\vxk}) \}$ OR
        \item If
        \begin{align}
            \scal{\mfk}(\vxk) -
            \scal{\mfk}(\vxSk)
            &\ge {\kappa}_{\theta} \thetak^\psi
            \label{eqn:algo_model_dec}
            \shortintertext{AND}
            {\rho}\KK :=
            \frac{%
            \scal{\vf}(\vxk) -
            \scal{\vf}(\vxSk)}{%
            \scal{\mfk}(\vxk) -
            \scal{\mfk}(\vxSk)}
            &< \CCaccept,
            \label{eqn:algo_ratio_test}
        \end{align}
    \end{itemize}
    keep $\vxNext \leftarrow \vxk, \vn^{(\kk+1)} \leftarrow\vnk$,
    choose $\DeltaPreNext\in [\gammass \Deltak, \gammas \Deltak]$, increment $\kk$ and go to step~\ref{algo:updates}.
    \item\label[step]{algo:filter_test} \textbf{Filter Test:}
    If~\eqref{eqn:algo_model_dec} fails, include $\vxk$ in the filter (${\theta}$ -iteration).
    \item\label[step]{algo:iterate_update} \textbf{Iterate Updates:}
    Set $\vxNext \leftarrow \vxSk$ and choose
    \[
      \DeltaPreNext \in
      \begin{cases}
      [\gammass \Deltak, \gammas \Deltak] & \text{if ${\rho}\KK < \CCsuccess$}, \\
      [\Deltak, \max\{\gammag \Deltak, \DeltaMax\}] & \text{if ${\rho}\KK\ge \CCsuccess$.}
      \end{cases}
    \]
    \item\label[step]{algo:updates} \textbf{Model Updates:}
      Update the surrogates for $\vxNext$ and $\DeltaPreNext$.
      If the models are dependent on the trust-region radius, then set $\vn^{(\kk+1)}
      \leftarrow\None$.
      Finally, let $\kk \leftarrow \kk + 1$ and go to step~\ref{algo:compat_test}.
\end{enumerate}

Note, that in contrast to regular trust-region methods, a step that does not pass the test~\eqref{eqn:algo_ratio_test} may still be accepted as the next iterate, namely in a ${\theta}$-iteration.
The restoration step $\vx\itCmd{\kk + 1} = \vxk + \ve r\KK$ aims at reducing the
constraint violation and is usually calculated by solving
\begin{equation}
    \min_{\vx } {\theta}\left( \vx \right).
    \tag{R}
    \label{eqn:restoration_problem}
\end{equation}
For a detailed discussion of the role of~\eqref{eqn:algo_model_dec} and the restoration procedure,
we refer to~\cite{fletcher_global_2002}.
We only want to note a few important properties:
\begin{itemize}
    \item No feasible iterate is ever added to the filter.
    \item Unless noted otherwise, we assume that the restoration procedure is always able to find
    a suitable restoration step.
    \item There are no two successive restoration iterations.
    \item Under \cref{ass:normal_step_exists_bounded}, there can only ever be a finite number of
    sub-iterations in the restoration procedure until a compatible normal step exists.
\end{itemize}
\subsection{Criticality Routine}%
\label{section:crit_routine}

The criticality routine is provided with the current models $\mfk,\mgk,\mhk$
and the preliminary trust-region radius $\DeltaPrek$.
All the models as well as the normal step $\vnk$
and the inexact criticality value may be modified within this sub-routine.\\
When the Criticality Routine starts, we know that
$\DeltaPrek > \CCMu \chik$ and that $\vnk$ is compatible for $\DeltaPrek$.
\begin{enumerate}
    \item\label{critAlgo:init}
    Let $j\leftarrow 0$
    and set
    $\delta_j \leftarrow \DeltaPrek$,
    ${\chi}_j \leftarrow \chikPre$ and $\vn_j \leftarrow \vnk$.
    \item While
    %$\delta_j^{\frac{1}{1+{\mu}}} > \left(\CCDelta \CCmu\right)^{\frac{1}{1+{\mu}}} \CCMu {\chi}_j$
    ${\delta}_j > \CCMu {\chi}_j$:
    \begin{enumerate}
        \item\label{critAlgo:radius}
        Tentatively decrease the radius:
        ${\delta}_j^+ \leftarrow \CCalpha \delta_{j} = \CCalpha^j \DeltaPrek$.
        \item\label{critAlgo:model_updates}
        Make the models fully linear for the new radius ${\delta}_j^+$.
        \item\label{critAlgo:compat_test}
        Solve~\eqref{eqn:itr_normal} to get $\vn_j^+$ for the new models and ${\delta}_j^+$.\\
        If $\vn_j^+$ is not compatible for ${\delta}_j^+$, then \texttt{BREAK}.
        \item\label{critAlgo:index}
        Set $j\leftarrow j+1$, ${\delta}_j\leftarrow {\delta}_{j-1}^+$ and keep
        $\vn_j \leftarrow \vn_{j-1}^+$.
        \item\label{critAlgo:crit_update}
        Solve~\eqref{eqn:itr_tangential}  to get ${\chi}_j$ for the new models and ${\delta}_j$.
    \end{enumerate}
    \item\label{critAlgo:finish}
    Keep the updated models as $\mfk, \mgk, \mhk$,
    set $\chik \leftarrow {\chi}_j$, $\vnk \leftarrow \vn_j$
    and choose
    \begin{equation}
        \Deltak \leftarrow \min\left\{ \max\{ {\delta}_j, \CCBeta \chik\}, \DeltaPrek \right\}.
    \label{eqn:crit_test_final_radius}
    \end{equation}
\end{enumerate}

There are two ways the routine can stop.
If the radius is sufficiently small compared to the criticality, i.e.,
${\delta}_j \le \CCMu {\chi}_j$,
or if the next prospective normal step $\vn_j^+$ is no longer compatible for
the smaller radius ${\delta}_j^+$.
The second criterion ensures (inductively), that after the routine has stopped
finitely the normal step will be compatible for ${\delta}_j$ and for any radius $\Delta\ge {\delta}_j$.
Because it always holds that ${\delta}_j\le \DeltaPrek$ and $\max\{ {\delta}_j, \CCBeta \chik\} \ge {\delta}_j$,
it follows from~\eqref{eqn:crit_test_final_radius} that $\vnk$ is compatible for
$\Deltak$.

One reason for the Criticality Test (and the Criticality Routine)
is to have a lower bound on the inexact criticality
if the constraint violation is sufficiently small, i.e.,
$\thetak \le \CCcritTesttheta \le \CCdeltan$ and the routine stops because
of the first condition.
From $\CCBeta \le \CCMu$ as well as $\bar{{\delta}_j} \le \DeltaPrek$,
it follows
that after the Criticality Routine has finished due to
$
{\delta}_j \le \CCMu {\chi}_j,
$
we have
\(
  \chik = {\chi}_j \ge \nicefrac{\Deltak}{\CCMu}.
\)
If $\chikPre \ge \CCcritTestchi$, then the criticality loop is not entered and
$\chik = \chikPre \ge \CCcritTestchi$.
Thus, if the second stopping condition does not apply,
it follows from \(
  \thetak \le \CCcritTesttheta
\)
that 
\(
  \chik
  \ge \min
  \left\{\nicefrac{\Deltak}{\CCMu}, \CCcritTestchi\right\}
\).

Finally, note that the models $\mfk, \mgk, \mhk$ are fully linear when the
Criticality Routine has finished after a finite number of iterations
-- even if the radius is slightly increased above ${\delta}_j$ --
because of \cref{ass:fully_linear_models}.

\section{Convergence to Quasi-Stationary Points}%
\label{sec:quasi_convergence}

Note that the problem~\eqref{eqn:itr_tangential} defining
%${\omega}\left(  \vxk; \mfk, \linFeasApproxk, \normIt{{\bullet}} \right)$
$\chik$
is similar to the criticality problem~\eqref{eqn:DemandIdenticalOval},
but we compute the descent step starting
at $\vxNk$ instead of $\vxk$,
the surrogate model gradients are used to determine a model descent step
within the approximated linearized feasible set and
a variable vector norm $\normIt{{\bullet}}$ bounds the problem (which must not necessarily equal the trust-region norm $\normTrIt{{\bullet}}$).
% \begin{enumerate}
%     \item we compute the descent step starting
%     at $\vxNk$ instead of $\vxk$,
%     \item the surrogate models gradients are used to determine a model descent step
%     \item within the approximated linearized feasible set and
%     \item a variable vector norm $\normIt{{\bullet}}$ bounds the problem (which must not necessarily equal the
%     trust-region norm $\normTrIt{{\bullet}}$.)
% \end{enumerate}
For the convergence analysis, the additional uncertainties are dealt with in two steps:
The first part (the remainder of this section)
is concerned with proving that there
is a subsequence $\{\vxkl\}$ of iterates converging to a quasi-stationary point.
That is, it holds that
\[
  \lim_{{\ell}\to \infty} {\theta}( \vxkl ) = 0
  \quad\text{and}\quad
  \lim_{{\ell}\to \infty}
  {\omega}\left(\vxNkl; \mfkl, \linFeasApproxkl, \normInf{{\bullet}}\right) = 0.
\]
Afterwards, \cref{section:kkt_convergence} builds on this result to show convergence
to actual KKT-points.

%\begin{myremark*}
\subsection{Comparison with Related Algorithms}\label{sec:other_algos}
Most results in this section are ``translated'' from their single-objective
pendants in~\cite{fletcher_global_2002} or~\cite{walther2016}.
The latter article is also concerned with inexact surrogates, which have
error bounds that are slightly different from ours and can become exact eventually.
Hence, we cannot simply apply the maximum-scalarization and be done
(unfortunately).
We also have to take care of the Criticality Routine and other, more subtle
differences, such as the iteration-dependent norms.
That is why we have decided to cite the results
from~\cite{fletcher_global_2002,walther2016} explicitly
whenever we introduce a similar one.
If we deem the proofs very similar and easily transferable, we have stated
so.
Otherwise, hints are given on how to adapt them or the proofs are provided wholly for the 
sake of completeness.

Special mention has to be made of the single-objective algorithm in~\cite{eason_trust_2016}
that is further detailed in the dissertation~\cite{eason_trust_2018}.
This algorithm also employs a Criticality Routine and uses fully linear models.
By encoding the surrogate modeling of black-box components as additional constraints
to the original problem and using a nonlinear (even non-quadratic) subproblem
for the computation of the normal step, the convergence analysis becomes easier
due to the resulting inexactness bounds. 
At the time of writing, the nonlinear normal step computation did not suit our particular needs,
but we think it easily possible and very worth-wile to also transfer their approach to 
the multi-objective case.

%\end{myremark*}

\subsection{Final Definitions and Requirements}
To actually investigate limit points of algorithmic iteration sequences we need
\cref{ass:closed_and_bounded}, which guarantees their existence:
\begin{assumption}%
  \label{ass:closed_and_bounded}
  The set $\covering$ is contained in a closed and bounded set.
  \end{assumption}
  A detailed discussion of~\cref{ass:closed_and_bounded} and alternatives
  is given in~\cite{fletcher_global_2002}.
  If the iterates are contained in a bounded domain,
  the true functions (which are continuous according to~\cref{ass:lipschitz_gradients_true})
  are bounded and so are their Lipschitz-continuous gradients.
  With~\cref{ass:fully_linear_models} the models are fully linear and
  have to be bounded at the iterates as well, due to $\DeltaPrek\le\DeltaMax<\infty$:
  \begin{mycorollary}%
      \label{thm:bounded_model_gradients}
      If~\cref{%
      ass:closed_and_bounded,%
      ass:lipschitz_gradients_true,%
      ass:fully_linear_models}
      and $\DeltaPrek \le \DeltaMax < \infty$ hold,
      then the norm of all function and model gradients is uniformly bounded above,
      i.e., there is a positive constant $\CCubj > 0$ such that for
      \(
        {\varphi} \in \left\{
          \fL_1, \ldots , \fL_\fDim,
          \gL_1, \ldots , \gL_\gDim,
          \hL_1, \ldots , \hL_\hDim,
        \right\}
     \),
     for all $\kk\in\mathbb{N}_0$ and
     \(
       {\varphi}\KK \in \left\{
         \mfL\KK_1, \ldots , \mfL\KK_\fDim,
         \mgL\KK_1, \ldots , \mgL\KK_\gDim,
         \mhL\KK_1, \ldots , \mhL\KK_\hDim,
       \right\}
     \)
     % \[
    % \begin{aligned}
    % {\varphi} &\in \left\{
    %     \fL_1, \ldots , \fL_\fDim,
    %     \gL_1, \ldots , \gL_\gDim,
    %     \hL_1, \ldots , \hL_\hDim,
    % \right\}
    % &\text{and}
    % \\
    % {\varphi}\KK &\in \left\{
    %     \mfL\KK_1, \ldots , \mfL\KK_\fDim,
    %     \mgL\KK_1, \ldots , \mgL\KK_\gDim,
    %     \mhL\KK_1, \ldots , \mhL\KK_\hDim,
    %     \right\}
    % \end{aligned}
    % \]
    it holds that
    \begin{equation*}
        \norm{\grad {\varphi}(\vxk)}\le \CCubj
        \quad\text{and}\quad
        \normInf{ \grad {\varphi}\KK(\vxk) } \le \CCubj.
        \label{eqn:def_CCubj}
  \end{equation*}
\end{mycorollary}

Throughout this section, we consider the following iteration index sets:
\begin{mydefinition}
    The set of ``successful'' iteration indices is
    \(
       \acceptIndices = \{ \kk: \vxNext = \vxSk \}.
    \)
    The set of restoration indices is defined as
    \(
    \restoreIndices = \{ \kk:
        \vnk \text{ does not exist
        or~\eqref{eqn:compatible_normal_step}
        fails, i.e., $\vnk$ is incompatible}
    \}.
    \)
    Finally, the set of filter indices (indices of iterations that modify the filter) is
    \(
    \filterIndices = \{
        \kk: \vxk \text{ is added to the filter }
    \}
    \supseteq \restoreIndices.
    \)
\end{mydefinition}

Whenever we require ``the Criticality Routine to finish finitely'' in subsequent
statements, then we want it to finish after a finite number of iterations
and explicitly include the case that the routine is not even entered due to
the Criticality Test failing in step~\ref{algo:crit_test} of the algorithm.

\subsection{Convergence Analysis}
Similar to~\cite{walther2016} we can state accuracy requirements that bound
the linearized constraint violation of the steps $\vnk$ and $\vsk$ by the
trust-region radius $\Deltak$:
\begin{mylemma}[Accuracy Requirements similar to {\cite[A. 2.4]{walther2016}}]%
    \label{thm:accuracy_requirements}
    Suppose~\cref{%
    ass:covering,%
    ass:lipschitz_gradients_true,%
    ass:fully_linear_models,%
    ass:equivalent_norms%
    } hold and that $\kk\in \mathbb{N}_0$ is an iteration index.
    Then there is a constant $\CCaccN > 0$ such that,
    if $\vnk$ exists as the solution to~\eqref{eqn:itr_normal},
    it holds that
    \begin{equation}
        \max \left\{
          %\begin{aligned}
          %&
          \max_{\ell} \left| \hL_{\ell}(\vxk) + \dh_{\ell}(\vxk) \cdot \vnk \right|,
          %\\
          %&
          \max_{\ell} \gL_{\ell}(\vxk) + \dg_{\ell}(\vxk) \cdot \vnk
          % \end{aligned}
        \right\}
        \le
        %\CCnorm \max\left\{ \CCflDiffG, \CCflDiffH \right\}
        \CCaccN
        \DeltaPrek \normTrIt{\vnk}.
        \label{eqn:recantsKnuts}
    \end{equation}
    For any $\kk \in \mathbb{N}_0$
    for which $\vnk$ exists and satisfies
    $\normTrIt{\vnk}\le \DeltaPrek$
    it also holds that
    \begin{align}
        \max \left\{
            \max_{\ell} \left| \hL_{\ell}(\vxk) + \dh_{\ell}(\vxk) \cdot \vnk \right|,
            \max_{\ell} \gL_{\ell}(\vxk) + \dg_{\ell}(\vxk) \cdot \vnk
        \right\}
        &\le \CCaccN \DeltaPrek^2,
        \label{eqn:accuracy_normal}
        \intertext{%
        and if the Criticality Routine finishes finitely and
        $\vsk$ is the step $\vnk + \stepsizek \vdk$ with $\vdk$ computed
        using~\eqref{eqn:itr_tangential} and $\normTrIt{\vsk} \le \Deltak$,
        we have}
        \max \left\{
            \max_{\ell} \left| \hL_{\ell}(\vxk) + \dh_{\ell}(\vxk) \cdot \vsk \right|,
            \max_{\ell} \gL_{\ell}(\vxk) + \dg_{\ell}(\vxk) \cdot \vsk
        \right\}
        &\le \CCaccS \Deltak^2.
        \label{eqn:accuracy_step}
    \end{align}
\end{mylemma}
\begin{proof}
    The constraints of~\eqref{eqn:itr_normal} ensure that
    \begin{equation}
    \mhk(\vxk) + \DmH(\vxk) \cdot \vnk = \vZ
    \quad\text{and}\quad
    \mgk(\vxk) + \DmG(\vxk) \cdot \vnk \le \vZ.
    \label{eqn:centimetersHectic}
    \end{equation}
    %and we have assumed that
    %$\normInf{\vnk} \le \CCnorm \normTrIt{\vnk} \le \CCnorm \DeltaPrek$.\\
    Because the models are interpolating (\cref{ass:fully_linear_models})
    we also have 
    \(\mhk(\vxk) = \vh(\vxk)\)
    and 
    \(\mgk(\vxk) = \vg(\vxk)\).
    Hence, additionally using the Cauchy-Schwartz inequality and the error-bounds
    of the models (for the preliminary radius $\DeltaPrek$),
    we obtain
    \[
    \begin{aligned}
    \max_{\ell} \abs{
        \hL_{\ell} + \grad h_{\ell} \cdot \vnk
    }
    &=
    \abs{
        \hL_{\ell} + \grad h_{\ell} \cdot \vnk - 0
    }
    \stackrel{\eqref{eqn:centimetersHectic}}=
    \abs{
        (\dh_{\ell} - \dmh\KK_{\ell}) \cdot \vnk
    }
    \\
    &%
    \le
    \normTwo{
        (\dh_{\ell} - \dmh_{\ell}\KK)
    } 
    \normTwo{\vnk}
    \le
    \CCflDiffH \DeltaPrek \CCnorm \normTrIt{\vnk},
    \end{aligned}
    \]
    where the argument $\vxk$ was dropped to improve readability.
    Similarly, we find that
    \[
    \max_{\ell} \gL_{\ell}(\vxk) + \dg_{\ell}(\vxk) \cdot \vnk
    \le \CCflDiffG \DeltaPrek \CCnorm \normTrIt{\vnk}
    \]
    and~\eqref{eqn:recantsKnuts} follows with
	$\CCaccN = \CCnorm \max\{ \CCflDiffH, \CCflDiffG \}$.
    Moreover, $\normTrIt{\vnk} \le \DeltaPrek$ leads to~\eqref{eqn:accuracy_normal}.
    \\
    The second inequality~\eqref{eqn:accuracy_step} is derived analogously,
    respecting the fact that the step-size is chosen
	  so that $\normTrIt{\vsk} \le \Deltak$ and that $\vxk + \vsk$ also satisfies
	  the approximated linearized constraints.
\end{proof}

From the preceding accuracy results, a bound on the constraint violation can be derived:
\begin{mylemma}[{\cite[Lemma 4.3]{walther2016},\cite[Lemma 3.1]{fletcher_global_2002}}]%, Lemma 3.1 in~\cite{fletcher_global_2002}]%
    \label{thm:bound_on_constraint_violation1}
    Assume that the algorithm is applied to~\eqref{eqn:mop_nonlin} and that
    \cref{%
    ass:covering,%
    ass:closed_and_bounded,%
    ass:fully_linear_models,%
    ass:normal_step_exists_bounded,%
    ass:lipschitz_gradients_true,%
    %ass:hess_norm_bounded,%
    %ass:sufficient_decrease,
    ass:equivalent_norms%
    } hold.
    %If $\thetak \le \CCdeltan$,
    %then the normal step exists.
    There is a constant $\CCubj > 0$
    (independent of $\kk$) such that, if the normal step exists,
    it holds that
    \begin{equation}
        \thetak \le \left( \CCaccN \DeltaPrek + \CCubj \CCnorm \right) \normTrIt{\vnk}.
        \label{eqn:balladedAdzukis}
    \end{equation}
    %\begin{equation}
    %    \thetak \le \CCaccN \DeltaPrek^2 + \CCubj \normTwo{\vnk}.
    %    \label{eqn:bound_on_constraint_violation1}
    %\end{equation}
\end{mylemma}
\begin{proof}
  The proof works very similar to that of Lemma 4.3
  in~\cite{walther2016}.
  We only need to take~\cref{thm:accuracy_requirements} into account for the
  different accuracy requirements.
\end{proof}

When the normal step is compatible, i.e., for iterations with $\kk \notin \restoreIndices$,
we can further refine the bound on $\thetak$:
\begin{mylemma}[Lemma 4.4 in~\cite{walther2016}, Lemma 3.4 in~\cite{fletcher_global_2002}]%
    \label{thm:bound_on_constraint_violation2}
    Suppose that the algorithm is applied to~\eqref{eqn:mop_nonlin} and that
    \cref{%
        ass:covering,%
        ass:closed_and_bounded,%
        ass:lipschitz_gradients_true,%
        ass:fully_linear_models,%
        ass:normal_step_exists_bounded,%
        ass:lipschitz_gradients_true,%
        ass:equivalent_norms%
    } hold.
    Suppose further that $\kk\notin\restoreIndices$ and
    that~\eqref{eqn:balladedAdzukis} holds.
    Then there is a constant $\CCthetaub > 0$ such that
    \begin{equation}
        \thetak \le \CCthetaub \max \left\{ \DeltaPrek^{1+{\mu}}, \DeltaPrek^2 \right\}
        \quad\text{and}\quad
        {\theta}\left( \vxk + \vsk \right) \le \CCthetaub \Deltak^2.
        \label{eqn:shoalinessesDespicablenesses}
    \end{equation}
\end{mylemma}

\begin{proof}
  Since $\kk\notin\restoreIndices$, one obtains
  from~\eqref{eqn:balladedAdzukis} and
  from~\eqref{eqn:compatible_normal_step} that
  \begin{align*}
    \thetak
    &\le
    \left( \CCaccN \DeltaPrek + \CCubj \CCnorm \right) \normTrIt{\vnk}
    \le
    \left( \CCaccN \DeltaPrek + \CCubj \CCnorm \right)
    \CCDelta \DeltaPrek \min \left\{ 1 , \CCmu \DeltaPrek^\mu \right\}
    \\
    &\le
    \CCDelta
    \left( \CCaccN \DeltaPrek^2 + \CCubj \CCnorm \CCmu \DeltaPrek^{1+{\mu}}\right)
%    \\
%    &%
    \le
    \underbrace{\CCDelta
    \left(
        \CCaccN
        + \CCubj \CCnorm \CCmu
    \right)}_{=:\CCthetaub^*\text{ (const.)}} \max \left\{ \DeltaPrek^{1+\mu}, \DeltaPrek^2 \right\}.
  \end{align*}
  For the second bound let $c_{\ell}$ be any constraint component from $\vh$ or $\vg$.
  Due to~\cref{ass:covering,ass:lipschitz_gradients_true} we can
  construct a Taylor polynomial and use the Mean-Value-Theorem
  to obtain
  \begin{equation}
    c_{\ell}( \vxk + \vsk )
    =
    c_{\ell}(\vxk) + \tran{\grad c_{\ell}(\vxk)} \vsk +
      \frac{1}{2}\langle \hess c_{\ell}(\ve {\xi}) \vsk, \vsk\rangle
    \label{eqn:proof_bound_on_constraint_violation2_1}
  \end{equation}
  for some $\ve {\xi}$ on the line segment from $\vxk$ to $\vxk + \vsk$.
  From~\eqref{eqn:accuracy_step} we know that
  \(\abs{ c_{\ell}(\vxk) + \tran{\grad c_{\ell}(\vxk)} \vsk}
    \le
    \CCaccS \Deltak^2.
  \)
  With~\cref{ass:closed_and_bounded,ass:lipschitz_gradients_true}
  the norm off all constraint Hessians in
  \cref{eqn:proof_bound_on_constraint_violation2_1}
  can be bounded from above globally (say, by $2{\kappa}_{H} > 0$) so that
  the Cauchy-Schwartz inequality gives
  \[
    \frac{1}{2}\langle \hess c_{\ell}(\ve {\xi}) \vsk, \vsk\rangle
    \le
    \frac{1}{2}\normTwo{\vsk}^2 \normTwo{\hess c_{\ell}(\ve {\xi})}
    \le
    {\kappa}_{H} \Deltak^2.
  \]
  Hence,~\eqref{eqn:shoalinessesDespicablenesses} follows with
  $
  \CCthetaub =: \max \left\{
      \CCthetaub^*,
      \CCaccS + {\kappa}_{H}
      \right\}.
  $
\end{proof}

The requirements for the accuracy requirements to hold are also met within the
Criticality Routine,
allowing for the following corollary:
\begin{mycorollary}
    Let $\kk \in \mathbb{N}_0$ be an iteration index and suppose the same requirements
    as in~\cref{thm:bound_on_constraint_violation1} hold.
    If the Criticality Routine is entered,
    then for any sub-iteration index $j\in \mathbb N_0$, it holds that
    % \begin{equation}
    %     \max \left\{
    %         \max_{\ell} \left| \hL_{\ell}(\vxk) + \dh_{\ell}(\vxk) \cdot \vnk_j \right|,
    %         \max_{\ell} \gL_{\ell}(\vxk) + \dg_{\ell}(\vxk) \cdot \vnk_j
    %     \right\}
    %     \le \CCaccN {\delta}_j^2
%        \label{eqn:accuracy_normal_inloop}
    %  \quad\text{and}\quad        
    \begin{equation}
      \thetak
        \le
        \CCaccN
            {\delta}_j^2
            +
            \CCubj \normTwo{\vnk_j},
        \label{eqn:bound_on_constraint_violation1_inloop}
    \end{equation}
    If the Criticality Routine is not entered or if it has
    completed finitely, then it also holds that
    \begin{align}
        \thetak
        &\le
        \CCaccN \Deltak^2 + \CCubj \normTwo{\vnk},
        \label{eqn:bound_on_constraint_violation1_finCrit}
        \shortintertext{and if $\kk\notin\restoreIndices$, then}
        \thetak
        &\le
        \CCthetaub \max \left\{ \Deltak^{1+{\mu}}, \Deltak^2 \right\}.
        \label{eqn:shoalinessesDespicablenesses_finCrit}
    \end{align}
\end{mycorollary}

\subsubsection{Convergence in the Criticality Routine}

Our first important observation is that quasi-criticality is approximated if the
criticality loop runs infinitely:
\begin{mylemma}[{\cite[Lemma 8]{ours}}]%
    \label{thm:infinite_crit_routine}
    Assume that
    \cref{%
    ass:covering,%
    ass:closed_and_bounded,%
    ass:fully_linear_models,%
    ass:normal_step_exists_bounded,%
    ass:lipschitz_gradients_true,%
    ass:equivalent_norms%
    }
    hold.
    Denote by
    ${\chi}_j$ the inexact criticality value
    from~\cref{thm:ScaryPitchOld} in the criticality subroutine.
    If the criticality routine runs infinitely ($j\to \infty$) at $\vxk$,
    then
    \[{\theta}\left( \vxk \right) = 0
    \text{  and  }
    \lim_{j\to \infty} {\delta}_j = 0
    \text{  and  }
    \lim_{j\to \infty}{\chi}_j = 0.
    \]
	\end{mylemma}

\begin{proof}
    As in the description of the Criticality Routine, let $j$ be the sub-iteration index
    and $\delta_j^+ := \CCalpha^j \DeltaPrek$,
    and $\delta_j = {\delta}_{j-1}^+$ after the index $j$ has been increased.
    Also, denote by ${\chi}_j$ the updated doubly inexact criticality measure from
    step~\ref{critAlgo:crit_update}.
    There are two logically disjunct stopping criteria for the Criticality Routine.
    \\
    Firstly, the Criticality Routine may stop if
    \(
    \CCMu {\chi}_j
    \ge
    {\delta}_j
    \)
    for some $j\in \mathbb{N}_0$.
    Conversely, if the routine loops infinitely, then it must hold that
    \[
    {\chi}_j
	<
	\frac{{\delta}_j}{\CCMu_j}
    =
    \CCalpha^j \frac{\DeltaPrek}{\CCMu}
	\qquad \text{for all $j\in \mathbb{N}$.}
    \]
    Because of $\CCalpha\in (0,1)$ the right side goes to zero
    and then the second limit in~\cref{thm:infinite_crit_routine} also follows.
    \\
    Secondly, the routine stops if no compatible step $\vn_j$ for the
    preliminary radius ${\delta}_j^+$ can be found anymore.
    Vice versa, if the routine runs infinitely, there always is a compatible step with
    \[
    \normTwo{\vn_j}
    \le
    \CCnorm
    \CCDelta {\delta}_j
    \min \left\{
        1,
        \CCmu {\delta}_j^{\mu}
    \right\}.
    \]
    This implies $\vnk_j \to 0$ for $j\to \infty$.
    From~\eqref{eqn:bound_on_constraint_violation1_inloop}
    it then follows that it must already hold that ${\theta}(\vxk) = 0$.
\end{proof}

\subsubsection{Infinitely Many Filter Iterations}

Outside of the criticality loop, we can show that feasibility is approached if the number of
filter iterations is infinite:
\begin{mylemma}[{\cite[Lemma 3.3]{fletcher_global_2002}}]%
    \label{thm:infinite_filter_iters_imply_feasibility}
    Suppose that the algorithm is applied to~\eqref{eqn:mop_nonlin} and that
    \cref{%
        ass:covering,%
        ass:lipschitz_gradients_true,%
        ass:closed_and_bounded,%
    }
    hold.
    If $\abs{\filterIndices} = \infty$, then
    \(
    \lim_{\substack{\kk \to \infty,\\ \kk\in \filterIndices}} \thetak = 0.
    \)
\end{mylemma}
\begin{proof}
    The proof is exactly the same as in~\cite{fletcher_global_2002},
    only $f({\bullet})$ has to be substituted by $\scalf{{\bullet}}$.
    It also works with a $\fDim + 1$ dimensional filter that uses $\vf$ instead of $\scalf{{\bullet}}$.
\end{proof}

The next result shows that feasibility is approached for 
\emph{any} subsequence of iteration indices.
The result is a bit stronger than that found in~\cite{fletcher_global_2002} and not
necessarily needed for the convergence proofs.
A similar theorem can be found in~\cite{wachter_line_2005}.
\begin{mylemma}%
    \label{mythm:infinite_filter_iters_imply_feasibility}
    Suppose that the algorithm is applied to~\eqref{eqn:mop_nonlin} and that
    \cref{ass:covering,ass:closed_and_bounded,ass:lipschitz_gradients_true}
    hold.
    If $\abs{\filterIndices} = \infty$, then
    \(
    \lim_{\kk \to \infty} \thetak = 0.
    \)
\end{mylemma}
\begin{proof}
    For a contradiction, assume that there is a subsequence of indices $\{\kk_{\ell}\}$
    for which the constraint violation is bounded below
    \begin{equation}
    {\theta}\itthetaCmd{\kk_{\ell}} \ge {\varepsilon} > 0\qquad \forall k.
    \label{eqn:perturbationsBibliopegy}
    \end{equation}
    From~\cref{thm:infinite_filter_iters_imply_feasibility} we know that there is some
    $\kk_0$ such that
    \begin{equation}
    {\theta}\itthetaCmd{\kk} < {\varepsilon}
    \quad\text{for all $%
    \kk\in \filterIndices$ with
    $\kk\ge \kk_0$.}
    \label{eqn:humanizesTemporizers}
    \end{equation}
    Because there are infinitely many indices both in $\{\kk_{\ell}\}$ and
    $\filterIndices$, each $\kk\in \{\kk_{\ell}\}$
    (except maybe the first)
    must lie between a largest filter index ${\kappa}_1\in \filterIndices$ and a
    smallest index ${\kappa}_2\in\filterIndices$, respectively,
    ${\kappa}_1(\kk) \le \kk \le {\kappa}_2(\kk)$,
    and from~\eqref{eqn:perturbationsBibliopegy} and~\eqref{eqn:humanizesTemporizers}
    we can deduce that
    there is a $\kk_1 \ge \kk_0$ such that
    \[
    \kk_0 < {\kappa}_1(\kk) < \kk < {\kappa}_2(\kk)
    \quad \text{for all $\kk\in \{\kk_{\ell}\}$ with
    $\kk \ge \kk_1$.}
    \]
    Let $\kk\in \{\kk_{\ell}\}$ be an index with $\kk \ge \kk_1$.
    For all indices ${\kappa}$ with ${\kappa}_1(\kk) < {\kappa} < {\kappa}_2(\kk)$,
    the iteration is not a filter iteration and hence the function $\scalShort$ is
    ``monotonic'':
    \[
    \scalShort\left(\vx\itCmd{{\kappa}}\right)
    -
    \scalShort\left(\vx\itCmd{{\kappa}+1}\right)
    \ge 0.
    \]
    % There then must be a largest successful index ${\kappa}_1^*(\kk)$ with
    % ${\kappa}_1(\kk) \le {\kappa}_1^*(\kk) < \kk$,
    % such that the constraint violation is increased from below ${\varepsilon}$ to above ${\varepsilon}$,
    % \[
    % {\theta}\itthetaCmd{{\kappa}_1^*(\kk)}
    % <
    % {\varepsilon}
    % \le
    % {\theta}\itthetaCmd{{\kappa}_1^*(\kk) + 1},
    % \]
    % and from the filter mechanism we now that then the objective reduction has
    % to be significant:
    % \[
    % \scalShort\left( \vx\itCmd{{\kappa}_1^*(\kk) + 1} \right)
    % \le
    % \scalShort\left( \vx\itCmd{{\kappa}_1^*(\kk)} \right)
    %     - \filterConst {\theta}\itthetaCmd{{\kappa}_1^*(\kk)}.
    % \]
    Moreover, there must be some smallest successful index ${\kappa}_2^*(\kk)$
    with
    $\kk \le {\kappa}_2^*(\kk) < {\kappa}_2(\kk)$,
    reducing the constraint violation from above ${\varepsilon}$ to below ${\varepsilon}$ again.
    ${\kappa}_2^*(\kk)$ is not a filter index, so~\eqref{eqn:algo_model_dec}
    succeeds and ${\rho}\itCmd{{\kappa}_2^*(\kk)} \ge \CCaccept > 0$,
    resulting in
    \begin{equation}
    \scalShort\left( \vx\itCmd{{\kappa}_2^*(\kk)} \right)
    -
    \scalShort\left( \vx\itCmd{{\kappa}_2^*(\kk) + 1} \right)
    \ge
    \CCaccept
    {\kappa}_{\theta}
    {\theta}\itthetaCmd{{\kappa}_2^*(\kk)}^{\psi}
    \ge
    \CCaccept
    {\kappa}_{\theta}
    {\varepsilon}^{\psi} > 0.
    \label{eqn:samovarDemiurgic}
    \end{equation}
    The function $\scalShort$ might be increased eventually
    in the filter iteration ${\kappa}_2(\kk)$ (when accepting
    $\vx\itCmd{{\kappa}_2(\kk) +1}$), but
    due to the monotonicity of the intermediate iterations it holds that
    \(
      \scalShort\left(\vx\itCmd{{\kappa}_2^*(\kk) + 1}\right)
      -
      \scalShort\left(\vx\itCmd{{\kappa}_2(\kk)}\right)
      \ge 0
    \)
    and thus
    \begin{equation}
    %\stackrel{\eqref{eqn:samovarDemiurgic}}\Rightarrow\quad
    \scalShort\left( \vx\itCmd{{\kappa}_2^*(\kk)} \right)
    -
    \scalShort\left( \vx\itCmd{{\kappa}_2(\kk)} \right)
    \ge
    \CCaccept
    {\kappa}_{\theta}
    {\varepsilon}^{\psi} > 0,
    \label{eqn:bulgineOped}
    \end{equation}
    and $\vx\itCmd{{\kappa}_2(\kk)}$ is added to the filter.

    For $\kk$ from above, let $\kk_+ \in \{\kk_{\ell}\}$ be
    the smallest iteration index following $\kk$
    such that it holds for the enclosing filter indices that
    ${\kappa}_2(\kk) \le {\kappa}_1(\kk_+)$.
    There then must be a largest successful index ${\kappa}_1^*(\kk_+)$ with
    ${\kappa}_1(\kk_+) \le {\kappa}_1^*(\kk_+) < \kk_+$,
    such that the constraint violation is increased from below ${\varepsilon}$ to above ${\varepsilon}$,
    \[
    {\theta}\itthetaCmd{{\kappa}_1^*(\kk_+)}
    <
    {\varepsilon}
    \le
    {\theta}\itthetaCmd{{\kappa}_1^*(\kk_+) + 1},
    \]
    and from the filter mechanism we know that then the objective reduction has
    to be significant compared to \emph{all} points in the filter with a smaller constraint violation,
    including that with index ${\kappa}_2(\kk)$ and ${\theta}\itthetaCmd{{\kappa}_2(\kk)}<{\varepsilon}$:
    \[
    \scalShort\left( \vx\itCmd{{\kappa}_2(\kk)} \right)
    -
    \scalShort\left( \vx\itCmd{{\kappa}_1^*(\kk_+) + 1} \right)
    \ge
    \filterConst {\theta}\itthetaCmd{{\kappa}_2(\kk)}
    \ge 0.
    \]

    Using the same notation as above,
    (i.e., ${\kappa}_2^*(\kk_+)$ as in~\eqref{eqn:samovarDemiurgic}),
    we deduce from the monotonicity of the intermediate iterations that
    \( 
      \scalShort\left( \vx\itCmd{{\kappa}_2(\kk)} \right)
      -
      \scalShort\left( \vx\itCmd{{\kappa}_2^*(\kk_+)} \right)
      \ge 0
    \)
        %\nonumber
        %\shortintertext{and thus, with {\eqref{eqn:bulgineOped}},}
    and thus, with {\eqref{eqn:bulgineOped}},
    \begin{equation}
        \scalShort\left( \vx\itCmd{{\kappa}_2^*(\kk)} \right)
        -
        \scalShort\left( \vx\itCmd{{\kappa}_2^*(\kk_+)} \right)
        \ge
        \CCaccept
        {\kappa}_{\theta}
        {\varepsilon}^{\psi} > 0.
        \label{eqn:aluminousVixens}
    \end{equation}

    By repeating the above procedure,
    we see that it is possible to construct an infinite subsequence
    $\{{\kappa}_j\}$ of
    iteration indices from the ${\kappa}_2^*$ values
    (of indices from $\{\kk_{\ell}\}$),
    for which $\scalShort(\vx\itCmd{{\kappa}_j})$ is
    strictly monotonically decreasing with a guaranteed (constant) objective
    reduction~\eqref{eqn:aluminousVixens}.
    This is a contradiction to $\scalShort$ being bounded below as per
    \cref{ass:closed_and_bounded,ass:lipschitz_gradients_true}.
\end{proof}

What follows next is a series of auxiliary lemmata to finally show
convergence of the inexact criticality measure when $\abs{\filterIndices} = \infty$.
First, we have a bound on the surrogate objective change along the normal step:
\begin{mylemma}[{Similar to a bound in~\cite[Lemma 3.5]{fletcher_global_2002}}]%
    \label{thm:normal_step_decrease_ub}
    Suppose the algorithm is applied to~\eqref{eqn:mop_nonlin}, that
    \cref{%
    ass:covering,%
    ass:closed_and_bounded,%
    ass:lipschitz_gradients_true,%
    ass:fully_linear_models,%
    ass:equivalent_norms,%
    ass:hess_norm_bounded,%
    %ass:sufficient_decrease%
    } hold and that the Criticality Routine does not run infinitely.
    Suppose further that $\kk\notin\restoreIndices$.
    Then it holds that
    \begin{align}
        \abs{
            \scalk(\vxk) - \scalk(\vxNk)
        }
        =
        \abs{
            \scalk(\vxNk) - \scalk(\vxk)
        }
        &\le
        \CCubj \normTwo{{\vnk}}
        + \frac{1}{2}\CChessF \normTwo{{\vnk}}^2
        \label{eqn:deafenedPonytailed}.
    \end{align}
\end{mylemma}

\begin{proof}
    There are maximizing indices ${\ell},j\in\{1,\ldots ,\fDim\}$ such that
    \[
    \scalk(\vxNk) - \scalk(\vxk)
    =
    \mfL_{\ell}\KK(\vxNk) - \mfL_j\KK(\vxk)
    %\]
    %and because the indices are maximizing at $\vxNk$ and $\vxk$, respectively,
    %it holds that
    %\begin{align*}
    \le
      \scalk(\vxNk) - \scalk(\vxk)
      \le
      \mfL_{\ell}\KK(\vxNk) - \mfL_{\ell}\KK(\vxk)
      .
    \]
    Using a 2nd degree Taylor approximation of $\mfL_{\ell}\KK$ around $\vxk$
    at $\vxNk$
    %,
    % \[
    % \mfL_{\ell}\KK(\vxNk) =
    % \mfL_{\ell}\KK(\vxk)
    % + \langle
    % \grad \mfL_{\ell}\KK(\vxk),
    % \vnk
    % \rangle
    % +
    % \frac{1}{2}
    % \langle
    % \hess \mfL_{\ell}\KK(\vxk) \vnk,
    % \vnk
    % \rangle,
    % \]
    results in
    \begin{align*}
        \scalk(\vxNk) - \scalk(\vxk)
        &\le
        \mfL_{\ell}\KK(\vxNk) - \mfL_{\ell}\KK(\vxk)
        %\nonumber \tag{!!\theequation}\stepcounter{equation}
        \\
        &\le
        \abs{
          \langle
          \grad \mfL_{\ell}\KK(\vxk),
          \vnk
          \rangle
        }
          +
        \abs{
          \frac{1}{2}
          \langle
          \hess \mfL_{\ell}\KK(\vxk) \vnk,
          \vnk
          \rangle
        }
        % \nonumber \tag{!!\theequation}\stepcounter{equation}
        \\
        &\le
        \normTwo{\grad \mfk_{\ell}(\vxk)}\normTwo{\vnk}
        + \frac{1}{2}\normTwo{\hess \mfk_{\ell}(\vxk)} \normTwo{{\vnk}}^2
        %\nonumber \tag{!!\theequation}\stepcounter{equation}
        \\
        &%
        %\stackrel{\mathclap{%
        %    \text{\cref{ass:hess_norm_bounded}},%
        %    \text{\cref{thm:bounded_model_gradients}}%
        %}}%
        \le
        \CCubj \normTwo{{\vnk}}
        + \frac{1}{2}\CChessF \normTwo{{\vnk}}^2,
        % \tag{!!\theequation}\stepcounter{equation}
    \end{align*}
    where the last inequality comes
    from~\cref{ass:hess_norm_bounded,thm:bounded_model_gradients}.
    Analogously, we can show
	\begin{equation*}
	%&\mfL_j\KK(\vxNk) - \mfL_j\KK(\vxk) \le \scalk(\vxNk) - \scalk(\vxk)\\
	%\Leftrightarrow\quad
	%&
  -\left(\scalk(\vxNk) - \scalk(\vxk)\right) 
  \le 
  \mfL_j\KK(\vxk) - \mfL_j\KK(\vxNk)
	\le
    \CCubj \normTwo{{\vnk}}
     + \frac{1}{2}\CChessF \normTwo{{\vnk}}^2.
  \end{equation*}
  %The bound~\eqref{eqn:deafenedPonytailed} follows.
\end{proof}

The next lemma provides a sufficient decrease bound in case that the doubly
inexact criticality measure is bounded below, and the radius is sufficiently small:
\begin{mylemma}[{\cite[Lemma 3.5]{fletcher_global_2002}}]%
    \label{thm:sufficient_decrease_when_crit_bounded}
    Suppose the algorithm is applied to~\eqref{eqn:mop_nonlin}, that
    \cref{%
    ass:covering,%
    ass:closed_and_bounded,%
    ass:lipschitz_gradients_true,%
    ass:fully_linear_models,%
    ass:equivalent_norms,%
    ass:hess_norm_bounded,%
    ass:sufficient_decrease%
    } hold and that the Criticality Routine does not run infinitely.
    Suppose further, that $\kk\notin\restoreIndices$,
    %that $\vnk$ is compatible,
    that
    \begin{equation}
    \chik \ge \frac{1}{\CCcrit} \chik_2 > \frac{1}{\CCcrit} \CCepschi
    \tag{LB}
    \label{eqn:porcellanizingTrichromat}
    \end{equation}
    for some $\CCepschi>0$ and that
    \begin{equation}%
      \label{eqn:pentanolFice}
      \Deltak
      \le
      {\delta}_m :=
      \min \Bigg\{ \;
            \frac{
                \CCepschi
            }{
              \CCcritDenom
            },
          {\left(
              \frac{
                  2 \CCubj
              }{
                  \CChessF\CCDelta\CCmu
              }
          \right)}^{\frac{1}{1+\mu}},
          \left(
              \frac{
                  \CCsd \CCepschi
              }
              {
                  4\CCubj \CCDelta \CCmu
              }^{\frac{1}{1+\mu}}
          \right)
      \; \Bigg\}.
    \end{equation}
    Then
    \begin{equation}
      \scalk(\vxk) - \scalk(\vxk + \vsk)
      \ge
      \frac{1}{2}
      \CCsd \CCepschi \Deltak.
    \label{eqn:foreversTalced}
    \end{equation}
\end{mylemma}
\begin{proof}
  The proof is similar to that given in~\cite{fletcher_global_2002}.
  The bounds~\eqref{eqn:porcellanizingTrichromat} and~\eqref{eqn:pentanolFice}
  are plugged into the sufficient decrease equation
  of~\cref{ass:sufficient_decrease}.
  The bound~\eqref{eqn:foreversTalced} follows
  from~\cref{thm:normal_step_decrease_ub},~\eqref{eqn:pentanolFice}
  and the fact that $\vnk$ is compatible with some tedious algebra.
\end{proof}

Under similar conditions as in~\cref{thm:sufficient_decrease_when_crit_bounded} the iteration will
be successful:
\begin{mylemma}[{\cite[Lemma 3.6]{fletcher_global_2002}}]%
    \label{thm:success_when_crit_bounded}
    Suppose the algorithm is applied to~\eqref{eqn:mop_nonlin}, that
    \cref{%
    ass:covering,%
    ass:closed_and_bounded,%
%    ass:normal_step_exists_bounded,%
    ass:lipschitz_gradients_true,%
    ass:fully_linear_models,%
    ass:equivalent_norms,%
    ass:hess_norm_bounded,%
    ass:sufficient_decrease%
    } hold and that the Criticality Routine does not run infinitely.
    Suppose further, that $\kk\notin\restoreIndices$,
    that~\eqref{eqn:porcellanizingTrichromat} holds again and that
    \begin{equation}
        \Deltak \le \delta_{\rho} :=
        \min
        \left\{
        \delta_m,
        \frac{
            (1-\CCsuccess) \CCsd \CCepschi
        }
        {
            2 \CCflF
        }
        \right\}.
        \label{eqn:yogurt_plant_cup}
    \end{equation}
    Then $\KK$ is successful, that is,
    \(
    \rho\KK \ge \CCsuccess.
    \)
\end{mylemma}
\begin{proof}
  We can again follow the proof in~\cite{fletcher_global_2002},
  because from~\cref{ass:fully_linear_models} we can conclude 
  (cf.~\cite[Lemma 4.16]{thomann_paper})
  that the model error bound holds also for the scalarization:
  \begin{equation}
    \abs{
      \scal{\vf}(\ve {\xi}) -
      \scal{\mfk}(\ve {\xi})
    }
    =
    \abs{
      \scalShort(\ve {\xi}) -
      \scalk(\ve {\xi})
    } \le \CCflF \Deltak^2 \quad \forall \ve {\xi} \in \trk.
    \label{eqn:fully_linear_scalarization}
  \end{equation}
\end{proof}
Further, if the radius is sufficiently small, then the test~\eqref{eqn:algo_model_dec} will succeed
(prohibiting a ${\theta}$-iteration):
\begin{mylemma}[{\cite[Lemma 3.7]{fletcher_global_2002}}]%
    \label{thm:model_dec_test_succeeds}
    Suppose the algorithm is applied to~\eqref{eqn:mop_nonlin}, that
    \cref{%
    ass:covering,%
    ass:closed_and_bounded,%
    ass:normal_step_exists_bounded,%
    ass:fully_linear_models,%
    ass:lipschitz_gradients_true,%
    ass:equivalent_norms,%
    ass:hess_norm_bounded,%
    ass:sufficient_decrease%
    } hold and that the Criticality Routine does not run infinitely.
    Suppose further, that $\kk\notin\restoreIndices$,
    that~\eqref{eqn:porcellanizingTrichromat} holds again,
    that $\thetak \le \CCdeltan$,
    %\eqref{eqn:shoalinessesDespicablenesses_finCrit} holds
	and that
    \begin{equation}
        \Deltak \le
        {\delta}_f :=
        \min
        \left\{
            {\delta}_m,
            1,
            {\left(
				\frac{
					\CCsd \CCepschi }{
					2{\kappa}_{\theta} \CCthetaub^{{\psi}}
				}
            \right)}^{\frac{1}{{\psi}(1 + \mu)-1}}
        \right\}.
        \label{eqn:edison_tesla_boxfight}
    \end{equation}
    Then~\eqref{eqn:algo_model_dec} succeeds, i.e.,
    \(
    \scalk(\vxk) -
        \scalk(\vxSk)
    \ge
    {\kappa}_{\theta} \thetak^{\psi}.
    \)
\end{mylemma}
% \todo[inline]{%
%   The proof of~\cite{fletcher_global_2002} has to be slightly modified because
%   of the $\max$ in~\eqref{eqn:shoalinessesDespicablenesses}.
%   That is where $1$ comes from in~\eqref{eqn:edison_tesla_boxfight}.
% }
\begin{proof}
    \cref{thm:sufficient_decrease_when_crit_bounded} provides
    \begin{equation*}
    \scalk(\vxk) -
        \scalk(\vxSk)
    =
    \scalk(\vxk) - \scalk(\vxk + \vsk)
    \ge
    \frac{1}{2}
    \CCsd \CCepschi \Deltak.
    \end{equation*}
    The small constraint violation $\thetak \le \CCdeltan$ and $\kk\notin \restoreIndices$
    imply~\eqref{eqn:shoalinessesDespicablenesses_finCrit}.
    Because of $\Deltak \le 1$ and ${\mu} \in (0,1)$ the
    bound~\eqref{eqn:shoalinessesDespicablenesses_finCrit} simplifies to
    $\thetak \le \Deltak^{1+{\mu}}$, and it follows that
    \(
    \CCthetaub^{{\psi}} \Deltak^{{\psi}(1+{\mu})}
    \ge \thetak^{\psi}.
    \)
    The final inequality in~\eqref{eqn:edison_tesla_boxfight} then leads to
    \[
    \frac{1}{2}
    \CCsd \CCepschi \Deltak
    \ge
    {\kappa}_{\theta} \CCthetaub \Deltak^{{\psi}(1+\mu)},
    % \quad \Leftrightarrow \quad
    % \frac{1}{2}
    % \frac{
    %     \CCsd \CCepschi }{
    %     {\kappa}_{\theta} \CCthetaub^{{\psi}}
    % }
    % \ge
    % \Deltak^{{\psi}(1+\mu)-1}
    \]
    which then implies~\eqref{eqn:algo_model_dec}.
\end{proof}

If additionally the constraint violation is small enough, then the filter pair at the
trial point will
not be filter-dominated by the pair at $\vxk$:
\begin{mylemma}[{\cite[Lemma 3.8]{fletcher_global_2002}}]%
    \label{thm:small_radius_small_constraint_ensure_filter_acceptance}
 Suppose the algorithm is applied to~\eqref{eqn:mop_nonlin}, that
    \cref{%
    ass:covering,%
    ass:closed_and_bounded,%
    ass:normal_step_exists_bounded,%
    ass:fully_linear_models,%
    ass:lipschitz_gradients_true,%
    ass:equivalent_norms,%
    ass:hess_norm_bounded,%
    ass:sufficient_decrease%
    } hold and that the Criticality Routine does not run infinitely.
    Suppose further, that $\kk\notin\restoreIndices$,
    that~\eqref{eqn:porcellanizingTrichromat} holds again, that $\Deltak \le \delta_\rho$ as
    in~\eqref{eqn:yogurt_plant_cup} and
    that $\thetak \le \CCdeltan$
    %and hence~\eqref{eqn:shoalinessesDespicablenesses_finCrit}
    and that
    \begin{equation}
        \thetak
        \le
        \delta_{{\theta}}
        :=
        \min \left\{
        \frac{1}{\CCthetaub}
        {\left(
            \frac{
                \CCsuccess\CCsd\CCepschi
            }{
                2 \filterConst
            }
        \right)}^{2},
        \frac{1}{\CCthetaub^{\frac{1}{\mu}}}
        {\left(
            \frac{
                \CCsuccess\CCsd\CCepschi
            }{
                2 \filterConst
            }
        \right)}^{\frac{1 + \mu}{\mu}}
        \right\}.
        \label{eqn:fan_stapler_palm}
    \end{equation}
    Then
    \(
        \scalShort\left(\vxk\right)
        -
        \scalShort\left(\vxk + \vsk\right)
        \ge
        \filterConst \thetak.
    \)
\end{mylemma}
\begin{proof}
From~\cref{thm:bound_on_constraint_violation2} we deduce that 
\(
     \frac{
         \thetak
         }{
         \CCthetaub
     }
     \le
     \Deltak^{1+\mu}
\)
or 
\(
     \frac{
         \thetak
         }{
         \CCthetaub
     }
     \le \Deltak^2
\) and thus
\begin{equation}
\Deltak
    \ge
    \min\left\{
        {\left(
            \frac{
                \thetak
            }{
                \CCthetaub
            }
         \right)}^{\frac{1}{1+\mu}},
         {\left(
            \frac{
                \thetak
            }{
                \CCthetaub
            }
         \right)}^{\frac{1}{2}}
    \right\}.
    \label{eqn:fan_stapler_palm1}
\end{equation}
With~\cref{thm:success_when_crit_bounded}, it follows that
\begin{align*}
    \scalShort\left(\vxk \right)
        -
        \scalShort\left(\vxk + \vsk\right)
    &%
    \ge
    \CCsuccess \left(
        \scalk\left(\vxk \right)
        -
        \scalk\left(\vxk + \vsk\right)
    \right)
    %\\
    %&
    \stackrel{\mathclap{\eqref{eqn:foreversTalced}}}
    \ge
    \frac{1}{2}
    \CCsuccess
    \CCsd \CCepschi \Deltak
    \\
    &
    \stackrel{\mathclap{\eqref{eqn:fan_stapler_palm1}}}
    \ge
    \frac{1}{2}
    \CCsuccess
    \CCsd \CCepschi
    \min\left\{
        {\left(
            \frac{
                \thetak
            }{
                \CCthetaub
            }
         \right)}^{\frac{1}{1+\mu}},
         {\left(
            \frac{
                \thetak
            }{
                \CCthetaub
            }
         \right)}^{\frac{1}{2}}
    \right\}
    %\\
    %&
    \stackrel{\mathclap{\eqref{eqn:fan_stapler_palm}}}
    \ge
    \filterConst \thetak.
\end{align*}
\end{proof}

In the preceding lemmata, it has always been assumed that the iteration index does not belong to
the set of restoration indices $\restoreIndices$.
This is ensured if both the radius and the constraint violation are small enough:
\begin{mylemma}[{\cite[Lemma 3.9]{fletcher_global_2002}}]%
    \label{thm:non_restore_iteration}
    Suppose the algorithm is applied to~\eqref{eqn:mop_nonlin}, that
    \cref{%
    ass:covering,%
    ass:closed_and_bounded,%
    ass:normal_step_exists_bounded,%
    ass:fully_linear_models,%
    ass:lipschitz_gradients_true,%
    ass:equivalent_norms,%
    ass:hess_norm_bounded,%
    ass:sufficient_decrease%
    } hold and that the Criticality Routine does not run infinitely.
    Suppose that~\eqref{eqn:porcellanizingTrichromat} holds again.
    Suppose further that
    \begin{equation}%
      \label{eqn:GrabWoodenEvent}
        \Deltak \le
        \delta_{\restoreIndices}
        :=
        \min \Bigg\{
        \gammass \delta_\rho,
        {\left(
            \frac{1}{
                \CCmu
            }
        \right)}^{\frac{1}{\mu}},
%        \label{eqn:overopinionatedBeachwear}
        {\left(
            \frac{
                (1 - \filterConst)\gammass^2\CCDelta\CCmu
            }{
                \CCubn \CCthetaub
            }
        \right)}^{\frac{1}{1-\mu}},
      \Bigg\}
    \end{equation}
    and that
    \begin{equation}
        \thetak \le \min \left\{ \delta_{{\theta}},
        \CCdeltan
        %\CCcritTesttheta
        \right\}.
        \label{eqn:stoopeAquadromes}
    \end{equation}
    If $\kk > 0$, then $\kk \notin \restoreIndices$.
\end{mylemma}
\begin{proof}
  The proof works exactly as in~\cite{fletcher_global_2002}, because the
  Criticality Routine is not entered restoration iterations.
\end{proof}

The previous lemmata allow us to investigate two mutually exclusive cases defined
by the number of filter iterations.
The first convergence result is obtained for the case that there are
infinitely many such iterations.

\begin{mylemma}
  %[{\cite[Lemma 12]{ours}}]%
    \label{mythm:del_to_zero}
    Suppose the algorithm is applied to~\eqref{eqn:mop_nonlin}, that
    \cref{%
    ass:covering,%
    ass:closed_and_bounded,%
    ass:normal_step_exists_bounded,%
    ass:fully_linear_models,%
    ass:lipschitz_gradients_true,%
    ass:equivalent_norms,%
    ass:hess_norm_bounded,%
    ass:sufficient_decrease%
    } hold and that the Criticality Routine does not run infinitely.
	Suppose $\abs{\filterIndices} = \infty$.
    For any subsequence of iteration indices $\{\kk_{\ell}\}$ with
    $\abs{ \{ \kk_{\ell} \} \cap \filterIndices} = \infty$ it holds that
    \(
    \liminf_{{\ell} \to \infty} \Deltakl = 0.
    \)
    Especially, we have
    \(
     \lim_{\substack{\kk \to \infty, \\ \kk \in \filterIndices}}
      \Deltak = 0.
    \)
\end{mylemma}

%\todo[inline]{$\lim$ or $\liminf$ ???}

\begin{proof}
  The proof is along the lines of~\cite[Lemma
  3.10]{fletcher_global_2002},
  but additionally takes into account the Criticality Routine.
  First, let $\{\kk_{\ell}\}$ a sequence of indices containing
  infinitely many filter indices,
  $\abs{\{\kk_{\ell} \}\cap \filterIndices} = \infty$,
  and assume that it was bounded away from zero:
  \begin{equation}
    \Deltakl
    \ge
    \DeltaMin > 0
    \qquad \forall{\ell}\in \mathbb{N}_0.
    \label{eqn:outthrowGonadectomy}
  \end{equation}
  From
  \cref{mythm:infinite_filter_iters_imply_feasibility} and
  \cref{ass:normal_step_exists_bounded}
  it follows that
  \begin{equation}
    \lim_{{\ell}\to \infty} \normTrIt{ \vnkl }
    \stackrel{\text{\cref{ass:equivalent_norms}}}=
    \lim_{{\ell}\to \infty} \normTwo{ \vnkl }
    = 0.
    \label{eqn:siltsWeighted}
  \end{equation}
  Consequently, for sufficiently large ${\ell}$, the steps will become compatible and
  satisfy~\eqref{eqn:compatible_normal_step} and thus $\kk_{\ell}\notin\restoreIndices$.
  Hence, \cref{thm:normal_step_decrease_ub} applies and
  inequality~\eqref{eqn:deafenedPonytailed} holds again:
  \[
    \abs{
      \scalkl\left(\vxNkl \right)
      -
      \scalkl\left(\vxkl \right)
    }
    \le \CCubj \normTwo{{\vnkl}}
    + \frac{1}{2}\CChessF
    \normTwo{{\vnkl}}^2,
  \]
  so that~\eqref{eqn:siltsWeighted} gives
  \begin{equation}
    \lim_{{\ell}\to \infty}
    \abs{
      \scalkl\left( \vxkl \right)
      -
      \scalkl\left( \vxNkl \right)
    }
    = 0.
    \label{eqn:acceptedlyUnprofessed}
  \end{equation}
  % Step~\ref{algo:crit_test} and
  \cref{mythm:infinite_filter_iters_imply_feasibility}
  %ensures that, for large ${\ell}$, the Criticality Routine can affect the
  %inexact criticality value and the radius.
  enables the Criticality Routine.
  Additionally, from the boundedness of $\Deltakl$ as per~\eqref{eqn:outthrowGonadectomy}
  and from~\cref{ass:normal_step_exists_bounded},
  it also follows that, for large ${\ell}$, the Criticality Routine will not exit
  because no compatible normal step exists anymore.
  This ensures that
  \begin{equation}
    \chikl
    \ge \min \left\{ \CCcritTestchi, \frac{\Deltakl}{\CCMu} \right\}
    \ge \min \left\{ \CCcritTestchi, \frac{\DeltaMin}{\CCMu} \right\}
    =:
    \mathtt{z} > 0
    \label{eqn:limacineAgrestal}
  \end{equation}
  for large ${\ell}$.
  We can use this fact in
  \cref{ass:sufficient_decrease}
  % (together with $\CCMu \ge 1$ and $\CCcritDenom\ge 1$)
  and obtain
  \begin{equation}
    \scalkl \left(\vxNkl \right)
    -
    \scalkl\left(\vxSkl \right)
    \ge
    \underbrace{
      \CCsd
      % \min \left\{ \CCcritTestchi, \frac{\DeltaMin}{\CCMu} \right\}
      \mathtt z
      \min \left\{
        % \frac{\CCcritTestchi}{\CCcritDenom},
        % \frac{\DeltaMin}{\CCMu\CCcritDenom},
        \frac{\mathtt z}{\CCcritDenom},
        \DeltaMin
      \right\}
    }_{\text{const.}}
    > 0.
    \label{eqn:restrictivenessPassage}
  \end{equation}
  If we add zero
  we obtain the decomposition
  % ~\eqref{eqn:curculiosToutie}
  \begin{equation}
    %\begin{aligned}
      %&%
      \scalkl\left(\vxkl\right)
        -
        \scalkl\left(\vxSkl \right)
        =
      %\\
      %&\;
      \left(
        \scalkl\left(\vxkl\right)
        -
        \scalkl\left(\vxNkl \right)
      \right)
        + 
      \left(
        \scalkl\left(\vxNkl\right)
        -
        \scalkl\left(\vxSkl\right)
      \right).
    %\end{aligned}
    \label{eqn:curculiosToutie}
  \end{equation}
  By plugging in~\eqref{eqn:acceptedlyUnprofessed} for the left set of parentheses,
  we see that in the limit the difference of values is the same:
  \[
    \lim_{{\ell}\to \infty}
    % \left(
    \scalkl\left(\vxkl\right)
    -
    \scalkl\left(\vxSkl\right)
    % \right)
    =
    \lim_{{\ell}\to \infty}
    \scalkl\left(\vxNkl\right)
    -
    \scalkl\left(\vxSkl\right).
  \]
  Because of $\abs{\{\kk_{\ell}\}\cap\filterIndices} = \infty$,
  there is a subsequence $\{\kk_j\}\subseteq \{\kk_{\ell}\}$ of filter indices,
  $\kk_j\in \filterIndices$,
  for which it must hold that $\kk_j \in \restoreIndices$ or
  that~\eqref{eqn:algo_model_dec} fails.
  We have already shown that $\kk_j\notin\restoreIndices$ for large $j$.
  Hence, for $j$ sufficiently large, it must hold that
  %\begin{equation*}
  \(
		{\kappa}_{\theta} {\theta}\itthetaCmd{\kk_j}^{\psi}
    \ge
		\scalShort\itCmd{\kk_j}(\vx\itCmd{\kk_j})
		-
    \scalShort\itCmd{\kk_j}(\vx\itCmd{\kk_j} + \vs\itCmd{\kk_j})
  %\end{equation*}
  \)
  and~\cref{mythm:infinite_filter_iters_imply_feasibility} 
  implies that both sides must go to zero.
  % \[
  %   \lim_{j\to \infty}
  %   % \left(
  %   \scalShort\itCmd{\kk_j}(\vx\itCmd{\kk_j})
	% 	-
  %   \scalShort\itCmd{\kk_j}(\vx\itCmd{\kk_j} + \vs\itCmd{\kk_j})
  %   % \right)
  %   = 0
  %   .
  % \]
  This is a contradiction to~\eqref{eqn:restrictivenessPassage}.
  Hence, no infinite subsequence $\{\kk_{\ell}\}$
  with $\abs{\{\kk_{\ell}\}\cap \filterIndices} = \infty$ and~\eqref{eqn:outthrowGonadectomy} can
  exist.
\end{proof}

\begin{myremark}%
    \label{remark:durnedestStonen}
    Without further assumptions on the restoration step it does not seem possible to
    show $\lim_{\kk\to \infty} \Deltak = 0$ in the case that $\abs{\filterIndices} = \infty$.
    For any element of an index sequence $\{\kk_{\ell}\}$,
    bounded away from zero and with
    $\abs{\{\kk_{\ell}\}\cap\filterIndices} < \infty$,
    there might be a preceding
    restoration iteration with arbitrarily small radius itself, but possibly increasing
    it without any restriction,
    so that no contradiction can be
    derived.
    This difficulty has been observed in the literature before, see
    e.g.~\cite[Remark 7]{wachter_line_2005}.
\end{myremark}

\begin{mylemma}%[{\cite[Lemma 11]{ours}}]%
	\label{mythm:del_bounded}
    Suppose the algorithm is applied to~\eqref{eqn:mop_nonlin}, that
    \cref{%
    ass:covering,%
    ass:closed_and_bounded,%
    ass:normal_step_exists_bounded,%
    ass:fully_linear_models,%
    ass:lipschitz_gradients_true,%
    ass:equivalent_norms,%
    ass:hess_norm_bounded,%
    ass:sufficient_decrease%
    } hold and that the Criticality Routine does not run infinitely.
	Suppose that $\abs{\filterIndices} = \infty$ and
    that the doubly inexact criticality is bounded as
	per~\eqref{eqn:porcellanizingTrichromat} for all $\kk \in \filterIndices$.
    Then the trust-region radius is also bounded below, i.e., there is some
    $\DeltaMin > 0$ such that
    \(
    \Deltak \ge \DeltaMin 
    %\qquad\text{for all $k\in \filterIndices$.}
    \)
    for all $k\in \filterIndices$.
\end{mylemma}

\begin{proof}
    %In contrast to trust-region methods without a filter, the trust-region radius
    %may be changed not only by the usual updates but also by restoration steps
    %(which have indices in $\restoreIndices\subseteq\filterIndices$).
    If there are infinitely many filter indices $\filterIndices$,
    then we know from~\cref{mythm:infinite_filter_iters_imply_feasibility} that
    $\thetak \to 0$.
    %By~\cref{ass:normal_step_exists_bounded} it follows that for $\kk$ large enough
    %the normal step $\vnk$ exists and converges to zero (with~\cref{ass:equivalent_norms}).
    Suppose ${\{\Deltak\}}_{\kk \in \filterIndices}$ is not bounded away from zero.
    Then there is a subsequence $\{\kk_{\ell}\}\subseteq \filterIndices$ with
    $\lim_{{\ell}\to \infty} \Deltakl = 0$.
    Following the argumentation in~\cite[Lemma 3.10]{fletcher_global_2002}, we see that
    %there is a subsequence of indices $\{\kk_{\ell}\}\subseteq \filterIndices$,
    %such that
    for ${\ell}$ large enough,
    \cref{thm:non_restore_iteration} applies and
    guarantees that $\kk_{\ell}\notin\restoreIndices$.
    At the same time, \cref{thm:model_dec_test_succeeds} applies for
    large ${\ell}$,
    so that the test~\eqref{eqn:algo_model_dec} always succeeds for $\vxkl$.
    Hence, there is some ${\ell}_0\in \mathbb{N}_0$ such that for all
    $\kk_{\ell}$ with ${\ell}\ge {\ell}_0$,
    it holds that $\kk_{\ell} \notin \filterIndices$.
    This contradicts $\{\kk_{\ell}\}\subseteq \filterIndices$.
\end{proof}

We are finally able to state the convergence result for the case of infinitely many filter 
iterations as a cumulative corollary derived from 
~\cref{thm:infinite_filter_iters_imply_feasibility,mythm:del_to_zero,%
mythm:del_bounded}:
\begin{mycorollary}[{\cite[Lemma 3.10]{fletcher_global_2002}}]%
    \label{thm:infinite_filter_iterations_convergence}
    Suppose the algorithm is applied to~\eqref{eqn:mop_nonlin}, that
    \cref{%
    ass:covering,%
    ass:closed_and_bounded,%
    ass:normal_step_exists_bounded,%
    ass:fully_linear_models,%
    ass:lipschitz_gradients_true,%
    ass:equivalent_norms,%
    ass:hess_norm_bounded,%
    ass:sufficient_decrease%
    } hold and that the Criticality Routine does not run infinitely.
    Suppose further, that $\abs{\filterIndices} = \infty$.
	Then there is a subsequence $\{\kk_{\ell}\} \subseteq \filterIndices$
    of filter indices with
    \begin{equation*}
      \label{eqn:backpack_mug_phone1}
      \lim_{{\ell}\to \infty} {\theta}\itthetaCmd{\kk_{\ell}} = 0,
      \quad
      \lim_{\ell \to \infty} \Deltakl = 0 
    \quad \text{and} \quad
    \lim_{{\ell}\to \infty} {\chi}\itCmd{\kk_{\ell}} = 0.
    \end{equation*}
    %It also holds that $\lim_{{\ell}\to \infty} \Deltakl = 0$.
\end{mycorollary}

% \begin{proof}
%   \cref{mythm:del_to_zero} provides $\lim_{{\ell}\to \infty} \Deltakl = 0$.
%   The first part of \eqref{eqn:backpack_mug_phone1} follows from
%   \cref{thm:infinite_filter_iters_imply_feasibility}.

% 	To derive a contradiction for the second part, assume that the
% 	inexact criticality is bounded as per~\eqref{eqn:porcellanizingTrichromat},
% 	i.e.,
%     \begin{equation}
%         {\chi}\itCmd{\kk}
%         \ge
%         \frac{1}{\CCnorm}
%         {\chi}_2\itCmd{\kk}
%         \ge
%         \frac{1}{\CCnorm}
%         \CCepschi > 0 \qquad \forall\kk\in\filterIndices.
%         \label{eqn:backpack_mug_phone3}
%     \end{equation}
% 	From~\cref{mythm:del_bounded} we know that there then is $\DeltaMin > 0$ with
% 	\(
%     \Deltak \ge \DeltaMin > 0
%   \)
%   for all $\kk \in \filterIndices$.
%   But this contradicts $\lim_{{\ell}\to \infty} \Deltakl = 0$
%   and the bound in~\eqref{eqn:backpack_mug_phone3} cannot hold,
%   implying the existence of a suitable subsequence
%   $\{\kk_{\ell}\}\subseteq \filterIndices$.
% \end{proof}

\subsubsection{Finitely Many Filter Iterations}
We have shown that a quasi-stationary point is approached when there are
infinitely many filter indices.
We now concentrate on the case that there are only finitely many
filter indices (and note that this then implies that there are only
finitely many restoration indices as well).
\textbf{From now on, $\kk_0$ is the last iteration index for which
$\vx\itCmd{\kk_0 -1}$ is added to the filter.}

\begin{mylemma}[{\cite[Lemma 3.11]{fletcher_global_2002}}]%
    \label{thm:finite_filter_iter_implies_feasibility}
    Suppose the algorithm is applied to~\eqref{eqn:mop_nonlin}, that
    \cref{%
    ass:covering,%
    ass:closed_and_bounded,%
    ass:normal_step_exists_bounded,%
    ass:fully_linear_models,%
    ass:lipschitz_gradients_true,%
    ass:equivalent_norms,%
    ass:hess_norm_bounded,%
    ass:sufficient_decrease%
    } hold and that the Criticality Routine does not run infinitely.
    Suppose further, that $\abs{\filterIndices} < \infty$.
    Then
    \begin{equation}
    \lim_{\kk \to \infty} \thetak = 0.
    \label{eqn:delicious_gummy_bears0}
    \end{equation}
    Furthermore, $\vnk$
    satisfies~\eqref{eqn:balladedAdzukis}
    %\eqref{eqn:bound_on_constraint_violation1}
    for all $\kk\ge \kk_0$ large enough.
\end{mylemma}

\begin{proof}
    First consider the case that there are only finitely many successful
    indices $\acceptIndices$.
    Then $\Deltak \xrightarrow{\kk\to \infty} 0$. Because there are no
    restoration iterations for $\kk \ge \kk_0$, it follows from
    $\eqref{eqn:compatible_normal_step}$ that $\vnk \to 0$.
    Eq.~\eqref{eqn:bound_on_constraint_violation1_finCrit} finally gives
    $\thetak \to 0$.

    Now, suppose that $\abs{\acceptIndices} = \infty$ and
    consider any \emph{successful} iteration with $\kk\ge\kk_0$.
    Then $\vxk$ is not added to the filter, and it follows from step~\ref{algo:acceptance_test} of the algorithm
    that $\rho\KK \ge \CCaccept$ and from the definition of $\rho\KK$ and from~\eqref{eqn:algo_model_dec}
    that
    \begin{equation}
    %\begin{aligned}
    \scalShort\left(\vxk\right) -
    \scalShort\left(\vxk + \vsk\right)
    %&%
    \ge
    \CCaccept
    \left(
        \scalk\left(\vxk\right) -
        \scalk\left(\vxk + \vsk\right)
    \right)
    %\\
    %&%
    \ge
    \CCaccept {\kappa}_{\theta} \thetak^{\psi} \ge 0.
    \label{eqn:delicious_gummy_bears}
    %\end{aligned}
    \end{equation}
    The sequence ${\left\{\scalShort\left(\vxk\right)\right\}}_{\kk \ge \kk_0}$
    is bounded below due to~\cref{ass:closed_and_bounded,ass:lipschitz_gradients_true},
    and it is monotonically decreasing because~\eqref{eqn:algo_model_dec} always succeeds
    due to $\kk \ge \kk_0$.
    Thus, no
    new iterate is ever chosen with $\rho\KK < \CCaccept$.
    Hence, for the \ac{rhs} in~\eqref{eqn:delicious_gummy_bears} it follows that
    \begin{equation}
    \lim_{\substack{\kk\in\acceptIndices, \\ \kk\to \infty}}
    \scalShort\left(\vxk\right) -
    \scalShort\left(\vxk + \vsk\right)
    = 0.
    \label{eqn:delicious_gummy_bears1}
    \end{equation}
    The limit~\eqref{eqn:delicious_gummy_bears0} follows from the \ac{rhs} in~\eqref{eqn:delicious_gummy_bears}
    by noticing that ${\theta}\iterSubCmd{j} = \thetak$ for all non-successful
    indices $j\ge k$ with $\rho\itCmd{j}<\CCaccept$. 
    The bound~\eqref{eqn:balladedAdzukis} holds eventually because
    $\thetak \le \CCdeltan$ for large $\kk$
    and then~\cref{thm:bound_on_constraint_violation1} applies.
\end{proof}

The next auxiliary lemma shows that the trust-region radius is bounded below if the asymptotically feasible
iterates do not approach a quasi-stationary point. It is used afterwards to derive a contradiction.
\begin{mylemma}[{\cite[Lemma 3.12]{fletcher_global_2002}}]%
    \label{thm:finite_filter_iterations_bounded_radius}
    Suppose the algorithm is applied to~\eqref{eqn:mop_nonlin}, that
    \cref{%
    ass:covering,%
    ass:closed_and_bounded,%
    ass:normal_step_exists_bounded,%
    ass:fully_linear_models,%
    ass:lipschitz_gradients_true,%
    ass:equivalent_norms,%
    ass:hess_norm_bounded,%
    ass:sufficient_decrease%
    } hold and that the Criticality Routine does not run infinitely.
    Suppose further, that $\abs{\filterIndices} < \infty$ and that
    the doubly inexact criticality is bounded away from zero as in~\eqref{eqn:porcellanizingTrichromat} for all
    $\kk\ge\kk_0$.
    Then there is a $\DeltaMin > 0$ such that
    \(
    \Deltak \ge \DeltaMin 
    %\qquad \forall k \in \mathbb{N}_0.
    \)
    for all $k \in \mathbb{N}_0$.
\end{mylemma}
\begin{proof}
    %As in the proof of~\cref{thm:non_restore_iteration} we note that~\eqref{eqn:crit_test_final_radius}
    We sketch, how to adapt the proof in~\cite{fletcher_global_2002}:
    \cref{thm:finite_filter_iter_implies_feasibility} ensures that for large enough $\kk$
    we have $\thetak \le \CCcritTesttheta$ and the Criticality Routine affects the criticality
    value and the radius.
    For these $\kk$ it follows from~\eqref{eqn:crit_test_final_radius}
    and~\eqref{eqn:porcellanizingTrichromat} that
    \(
      \Deltak 
      \ge 
      \min\left\{ \CCBeta \omegak, \DeltaPrek \right\}
      \ge 
      \min\left\{ \CCBeta \CCepschi, \DeltaPrek \right\}
    \).
    Further, if $\kk \ge 1$ and $\kk -1\notin\restoreIndices$, then 
    \(
      \Deltak 
      \ge 
      \min\left\{ \CCBeta, \CCepschi, \gammass \Delta\itDeltaCmd{\kk-1}\right\}.
    \)
    % \begin{align}
    % \Deltak 
    % &% 
    % \ge
    % \min\left\{ \CCBeta \omegak, \DeltaPrek \right\}
    % \ge
    % \min\left\{ \CCBeta \CCepschi, \DeltaPrek \right\}
    % %\label{eqn:polydactyliesNondormant}
    % \nonumber
    % \shortintertext{%
    % and, if $\kk \ge 1$ with $\kk - 1\notin \restoreIndices$, that%
    % }
    % \Deltak &\ge
    % \min\left\{ \CCBeta \CCepschi, \gammass \Delta\itDeltaCmd{\kk - 1} \right\}.
    % \label{eqn:embrutingLove}
    % \end{align}
    %\todo[inline]{\eqref{eqn:polydactyliesNondormant} and~\eqref{eqn:embrutingLove}
    %only hold for large $\kk$.
    %}
    %\cref{thm:finite_filter_iter_implies_feasibility} applies, so suppose
    Suppose
    $\kk_1\ge\kk_0$ is large enough that $\thetak \le \CCcritTesttheta$
    and~\eqref{eqn:stoopeAquadromes} is fulfilled for all
    $\kk\ge \kk_1$, i.e., $\thetak \le \min\{ \CCcritTesttheta, {\delta}_{{\theta}} \}$
    for all $\kk\ge \kk_1$.
    For the purpose of deriving a contradiction, one may now assume
    that $j \ge \kk_1$ is the first index with
    \begin{equation}
        \Delta\itDeltaCmd{j}
        \le
        \gammass \min \left\{
        \delta_\rho,
        \sqrt{
            \frac{
                (1- \filterConst ){\theta}^{\filter}
            }{
                \CCthetaub
            }
        },
        \Delta\itDeltaCmd{\kk_1},
        %\frac{
            \CCBeta \CCepschi
        %}{\gammag}
        \right\}
        =:
        \gammass {\delta}_s,
        \label{eqn:nicenessesThreadiest}
    \end{equation}
    where ${\delta}_\rho$ is as defined in~\cref{thm:success_when_crit_bounded} and
    \(
    {\theta}^{\filter}
    :=
    \min_{i\in \filterIndices} {\theta}\iterSubCmd{i}
    \)
    and proceed as in~\cite{fletcher_global_2002}.
\end{proof}

\begin{mylemma}[{\cite[Lemma 3.13]{fletcher_global_2002}}]%
    \label{mythm:finite_filter_iterations_radius_to_zero}
    Suppose the algorithm is applied to~\eqref{eqn:mop_nonlin}, that
    \cref{%
    ass:covering,%
    ass:closed_and_bounded,%
    ass:normal_step_exists_bounded,%
    ass:fully_linear_models,%
    ass:lipschitz_gradients_true,%
    ass:equivalent_norms,%
    ass:hess_norm_bounded,%
    ass:sufficient_decrease%
    } hold and that the Criticality Routine does not run infinitely.
    Suppose further, that $\abs{\filterIndices} < \infty$.
    Then
    \(
    \liminf_{\kk \to \infty} \Deltak = 0.
    \)
\end{mylemma}

\begin{proof}
    If there are only finitely many successful indices, then the result follows
    immediately from the radius update rules.

    Thus, suppose that $\abs{\acceptIndices} = \infty$ and for a contradiction assume
    that there is a constant $\DeltaMin > 0$ such that
    \begin{equation}
        \Deltak \ge \DeltaMin > 0 \qquad \forall\kk\in \mathbb{N}_0.
        \label{eqn:jejunityThereanent}
    \end{equation}
    As before, we see that from~\cref{thm:finite_filter_iter_implies_feasibility}
    it follows for large $\kk$ that the normal step $\vnk$ exists and satisfies~\eqref{eqn:normal_step_exists_bounded} in~\cref{ass:normal_step_exists_bounded}.
    Then the same equations as in~\cref{thm:finite_filter_iter_implies_feasibility} hold,
    namely~\eqref{eqn:delicious_gummy_bears} and~\eqref{eqn:delicious_gummy_bears1}.
    %For the purpose of obtaining a contradicting assume that $\chik$ is bounded as in
    %\eqref{eqn:backpack_mug_phone3}.
    Very much like in the proof of~\cref{thm:infinite_filter_iterations_convergence},
    we can again decompose the model decrease via~\eqref{eqn:curculiosToutie}.
    For the first term, \eqref{eqn:deafenedPonytailed} applies again,
    % :
    % \[
    % \abs{
    % \scalk\left(\vxk\right)
    %    -
    %    \scalk\left(\vxNk \right)
    % }
    % \le
    %     \CCubj \normTwo{{\vnk}}
    %     + \frac{1}{2}\CChessF \normTwo{{\vnk}}^2
    % \]
    yielding
    \(
    \lim_{\kk \to \infty}
    \left(
    \scalk\left(\vxk\right)
       -
    \scalk\left(\vxNk\right)
    \right)
    = 0.
    \)
    From equations~\eqref{eqn:delicious_gummy_bears},~\eqref{eqn:delicious_gummy_bears1}
    %in the proof of~\cref{thm:finite_filter_iter_implies_feasibility}
    and from the model decrease
    decomposition it follows that
    \begin{equation}
    \lim_{
        \substack{
        \kk \to \infty, \\
        \kk \in \acceptIndices
        }
    }
        \left(
            \scalk\left(\vxNk\right)
            -
            \scalk\left(\vxSk \right)
        \right)
        =0.
        \label{eqn:slowingsRecreance}
    \end{equation}
    But the sufficient decrease condition (\cref{ass:sufficient_decrease}) still holds and
    (for large $\kk$ and assuming~\eqref{eqn:jejunityThereanent} holds)
    the Criticality Routine ensures
    \[
    \chik
    \ge
    \min\left\{ \CCcritTestchi, \frac{\Deltak}{\CCMu}\right\}
    \ge
    \min\left\{ \CCcritTestchi, \frac{\DeltaMin}{\CCMu}\right\}
    =: \mathtt{z} > 0,
    \]
    which can be plugged in and gives:
    \begin{align*}
    \scalk\left(\vxNk\right)
       -
    \scalk\left(\vxSk\right)
    \ge
    \CCsd
        %\min
        \mathtt{z}
        \min \left\{
        \frac{\mathtt{z}}{\CCcritDenom}
        ,
        \DeltaMin
        \right\} > 0,
    \end{align*}
    where the \ac{rhs} is constant, contradicting~\eqref{eqn:slowingsRecreance}.
\end{proof}

\begin{mycorollary}[{\cite[Lemma 3.13]{fletcher_global_2002}}]%
    \label{thm:infinite_filter_iterations_crit_min_convergence}
     Suppose the algorithm is applied to~\eqref{eqn:mop_nonlin}, that
    \cref{%
    ass:covering,%
    ass:closed_and_bounded,%
    ass:normal_step_exists_bounded,%
    ass:fully_linear_models,%
    ass:lipschitz_gradients_true,%
    ass:equivalent_norms,%
    ass:hess_norm_bounded,%
    ass:sufficient_decrease%
    } hold and that the Criticality Routine does not run infinitely.
    Suppose further, that $\abs{\filterIndices} < \infty$.
    Then
    \[
    \liminf_{\kk\to \infty} \chik = 0.
    \]
    For a subsequence $\{\kk_{\ell}\} \subseteq \{\kk\}$ with $\chikl \to 0$ it must also hold
    that $\lim_{{\ell}\to \infty} \Deltakl = 0$.
\end{mycorollary}

\begin{proof}
    %The proof works similar to that of~\cref{thm:infinite_filter_iterations_convergence}:
    To derive a contradiction, it is assumed that $\chik$ is bounded away
    from zero
    as per~\eqref{eqn:porcellanizingTrichromat}
    for all $\kk \in \mathbb{N}_0$.
    Then the radius must be bounded away from zero due to
    \cref{thm:finite_filter_iterations_bounded_radius},
    in contradiction to $\liminf_{\kk\to \infty} \Deltak= 0$, as assured by
    \cref{mythm:finite_filter_iterations_radius_to_zero}.
    If $\{\kk_{\ell}\}$ is such that $\chikl \to 0$, but we assume
    $\lim_{{\ell}\to \infty} \Deltakl = \DeltaMin > 0$, then the Criticality Routine again
    ensures~\eqref{eqn:limacineAgrestal} for large ${\ell}$, which
    is a contradiction.
\end{proof}

%We can now conclude the analysis with the main result of this section:
\begin{mytheorem}[Convergence to Quasi-Stationary Point]
    \label{thm:convergence1}
    Suppose the algorithm is applied to~\eqref{eqn:mop_nonlin} and that
    \cref{%
    ass:covering,%
    ass:closed_and_bounded,%
    ass:normal_step_exists_bounded,%
    ass:fully_linear_models,%
    ass:lipschitz_gradients_true,%
    ass:equivalent_norms,%
    ass:hess_norm_bounded,%
    ass:sufficient_decrease%
    } hold.
    Let $\{\vxk\}$ be the sequence of iterates produced by the algorithm.
    Then either the restoration procedure in step~\ref{algo:restoration} terminates unsuccessfully
    by converging to an infeasible, first-order critical point of
    problem~\eqref{eqn:restoration_problem} or there is a subsequence $\{\kk_{\ell}\}$ of indices
    for which
    \(
      \lim_{{\ell}\to \infty} \vxkl = \bar{\vx}
    \)
    and $\bar{\vx}$ is quasi-stationary, i.e.,
    \[
      {\theta}\left( \bar{\vx} \right) 
      = 
      \lim_{{\ell}\to \infty} \thetakl
      =
      0
      \text{  and  }
      \lim_{{\ell}\to \infty} \chikl = 0,
    \]
    and for which it also holds that $\lim_{{\ell}\to \infty} \Deltakl = 0$.
\end{mytheorem}
\begin{proof}
If the restoration procedure always terminates successfully, then the
convergence to a quasi-stationary point follows
from~\cref{ass:closed_and_bounded} and
\cref{thm:infinite_filter_iterations_convergence,%
thm:infinite_crit_routine,%
thm:infinite_filter_iters_imply_feasibility,%
thm:infinite_filter_iterations_crit_min_convergence%
}.
\end{proof}

\section{Convergence to KKT-Points}%
\label{section:kkt_convergence}

In this section, we conclude that a quasi-stationary point is also a KKT-point
under suitable constraint qualifications.
In the first step, we show that we can easily replace the 
\emph{objective} surrogates
with their true counterparts.
\begin{mylemma}%
	\label{thm:true_objective_criticality}
	Suppose \cref{%
    ass:covering,%
    ass:lipschitz_gradients_true,%
    ass:fully_linear_models,%
    ass:normal_step_exists_bounded,%
    ass:equivalent_norms}
  hold and suppose that it holds for a subsequence 
  $\{\vxk\}$ of iterates that
	\[
	\lim_{\kk \to \infty}\Deltak = 0
	\quad\text{and}\quad
	\lim_{\kk \to \infty} \omegak_2 =
	\lim_{\kk \to \infty} {\omega}\left(  \vxNk; \mfk, \linFeasApproxk, \normTwo{{\bullet}} \right) = 0.
	\]
	Then it also holds for the true objectives $\vf$ that
	\[
	\lim_{\kk \to \infty} {\omega}\KK_2 =
	\lim_{\kk \to \infty} {\omega}\left(  \vxNk; \vf, \linFeasApproxk, \normTwo{{\bullet}} \right) = 0.
	\]
\end{mylemma}

\begin{proof}
	Similar to~\cite[Lemma 7]{ours},
    we can show that for
	$\kk\in \mathbb{N}_0$ it must hold that
	\[
	\abs{
		{\omega}\left(  \vxNk; \mfk, \linFeasApproxk, \normTwo{\bullet}  \right)
		-
		{\omega}\left(  \vxNk; \vf, \linFeasApproxk, \normTwo{\bullet}  \right)
	}
	=
	\abs{
		\omegak_2 - {\omega}\KK_2
	}
	\le \mathtt{s} \Deltak
	\]
	for some constant $\mathtt{s} > 0$ (independent of $\kk$).
	The triangle inequality yields
	\[
	\abs{{\omega}\KK_2}
	\le
	\abs{{\omega}\KK_2 - \omegak_2} + \abs{\omegak_2} \le
	\mathtt{s}\Deltak + \abs{\omegak_2}
	\]
  and the \ac{rhs} goes to zero for $k\to \infty$.
\end{proof}

Note that \cref{thm:true_objective_criticality} applies both to a sequence of the main 
algorithm and to an infinite subsequence of the Criticality Routine 
(see \cref{thm:infinite_crit_routine}).
In our second step, we now replace the \emph{constraint} surrogates with the original functions.

\begin{mylemma}
Suppose \cref{%
    ass:covering,%
    ass:lipschitz_gradients_true,%
    ass:fully_linear_models,%
%    ass:equivalent_norms%
} hold.
Let $\{\vxk\}\subseteq \mathbb{R}^\varDim$ be an algorithmic sequence with
$\vxk\to \bar{\vx} \in \feas$,
${\omega}\left( \vxNk; \vf, \linFeasApproxk, \normInf{{\bullet}} \right) \to 0$
%,$\thetak \to 0$ 
and $\Deltak \to 0$.
Further, assume that the \ac{mfcq} hold at $\bar{\vx}$.
 % and that the \ac{mfcq} hold at $\vx$,
i.e.,
the rows of $\DH(\bar{\vx})$ are linearly independent,
and there is a direction $\vd \in \mathbb{R}^\varDim$
such that $\DH(\bar{\vx})\vd = \vZ$ and
$\tran{\vd} \dg_{\ell}(\bar{\vx}) < 0$ for all ${\ell}$ with
$\gL_\ell(\bar{\vx}) = 0$.
Then
$\bar{\vx}$ is a KKT-point of~\eqref{eqn:mop_nonlin}.
\end{mylemma}

\begin{proof}
% In the following, these notations are used:
% \begin{enumerate}
%     \item $\DmH$ is the Jacobian of the surrogate for $\ve h$,
%     evaluated at $\ve x\KK$.
%     \item $\DmG$ is the Jacobian of the surrogate for $\ve g$,
%     evaluated at $\ve x\KK$.
%     \item $\mgkSub$ is $\mgk(\vxk) + \DmG \cdot \vn\KK \le \vZ$.
% \end{enumerate}
In the following, we use $\DmH$ and $\DmG$ to denote the Jacobians of 
the surrogates $\mhk$ and $\mgk$, evaluated at $\vxk$.
Likewise, $\DFk$ is the Jacobian of the true objective function at 
$\vxk$, while $\mgkSub \le \vZ_\gDim$ is defined as $\mgk(\vxk) + \DmG \cdot \vn\KK$.\\
\Cref{ass:normal_step_exists_bounded} gives $\vnk \to \vZ$.
Because of this, and the error bounds in \cref{ass:fully_linear_models},
we have that
$\vxk \to \bar{\vx}$, 
$\Deltak \to 0$ as well as 
$\DmH \to \DH = \jac \vh (\bar{\vx}),
\DmG\to \DG = \jac \vh (\bar{\vx})$ and $\mgkSub\to \vg = \vg(\bar{\vx})$.\\
By assumption,
${\omega}\left( \vxNk; \vf, \linFeasApproxk, \normInf{{\bullet}} \right) \to 0$.
We thus have a sequence of linear programs,
\begin{equation*}
  {\omega}\left( \vxNk; \vf, \linFeasApproxk, \normInf{{\bullet}} \right)
  =
  %\begin{aligned}[t]
  \max_{\ve d \in \mathbb R^n, {\beta}^- \in \mathbb R}
  \begin{bmatrix}
      \tran{\ve 0}_{\varDim} & 1
  \end{bmatrix}
  \begin{bmatrix}
      \ve d \\ {\beta}^-
  \end{bmatrix}
  \quad\text{s.t.}\quad
  %&
  %\ve d \in \mathbb{R}^\varDim,\;
  %{\beta}^- \ge 0, \\
  %&
  \begin{bmatrix}
      -\eye{\varDim} & \ve 0_\varDim \\
      \eye{\varDim} & \ve 0_\varDim \\
      \DFk & \ve 1_\fDim \\
      \DmH & \ve 0_\hDim \\
      \DmG & \ve 0_\gDim
      %\DmG & -\DmG & \ve 0_{\gDim} \\
  \end{bmatrix}%
  \begin{bmatrix}
          \ve d \\ {\beta}^-
  \end{bmatrix}
  \begin{matrix}
  \le \\ \le \\ \le \\ = \\ \le
  \end{matrix}
  \begin{bmatrix}
          \ve 1_\varDim\\
          \ve 1_\varDim\\
          \ve 0_\fDim \\
          \ve 0_\hDim\\
          -\mgkSub \\
          %-\mgkSub
  \end{bmatrix}%
  ,
  %\end{aligned}
\end{equation*}
the values of which go to zero.
%{\color{lightgray}
The dual problems are
\begin{equation}
\begin{aligned}
&\min_{%
  \ve y^1, \ve y^2, \ve y^3, \ve y^5 \ge \vZ , \ve y^4 \in \mathbb R^\hDim%
  }
    \begin{bmatrix}
        \tran{\ve 1}_{\varDim} & \tran{\ve 1}_{\varDim} &  \tran{\ve 0}_{\fDim} & \tran{\ve 0}_{\hDim} & - \tran{\mgkSub}
    \end{bmatrix}
    \begin{bmatrix}
        \ve y^1 \\ \vdots \\ \ve y^5
    \end{bmatrix}
    \quad \text{s.t.}\quad
    %\\
    % &
    % \ve y^1 \ge \vZ,
    % \ve y^2 \ge \vZ,
    % \ve y^3 \ge \vZ,
    % \ve y^4 \in \mathbb{R}^\hDim,
    % \ve y^5 \ge 0, \\
    % 
    %&
    \begin{bmatrix}
        -\eye{\varDim} & \eye{\varDim} & \tran{\DFk} & \tran{\DmH} & \tran{\DmG}\\
        \tran{\ve 0}_{\varDim} & \tran{\ve 0}_{\varDim} & \tran{\ve 1}_{\fDim} & \tran{\ve 0}_{\hDim} & \tran{\ve 0}_{\gDim}
    \end{bmatrix}%
    \begin{bmatrix}
        \ve y^1 \\ \vdots \\ \ve y^5
    \end{bmatrix}
    =
    \begin{bmatrix}
        \ve 0_{\varDim} \\ 1
    \end{bmatrix}
    .
\end{aligned}
\label{eqn:dual_k_unscaled}
\tag{D\ensuremath{_k}}
\end{equation}
%For each of the primal problems, zero is a feasible solution.
%Because $\left\{ {\omega}\left( \vxk; \vf, \linFeasApproxk  \right) \right\}_{\kk}$ is bounded,
By strong duality, the dual problem for $k\in \mathbb{N}_0$
is also always feasible, 
and its optimal value equals the primal optimal value.
Suppose that $\left\{  (\ve y_k^1, \ldots , \ve y_k^5) \right\}$
is a sequence of dual optimizers.
By strong duality, it follows from
${\omega}\left( \vxNk; \vf, \linFeasApproxk, \normInf{{\bullet}} \right) \to 0$
%and $ - \mgkSub \to -\vg( \bar{\vx} ) \ge 0$ 
that
\begin{equation*}
  \lim_{k\to \infty} \ve y_k^1 = \vZ,
  \quad
  \lim_{k\to \infty} \ve y_k^2 = \vZ,
  \quad\text{and}\quad
  \lim_{k\to \infty} - \tran{\mgkSub} \ve y_k^5 = 0.
\end{equation*}
First, we show that it must also hold for $i\in \{3,4,5\}$ that the sequences
$\{ \ve y^i\}$
are bounded (see also~\cite{echebest_inexact_2017} for a similar idea).
To derive a contradiction,
assume that for some $i \in \{ 3, 4, 5\}$ the sequence $\{ \ve y_k^i\}$ 
is not bounded.
Define
\(
  {\nu}_k = \max \left\{ 1, \max_{i=3,4,5} \normInf{\ve y_k^i} \right\},
\)
and we get a bounded sequence of scaled variables 
$$
\tilde{y}_k^i = \frac{1}{\nu_k} y_k^i \qquad i = 1,\ldots,5.
$$
% \begin{equation}
%     \tilde{\ve y}_k^i = 
%     \begin{cases}
%       \ve y_k^i &\text{for $i\in \{1,2\}$},\\ 
%       \frac{1}{{\nu}_k} \ve y_k^i &\text{for $i\in \{3,4,5\}$}.
%     \end{cases}
%     \label{eqn:mezereonInspections}
% \end{equation}
Then, ${\nu}_k \xrightarrow{k\to \infty} \infty$,
but we can take a subsequence $\mathcal K$ of indices such
that for every $i\in\{3,4,5\}$
the sequence ${\{\tilde{\ve y}_k^i\}}_{k\in \mathcal K}$ converges to $\bar{\ve y}^i$
and so that for one $i \in \{3,4,5\}$ the limit $\bar{\ve y}^i$ is not zero.
%\todo[inline]{proof that one limit is not zero?}
Consequently, we have 
\begin{equation}
  \lim_{\substack{k\to\infty,\\ k\in \mathcal K}}
    (\tilde{\ve y}_k^1,\ldots, \tilde{\ve y}_k^5)
    = 
      (\bar{\ve y}^1, \ldots, \bar{\ve y}^5) 
    =
      (\vZ_{\varDim}, \vZ_{\varDim}, \bar{\ve y}^3, \bar{\ve y}^4, \bar{\ve y}^5),
  %\; 
  %\lim_{\substack{k\to\infty,\\ k\in \mathcal K}} \tilde{\ve y}_k^1 = \vZ,
  %\;
  %\lim_{\substack{k\to\infty,\\ k\in \mathcal K}} \tilde{\ve y}_k^2 = \vZ,
  \quad\text{and}\quad
  \lim_{\substack{k\to\infty,\\ k\in \mathcal K}} - \tran{\mgkSub} \tilde{\ve y}_k^5 =
  -\tran{\ve g} \ve{\bar{y}}^5 =
  0
  .
  \label{eqn:week_afford_income}
\end{equation}
% For $i = 5 $, it still holds that
% \begin{equation}
%     %\lim_{\substack{k\to \infty,\\ k \in \mathcal K}} \tilde{\ve y}_k^1 = \vZ, \;
%     %\lim_{\substack{k\to \infty,\\ k \in \mathcal K}} \tilde{\ve y}_k^2 = \vZ, \;
%     \lim_{\substack{k\to \infty,\\ k \in \mathcal K}} - \tran{\mgkSub} \tilde{\ve y}_k^5
%     =
%     -\tran{\vg} \bar{\ve y}^5
%     =
%     0
%     .
% \label{eqn:week_afford_income}
% \end{equation}
%If we substitute the scaled variables
%into~\eqref{eqn:dual_k_unscaled}, we see
For any dual feasible $[\ve y_k^1, \ldots, \ve y_k^5]$,
the vector $[\tilde{\ve y}_k^1, \ldots, \tilde{\ve y}_k^5]$ is feasible for
the problem 
\begin{equation*}
    \begin{aligned}[t]
        &
        %{\nu}_k \cdot 
        \min_{%
          \tilde{\ve y}^1,%
          \tilde{\ve y}^2,%
          \tilde{\ve y}^3,%
          \tilde{\ve y}^5 \ge \vZ ,%
          \tilde{\ve y}^4 \in \mathbb R^\hDim%
        }
        \begin{bmatrix}
            %\nicefrac{1}{{\nu}_k}
            \tran{\ve 1}_{\varDim}
            &
            %\nicefrac{1}{{\nu}_k}
            \tran{\ve 1}_{\varDim}
            &
            \tran{\ve 0}_{\fDim}
            &
            \tran{\ve 0}_{\hDim}
            &
            -
            \tran{\mgkSub}
        \end{bmatrix}
        \begin{bmatrix}
            \tilde{\ve y}^1 \\ \vdots \\ \tilde{\ve y}^5
        \end{bmatrix}
        \quad \text{s.t.}\quad
        %\\
        % &%
        % \tilde{\ve y}_1 \ge \vZ,
        % \tilde{\ve y}_2 \ge \vZ,
        % \tilde{\ve y}_3 \ge \vZ,
        % \tilde{\ve y}^4 \in \mathbb{R}^\hDim,
        % \tilde{\ve y}^5 \ge \vZ,
        % \\
        %&
        \begin{bmatrix}
            %\nicefrac{-1}{{\nu}_k} 
            \eye{\varDim} & 
            %\nicefrac{1}{{\nu}_k} 
            \eye{\varDim} & \tran{\DFk} & \tran{\DmH} & \tran{\DmG} \\
            \tran{\ve 0}_{\varDim} & \tran{\ve 0}_{\varDim} & \tran{\ve 1}_{\fDim} & \tran{\ve 0}_{\hDim} & \tran{\ve 0}_{\gDim}
        \end{bmatrix}%
        \begin{bmatrix}
            \tilde{\ve y}^1 \\ \vdots \\ \tilde{\ve y}^5
        \end{bmatrix}
        =
        \begin{bmatrix}
            \ve 0_{\varDim} \\ \nicefrac{1}{{\nu}_k}
        \end{bmatrix}
        .
    \end{aligned}
\end{equation*}
% Its limit problem (disregarding the objective factor) is
% \begin{equation}
% \begin{aligned}[t]
% &
% \min_{
%     \tilde{\ve y}^1,
%     \ldots ,
%     \tilde{\ve y}^5
%     }
% \begin{bmatrix}
%     \tran{\ve 0}_{\varDim} & \tran{\ve 0}_{\varDim} &  \tran{\ve 0}_{\fDim} & \tran{\ve 0}_{\hDim} & - \tran{\ve g}
% \end{bmatrix}
% \begin{bmatrix}
%     \tilde{\ve y}^1 \\ \vdots \\ \tilde{\ve y}^5
% \end{bmatrix}
% \qquad \text{s.t.}
% \\
% &\tilde{\ve y}_1 \ge \vZ,
% \tilde{\ve y}_2 \ge \vZ,
% \tilde{\ve y}_3 \ge \vZ,
% \tilde{\ve y}^4 \in \mathbb{R}^\hDim,
% \tilde{\ve y}^5 \ge \vZ,
% \\
% &\begin{bmatrix}
%    %\ve 0_{n,n} & \ve 0_{n,n} &
%    \tran{\DF} & \tran{\DH} & \tran{\DG} \\
%    %\tran{\ve 0}_{\varDim} & \tran{\ve 0}_{\varDim} &
%    \tran{\ve 1}_{\fDim} & \tran{\ve 0}_{\hDim} & \tran{\ve 0}_{\gDim}
% \end{bmatrix}%
% \begin{bmatrix}
%  \tilde{\ve y}^3 \\ \tilde{\ve y}^4 \\ \tilde{\ve y}^5
% \end{bmatrix}
% \begin{matrix}
% = \\ \ge
% \end{matrix}
% \begin{bmatrix}
%  \ve 0_{\varDim} \\ 0
% \end{bmatrix},
% \end{aligned}%
% \label{eqn:dual_scaled}
% \tag{D'}
% \end{equation}
In the limit, the last constraint becomes $\tran{\ve 1}_{\fDim} \bar{\ve y}^3 = 0$,
the limiting problem is 
\begin{equation}
    \begin{aligned}[t]
        &
        \min_{%
          \tilde{\ve y}^1,%
          \tilde{\ve y}^2,%
          \tilde{\ve y}^3,%
          \tilde{\ve y}^5 \ge \vZ ,%
          \tilde{\ve y}^4 \in \mathbb R^\hDim%
        }
        \begin{bmatrix}
            %\nicefrac{1}{{\nu}_k}
            \tran{\ve 1}_{\varDim}
            &
            %\nicefrac{1}{{\nu}_k}
            \tran{\ve 1}_{\varDim}
            &
            \tran{\ve 0}_{\fDim}
            &
            \tran{\ve 0}_{\hDim}
            &
            -\tran{\vg}
        \end{bmatrix}
        \begin{bmatrix}
            \tilde{\ve y}^1 \\ \vdots \\ \tilde{\ve y}^5
        \end{bmatrix}
        \quad \text{s.t.}\quad
        %\\
        % &%
        % \tilde{\ve y}_1 \ge \vZ,
        % \tilde{\ve y}_2 \ge \vZ,
        % \tilde{\ve y}_3 \ge \vZ,
        % \tilde{\ve y}^4 \in \mathbb{R}^\hDim,
        % \tilde{\ve y}^5 \ge \vZ,
        % \\
        %&
        \begin{bmatrix}
            %\nicefrac{-1}{{\nu}_k} 
            \eye{\varDim} & 
            %\nicefrac{1}{{\nu}_k} 
            \eye{\varDim} & \tran{\DF} & \tran{\DH} & \tran{\DG} \\
            \tran{\ve 0}_{\varDim} & \tran{\ve 0}_{\varDim} & \tran{\ve 1}_{\fDim} & \tran{\ve 0}_{\hDim} & \tran{\ve 0}_{\gDim}
        \end{bmatrix}%
        \begin{bmatrix}
            \tilde{\ve y}^1 \\ \vdots \\ \tilde{\ve y}^5
        \end{bmatrix}
        =
        \begin{bmatrix}
            \ve 0_{\varDim} \\ 0
        \end{bmatrix}
        .
    \end{aligned}
    \tag{D'}
    \label{eqn:dual_scaled}
\end{equation}
The primal thus has constant objective value:
\begin{equation}
\begin{aligned}[t]
\max_{\ve d \in \mathbb{R}^\varDim, {\beta}^- \in \mathbb R} 0
\qquad\text{s.t.}\qquad
\begin{bmatrix}
 -\eye{\varDim} & \vZ_{\varDim} \\
 \eye{\varDim} & \vZ_{\varDim} \\
 \DF & \ve 1_{\fDim} \\
 \DH & \ve 0_{\hDim} \\
 \DG & \ve 0_{\gDim}
\end{bmatrix}
\begin{bmatrix}
 {\ve d} \\ {{\beta}}^-
\end{bmatrix}
\begin{matrix}
  \le \\ \le \\ \le \\ = \\ \le
\end{matrix}
\begin{bmatrix}
 \ve 1_{\varDim} \\ \ve 1_{\varDim} \\ \ve 0_{\fDim} \\ \ve 0_{\hDim} \\ - \ve g
\end{bmatrix}.
\end{aligned}%
\label{eqn:primal_scaled}
\tag{P'}
\end{equation}
% For our feasible subsequence $\mathcal K$, we have
% $\tilde{\ve y}_k\to \bar{\ve y} =%
% (\vZ, \vZ, \bar{\ve y}^3, \bar{\ve y}^4, \bar{\ve y}^5)$.
By upper semi-continuity of the feasible set mappings 
(see \cite{wets_continuity_1985}),
$\bar{\ve y}$ is feasible for~\eqref{eqn:dual_scaled},
and we see from~\eqref{eqn:week_afford_income} and~\eqref{eqn:primal_scaled} 
that it is also optimal.
Let $\vd\in \mathbb R^n$ be a direction with $\normInf{\vd} \leq 1$ adhering to the 
\ac{mfcq}.
Furthermore, let $\beta$ be such that $(\vd, \beta)$ is feasible for~\eqref{eqn:primal_scaled}.
Then $(\vd, \beta)$ is optimal for~\eqref{eqn:primal_scaled}, and 
we make the following observations:
\begin{itemize}
    \item \textbf{$\bar{\ve y}^3$ must be zero,} because of the second constraint
    in~\eqref{eqn:dual_scaled}.
    \item \textbf{$\bar{\ve y}^5$ must be zero:}
    Because of the \ac{mfcq}
    %there is a primal optimal solution $\ve d$ of~\eqref{eqn:primal_scaled}, such that
    it holds that
    $\DG_{[\ell, :]} \cdot \ve d < 0$ whenever $\gL_\ell = 0$.
    By complementary slackness, it follows that for these indices $\ell$ that the entries 
    in $\bar{\ve y}^5$ must be $0$.
    From~\eqref{eqn:week_afford_income} it follows for the other indices too, 
    that they must be $0$.
    \item But with $\bar{\ve y}^4 \ne \vZ$ it then follows from the first
    constraint in~\eqref{eqn:dual_scaled} that $\tran{\DH} \cdot \bar{\ve y}^4 =
     \vZ$ in contradiction to the \ac{mfcq}.
\end{itemize}
Hence, for all $i \in \{3,4,5\}$ the sequence $\{\ve{y}_k^i\}$ of (unscaled) Lagrange 
multipliers must also be bounded!\\
We thus can take a subsequence $\mathcal K$ so that ${\{\ve y_k\}}_{k\in \mathcal K}$
converges to some $(\vZ, \vZ, \bar {\ve y}_k^3, \bar{\ve y}_k^4, \bar{\ve y}^5)$ 
and -- by upper semi-continuity of the feasible set mapping -- 
the limit point is feasible for the limiting problem of
~\eqref{eqn:dual_k_unscaled},
which happens to be~\eqref{eqn:buntingGeologically}.
% \begin{equation}
% \begin{aligned}
% &\min_{\ve y^1, \ldots , \ve y^5 }
% \begin{bmatrix}
%     \tran{\ve 1}_{\varDim} & \tran{\ve 1}_{\varDim} &  \tran{\ve 0}_{\fDim} & \tran{\ve 0}_{\hDim} & \tran{\ve g}% 
% \end{bmatrix}
% \begin{bmatrix}
%     \ve y^1 \\ \vdots \\ \ve y^5
% \end{bmatrix}
% \qquad \text{s.t.}
% \\
% &\ve y_1 \ge \vZ, \ve y_2 \ge \vZ, \ve y_3 \ge \vZ, \ve y^4 \in \mathbb{R}^\hDim, \ve y^5 \ge \vZ
% \\
% &\begin{bmatrix}
%     -\eye{\varDim} & \eye{\varDim} & \tran{\DF} & \tran{\DH} & \tran{\DG} \\
%     \tran{\ve 0}_{\varDim} & \tran{\ve 0}_{\varDim} & \tran{\ve 1}_{\fDim} & \tran{\ve 0}_{\hDim} & \tran{\ve 0}_{\gDim}
% \end{bmatrix}%
% \begin{bmatrix}
%     \ve y^1 \\ \vdots \\ \ve y^5
% \end{bmatrix}
% \begin{matrix}
%     = \\ \ge
% \end{matrix}
% \begin{bmatrix}
%     \ve 0_{\varDim} \\ 1
% \end{bmatrix}
% \end{aligned}
% \label{eqn:dual_limit_unscaled}
% \tag{D}
% \end{equation}
The optimal value at
$(\vZ, \vZ, \bar {\ve y}_k^3, \bar{\ve y}_k^4, \bar{\ve y}_k^5)$ is 0
due to strong duality.
Because the corresponding primal is~\eqref{eqn:legitimizersSannop},
it follows from \cref{thm:subproblems_implies_kkt}
that $\bar{\vx}$ is a KKT point of~\eqref{eqn:mop_nonlin}.
\end{proof}

With~\cref{thm:subproblems_implies_kkt}, it is now easy to derive the main result:
\begin{mytheorem}[Convergence to KKT-points]
  Suppose the same assumptions as in~\cref{thm:convergence1} hold 
  and that $\{\vxk\}$ is a quasi-stationary subsequence with limit point
  $\bar{\vx}$.
  If the \ac{mfcq} hold at $\bar{\vx}$, 
  then $\bar{\vx}$ is a KKT-point of~\eqref{eqn:mop_nonlin}.
\end{mytheorem}

%-----------------------sec05 ---------------------------------
\section{Numerical Examples}%
\label{sec:examples}

In this section, we provide numerical examples for which we have applied our
algorithm to two test problems.
The algorithm is implemented in the Julia language according to the pseudocode
in~\cref{section:algorithm}.
% without any focus on performance.
We also share our code as a Pluto~\cite{julia_pluto} notebook.
\footnote{\url{https://gist.github.com/manuelbb-upb/69582b2322346485333f01807a2c241c}}
Before describing the experiments, we would like to note that the simple
implementation neglects many interesting questions and details that could
(and should) be explored, such as optimized model construction algorithms.
Nonetheless, it serves our primary goal:
To demonstrate the general viability of the algorithm.

\subsection{Constrained Two-Parabolas Problem}

\begin{figure}
  \centering
  \includegraphics[width=\linewidth]{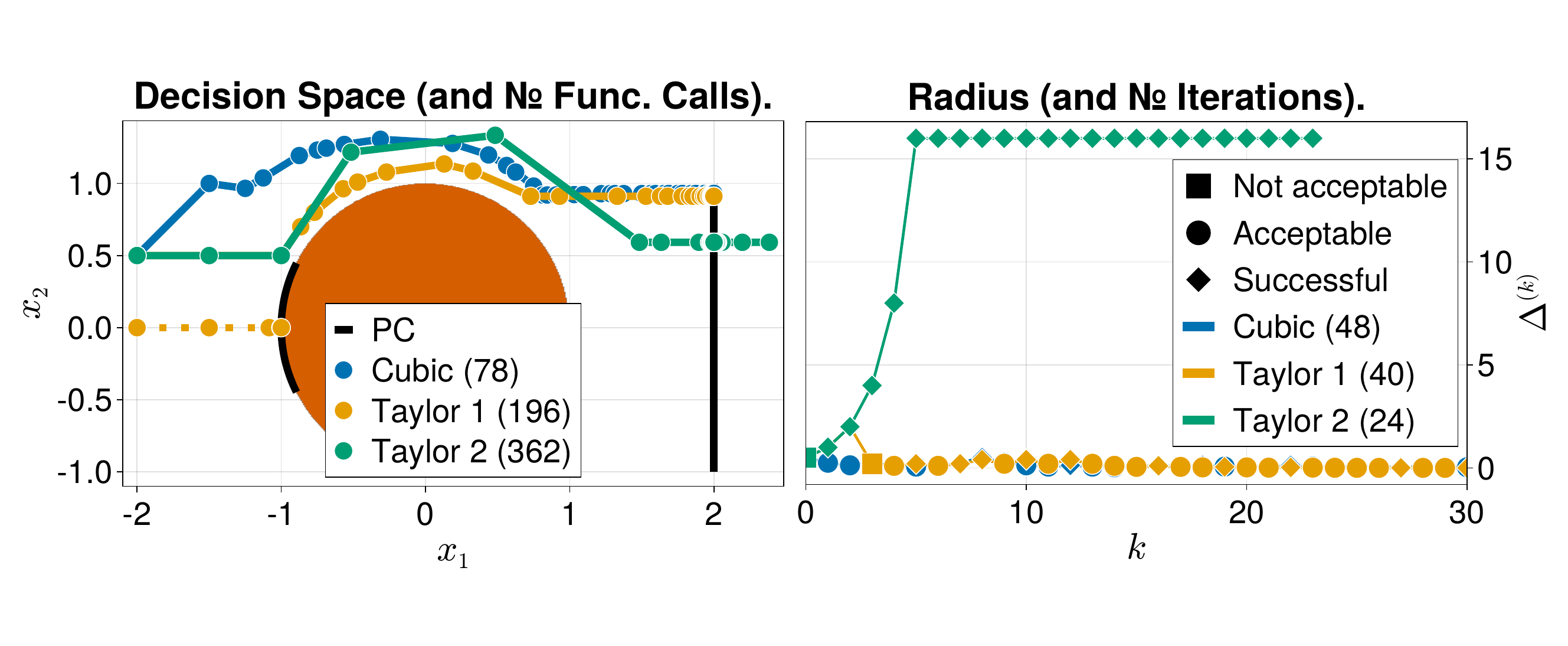}%
  \caption{Results for runs of our algorithm on~\eqref{eqn:bennets_problem}.
    Left: Iteration trajectories in decision space for different model types.
    The number of function calls for each model is given in parentheses.
    Right: Radius evolution for the different model types and iteration
    classification with respect to the metric $\rho\KK$.
    Plots created with Makie~\cite{julia_makie}.%
  }%
  \label{fig:example1}
\end{figure}

Variations of the two parabolas problem are popular test cases for \ac{moo}
methods.
The objectives are simply two $\varDim$-variate parabolic functions and
-- in the unconstrained case -- the Pareto Set is the line connecting their
respective minima.
The following constrained version is taken from~\cite{gebken_constraints}.
\begin{equation}
  \begin{aligned}
    &\min_{\vx\in \mathbb{R}^2}
      \begin{bmatrix}
        (x_1 - 2)^2 + (x_2 - 1)^2 \\
        (x_1 - 2)^2 + (x_2 + 1)^2
      \end{bmatrix}
    \quad\text{s.t.}\quad
    \gL(\vx) = 1 - x_1^2 - x_2^2 \le 0.
  \end{aligned}
  \label{eqn:bennets_problem}
  \tag{Ex. 1}
\end{equation}
The feasible set of \eqref{eqn:bennets_problem} is $\mathbb{R}^2$ without the
interior of the unit ball.
In~\cref{fig:example1} the infeasible area is shaded red.
The Pareto critical set is the line connecting $[2, -1]$ and $[2, 1]$ and the left boundary
of the unit ball:
\begin{equation*}
  \mathcal P_c =
  \left\{
    \begin{bmatrix}
      2\\
      s
    \end{bmatrix}:
    s\in[-1,1]
  \right\}
  \bigcup
  \left\{
    \begin{bmatrix}
      \cos t\\
      \sin t
    \end{bmatrix}
    :
    t \in [\pi - \theta, \pi + \theta], \theta = \arctan\left(\frac{1}{2}\right)
  \right\}.
\end{equation*}
We have plotted the critical set with black lines in~\cref{fig:example1}.

For this problem, we first applied our algorithm three times with the same
parameters and beginning at $\vx_0 = [-2, 0.5]$, but with different model types.
We compare \ac{rbf} models with first and second degree Taylor polynomials of
the objective and constraint functions.
The \ac{rbf} models use the cubic kernel and their construction is based on the work
in \cite{wild_orbit_2008}.
The derivatives for the Taylor models are approximated using finite differences
(to conform to the assumption that gradients are not available exactly).
The other parameters are
\begin{equation*}
  \begin{gathered}
    \DeltaInit = 0.5,
    \DeltaMax = 2^5 \cdot \DeltaInit,
    \gammass = 0.1,
    \gammas = 0.5,
    \gammag = 2.0,
    \CCaccept = 0.01,
    \CCsuccess = 0.9,\\
    \CCcritTestchi = 0.1,
    \CCcritTesttheta = 0.1,
    \kappa_{\theta} = 10^{-4},
    \psi = 2,
    \CCBeta = 1000,
    \CCMu = 3000,
    \CCalpha = 0.5,
    \CCDelta = 0.7
    \CCmu = 100,
    \mu = 0.01.
  \end{gathered}
\end{equation*}
These settings largely follow the recommendations
in~\cite{fletcher_global_2002}.
Additionally, we stop after 100 iterations or if the criticality loop has
executed once or if the trial point is accepted, and it holds that
$\|\vxk - \vxSk\| \le 10^{-5} \|\vxk\|$ or
$\|\vf(\vxk) - \vf(\vxSk)\| \le 10^{-5} \|\vf(\vxk)\|$.

The results are depicted in~\cref{fig:example1}.
As can be seen, all runs converge to the Pareto critical set and avoid
the infeasible area.
In terms of objective evaluations, the \ac{rbf} models require significantly
less function calls than the Taylor models.
(If exact gradients are used the numbers seem to be roughly equal.)
In terms of iterations, the second degree Taylor polynomials require the least.\\
Many gradient-based \ac{moo} algorithms suffer from bias towards
individual minima and, indeed, we also see that the final iterates are
close to the minimum $[2, 1]$.
There also is a more notable bias of the \ac{rbf} models towards that minimum
at the beginning of the optimization. This is very likely due to the infinity-norm 
being used and because of their construction:
In the first iteration, only 3 objective evaluations along the coordinate axes
around
$\vx\itCmd{0}$ are used in a large trust-region, so that the predicted descent
is very inaccurate.
In later iterations, the trust-regions are smaller, and the models become more
precise by including previous evaluations from a database.
The Taylor models of degree 2 can be considered exact for
\eqref{eqn:bennets_problem}.
This is why, for the green line, we see the typical zig-zag motion close to the
critical set.
Every trial point achieves the predicted objective decrease exactly and hence
every
iteration is deemed very successful, which makes the trust-regions have maximum
radius and practically reduces our method to ``vanilla''
gradient-descent.\\
Finally, it is worth-wile to compare the iteration trajectories to those
in~\cite{gebken_constraints}, where an active set strategy is used.
In contrast to the results presented there, we stay further away from
the infeasible area. This is because every step is computed with inexact
approximations for \emph{every} constraint function, instead of
only considering the active constraints.
If we start closer to the infeasible set, however, we also find a critical point on
its boundary: The trajectory depicted by the dashed orange line starts at $[-2,
0]$ and ends close to $[-1, 0]$.

\subsection{Non-Convex Test-Problem ``W3''}

\begin{figure}
  \centering
  \includegraphics[width=\linewidth]{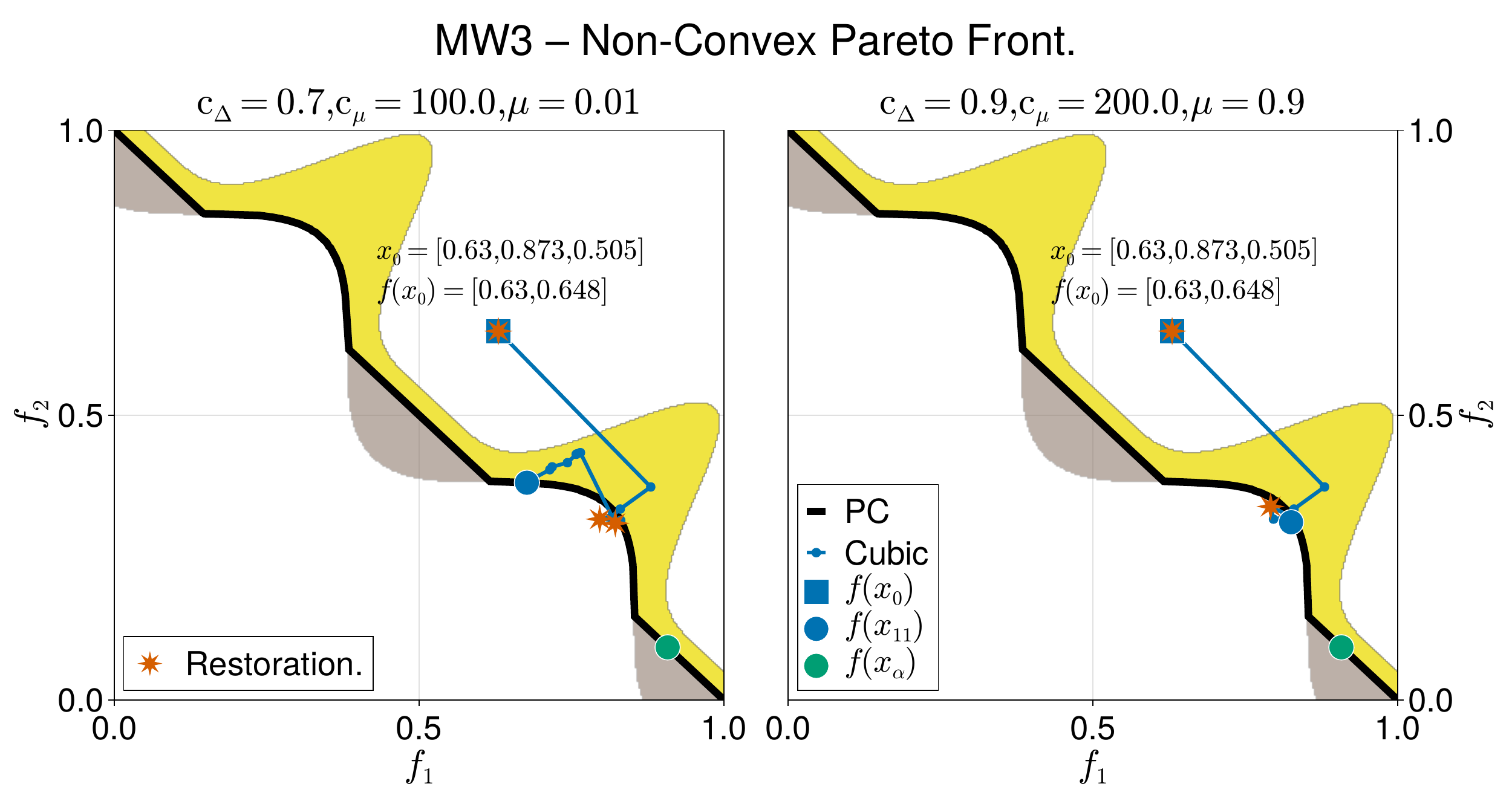}%
  \caption{Results for the MW3 problem and comparison with the weighted-sum
    approach.
    For the right graphic, the compatibility parameters were changed to avoid
    restoration iterations.
    In both plots, blue markers show our algorithm and the green marker is the
    weighted-sum result.
    The yellow area is attainable.%
  }%
  \label{fig:example2}
\end{figure}

The second example is problem ``W3'' taken from~\cite{ma_evolutionary_2019}.
The problem has box constraints and multiple non-linear inequality constraints.
Whilst the unconstrained Pareto Front is simply a line in objective space,
the constraints make it partially non-convex.
We have chosen the problem with $3$ variables and $2$ objectives to compare our
algorithm against a simple weighted-sum scalarization.
It reads:
\begin{equation}
\begin{aligned}
    &\min_{\vx \in \mathbb{R}^3}
      \begin{bmatrix}
        f_1(\vx) \\
        f_2(\vx)
      \end{bmatrix} =
      \min_{\vx \in \mathbb{R}^3}
      \begin{bmatrix}
        x_1\\
        d(\vx)\cdot\left(1 - \frac{f_1(\vx)}{d(\vx)}\right)
      \end{bmatrix}
      \qquad&\text{s.t.}
  \\
    &c_1(\vx) =  f_1(\vx) + f_2(\vx) - 1.05 - 0.45\sin(0.75 \pi \cdot l(\vx))^6 \le 0, \\
    &c_2(\vx) = -f_1(\vx) - f_2(\vx) + 0.85 + 0.3\sin(0.75 \pi \cdot l(\vx))^2 \le 0,\\
    &d(\vx) = 1 + 2( x_2 + (x_1 - 0.5)^2 - 1)^2 +2( x_3 + (x_2 - 0.5)^2 - 1)^2, \\
    &l(\vx) = \sqrt{2}(f_2(\vx) - f_1(\vx)).\\
\end{aligned}
\tag{Ex. 2}
\label{eqn:mw3}
\end{equation}
Because the constraints $c_1$ and $c_2$ describe the feasible set in terms of
the objective values we can show the attainable objective values in objective
space in~\cref{fig:example2}.
The critical set is again shown in black.

For the first run (shown on the left), we largely used the same parameters as in
the previous section, but relaxed the stopping criteria a bit.
For example, the relative stopping tolerances were reduced to $10^{-6}$ and up
to 500 iterations were allowed.
Whilst the step problems are modeled using JuMP~\cite{julia_jump}
and solved with any suitable LP or QP solver, like COSMO~\cite{julia_cosmo},
restoration uses the NLopt solver~\cite{nlopt} and its
local ``COBYLA'' algorithm.
The same non-linear solver and similar stopping criteria were used to solve a
weighted-sum scalarization of~\eqref{eqn:mw3} with
objective $\nicefrac{1}{2} f_1 + \nicefrac{1}{2} f_2$.
The blue trajectory in~\cref{fig:example2} shows that our method with \ac{rbf}
models
can find a critical point on the non-convex part of the Pareto Front while the
weighted-sum solution always must belong to the linear part.\\
Because the initial iterate was not feasible, a restoration step has to be
performed
right at the beginning. This can cause the final iterate to be far away from the
initial site.
We see can similarly observe in the left graphic, that there are two additional
restorations close to the concave knee of the Pareto Front, which again prevents
convergence.
For comparison, in a second run, we tried to relax the compatibility parameters.
The right plot shows that indeed one restoration can be avoided and the iterates
do not deviate that much.
%-----------------------sec06 ---------------------------------
\section{Conclusion and Outlook}%
\label{sec:conclusion}

In this article, we have presented an algorithm for non-linearly constrained
\acp{mop}.
The method does not necessarily need derivative information of the functions
but can use surrogate models instead.
Additionally, a Filter ensures convergence towards feasibility.
We have proven convergence of an algorithmic sub-sequence to Pareto-criticality
and confirm the theoretical results with numerical experiments.

The experiments have shown both promising features of the algorithm compared to
naive scalarizations and potential drawbacks.
\begin{itemize}
\item Like many gradient-based \ac{moo} algorithms, there often seems to be bias towards
  individual objective minima.
  Thus, an important task for future research lies in trying to remedy this
  behavior,
  maybe by using some gradient scaling.
\item In this regard, it would be desirable to be able to use alternative descent
  direction
  that allow guidance of the iterates or use momentum to accelerate convergence.
\item Moreover, we believe it possible to transfer the model construction
  optimizations from \cite{wild_orbit_2008,ours} into our algorithm.
  This would allow for ``model-improvement'' steps with models that are not fully
  linear and could potentially save expensive function evaluations.
\item Finally, our algorithm produces only one critical point.
  After a run has finished, can we leverage the surrogates we already have
  constructed to warm-start the optimization of a modified \ac{mop}
  to find additional critical
  points and obtain a covering of (parts) of the critical set?
\end{itemize}

A more theoretical question is whether or not our convergence results can be
strengthened. To this end, using a slanted filter like in~\cite{ferreira_global_2017} might prove
beneficial.
Furthermore, in~\cref{sec:other_algos} we have already talked about why an approach similar to that 
in~\cite{eason_trust_2016} could proof beneficial.

%--------------------------------appendix ---------------------------------------------tt
\section*{Appendix}

\subsection*{Equivalence of Inexact Criticality with Different Norms}

To prove \cref{thm:critical_values_equivalent} we transfer the corresponding single-objective
results from \cite{conn_global_1993} to the multi-objective case.
This works relatively straightforward, but to the best of our knowledge the multi-objective
results have not yet been published anywhere, so the proofs are given in detail.
First, we make use of the following auxiliary
result:

\begin{mylemma}[{\cite[Lemma 2]{conn_global_1993}}]
	\label{thm:delouserProthalamion}
	Suppose \cref{%
	ass:covering,%
	ass:lipschitz_gradients_true,%
	ass:fully_linear_models,%
	ass:equivalent_norms%
	}
	hold and that ${\omega}$ is defined by
	\begin{equation}
		\begin{aligned}
			{\omega}\left(  \vxNk; \mfk, \linFeasApproxk, \normIt{{\bullet}}, {\vartheta} \right)
			=-
			&
			\min_{\ve d \in \mathbb{R}^\varDim}
			\max_{\ell}
			\langle
				\dmf_{\ell}(\vxNk), \ve d
			\rangle
			\quad\text{s.t.}\quad
			%\\
			%&\DmF(\vxNk) \cdot \vd \le {\beta},
			\vd \in L_k,
			\;
			%\normTr{\vnk + \vd} \le \Deltak.
			%\normTr{\vnk + \vd} \le 1.
			\normIt{\vd} \le {\vartheta},
		\end{aligned}
		\label{eqn:itr_tangential_theta}
	\end{equation}
	where $L_\kk = \left(\linFeasApproxk - \vxNk\right)$.
	For any $\kk\in \mathbb{N}_0$, for which the normal step exits,
	the following statements hold:
	\begin{enumerate}
		\item The function
		${\vartheta} \mapsto {\omega}\left(  \vxNk; \mfk, \linFeasApproxk, \normIt{{\bullet}}, {\vartheta} \right)$ is continuous
		and non-decreasing for ${\vartheta}\ge 0$.
		\item The function
		${\vartheta} \mapsto \frac{
			{\omega}\left(  \vxNk; \mfk, \linFeasApproxk, \normIt{{\bullet}}, {\vartheta} \right)
		}{{\vartheta}}$
		is non-increasing for ${\vartheta}> 0$.
	\end{enumerate}
\end{mylemma}

\begin{proof}
	Let $\kk\in \mathbb{N}_0$ be an iteration index with algorithmic
	variables $(\vxNk, \linFeasApproxk, \normIt{{\bullet}})$.
\begin{enumerate}
	\item Same as in \cite{conn_global_1993}, the feasible set mapping
	$$
	{\vartheta} \mapsto
	D_\kk({\vartheta})
	=
	\left\{
	\vd \in \mathbb{R}^\varDim:
	\normIt{\ve d} \le {\vartheta},
	\vd \in \left(\linFeasApproxk - \vxNk\right)
	\right\}
	$$
	is continuous for all ${\vartheta}\ge 0$.
	The function
	$\Psi( {\vartheta}, \ve d) = \max_{\ell}
		\langle
			\dmf_{\ell}(\vxNk), \ve d
	\rangle$ is defined on $\mathbb{R}_{\ge 0} \times \mathbb{R}^\varDim$
	and for each ${\vartheta}\ge 0$ it is convex in its second argument.
	Hence, we can use results from \cite{fiacco}, which guarantee
	the continuity of the optimal value
	$$
	- \min_{\ve d \in D_\kk({\vartheta})} \Psi({\vartheta}, \ve d),
	=
	{\omega}\left(  \vxNk; \mfk, \linFeasApproxk, \normIt{{\bullet}}, {\vartheta} \right)$$
  which is non-decreasing with respect to $\vartheta$ by definition.
	\item Consider $0<{\vartheta}_1<{\vartheta}_2$ and $\vd_1, \vd_2\in \mathbb{R}^\varDim$ such that
	\begin{alignat}{3}
	{\omega}\left(  \vxNk; \mfk, \linFeasApproxk, \normIt{{\bullet}}, {\vartheta}_1 \right)
		&=
			-\max_{\ell} \langle \dmf_{\ell}(\vxNk), \vd_1 \rangle,
			\;
			&
				\normIt{\vd_1}
				&
					\le {\vartheta}_1,
			\;
			&
				\vd_1 \in L_\kk
	\label{eqn:derivativelySentimo}
	\\
	{\omega}\left(  \vxNk; \mfk, \linFeasApproxk, \normIt{{\bullet}}, {\vartheta}_2 \right)
		&=
			-\max_{\ell} \langle \dmf_{\ell}(\vxNk), \vd_2 \rangle,
			\;
			&
				\normIt{\vd_2}
				&
					\le {\vartheta}_2,
			\;
			&
				\vd_2 \in L_\kk
	\label{eqn:untrackedEtaerio}
	\end{alignat}
	Because the set $L_\kk$ is convex, it follows from
	$\nicefrac{{\vartheta}_1}{{\vartheta}_2} < 1$, $\vZ \in L_\kk$ and $\vd_2 \in L_\kk$ that
	also
	$
	\nicefrac{{\vartheta}_1}{{\vartheta}_2} \vd_2 \in L_\kk
	$.
	Moreover, it follows from \eqref{eqn:derivativelySentimo} that
	$$
	\normIt{ \frac{{\vartheta}_1}{{\vartheta}_2} \vd_2 }
	=
	\frac{{\vartheta}_1}{{\vartheta}_2} \normIt{\vd_2} \le {\vartheta}_1.
	$$
	Thus, $\nicefrac{{\vartheta}_1}{{\vartheta}_2} \vd_2$ is feasible for the problem in \eqref{eqn:untrackedEtaerio}.
	Consequently,
	% remembering that the maximum is positively homogenous and
	% the inner product is bi-linear:
	$$
	\frac{{\omega}\left(  \vxNk; \mfk, \linFeasApproxk, \normIt{{\bullet}}, {\vartheta}_1 \right)}{{\vartheta}_1}
	\ge
	-
	\frac{1}{{\vartheta}_1}
	\max_{\ell}
	\left\langle
	\dmf_{\ell}(\vxNk), \frac{{\vartheta}_1}{{\vartheta}_2} \vd_2
	\right\rangle
	=
	\frac{
	{\omega}\left(  \vxNk; \mfk, \linFeasApproxk, \normIt{{\bullet}}, {\vartheta}_2 \right)
	}{{\vartheta}_2}.
	$$
\end{enumerate}
\end{proof}

We can now proof~\cref{thm:critical_values_equivalent}, which states that
uniformly equivalent norms imply uniformly equivalent inexact criticality values,
i.e.,
\begin{equation}
	\frac{1}{\CCcrit}
	{\omega}\left(  \vxNk; \mfk, \linFeasApproxk, \normIt{{\bullet}} \right)
	\le
	{\omega}\left(  \vxNk; \mfk, \linFeasApproxk, \normTwo{{\bullet}} \right)
	\le
	\CCcrit \cdot
	{\omega}\left(  \vxNk; \mfk, \linFeasApproxk, \normIt{{\bullet}}, \right).
	\label{eqn:overclaimedOrchestrate}
\end{equation}

\begin{proof}
	The proof works similarly to the single-objective case \cite[Th. 4]{conn_global_1993}.
	Let $\kk\in \mathbb{N}_0$ be such that the normal step exists.
	We first make the following observations:
	\begin{enumerate}
		% \item \label{item:sturtsNodes}%
		% For the minimizer $\vdk$ in \Cref{eqn:itr_tangential} it holds that
		% $\normIt{\vdk} \le 1$. If equality holds, $\normIt{\vdk} = 1$,
		% then with \Cref{ass:equivalent_norms} it follows that
		% \begin{equation}
		% \frac{1}{\CCnorm}
		% \le
		% \normTwo{\vdk}
		% \le
		% \CCnorm.
		% \label{eqn:hypostasisingPlaydough}
		% \end{equation}
		\item \label{item:denominatingConcocted}%
		The ball defined by
		$\normTwo{\vd} \le \frac{1}{\CCnorm}$ is contained in the ball defined by
		$\normIt{\vd} \le 1$,
		% again
		due to~\cref{ass:equivalent_norms}.
		\item \label{item:noninterchangeablePhenoxy}%
		Likewise, the ball with $\normIt{\vd} \le 1$ is contained in the ball defined
		by $\normTwo{\vd} \le \CCnorm$.
	\end{enumerate}
	According to \eqref{eqn:itr_tangential_theta}, we define
	\begin{align*}
	{\omega}_\text{max} =
	{\omega}\left(  \vxNk; \mfk, \linFeasApproxk, \normTwo{{\bullet}}, \CCnorm \right)
  \quad\text{and}\quad
	{\omega}_\text{min} =
	{\omega}\left(  \vxNk; \mfk, \linFeasApproxk, \normTwo{{\bullet}}, \CCnorm^{-1} \right).
	\end{align*}
	From the second statement in~\cref{thm:delouserProthalamion} it then
	follows that
	\begin{equation}
	{\omega}_\text{max} \le \CCnorm^2 {\omega}_\text{min}.
	\label{eqn:heldentenorsAlineation}
	\end{equation}
	If
	$
	{\omega}\left(  \vxNk; \mfk, \linFeasApproxk, \normTwo{{\bullet}} \right)
	=
	{\omega}\left(  \vxNk; \mfk, \linFeasApproxk, \normIt{{\bullet}} \right)
	$
	there is nothing to show.
	Hence, first we assume that
	\begin{equation}
	{\omega}\left(  \vxNk; \mfk, \linFeasApproxk, \normTwo{{\bullet}} \right)
	<
	{\omega}\left(  \vxNk; \mfk, \linFeasApproxk, \normIt{{\bullet}} \right),
  \label{eqn:antiorganizationTalkinesses}
  ,
	\end{equation}
	and we again take the respective minimizers $\vdk, \vd_2 \in \mathbb{R}^\varDim$
	with
	\begin{alignat}{3}
	{\omega}\left(  \vxNk; \mfk, \linFeasApproxk, \normIt{{\bullet}} \right)
		&=
			-\max_{\ell} \langle \dmf_{\ell}(\vxNk), \vdk \rangle,
			\;
			&
				\normIt{\vdk}
				&
					\le 1,
			\;
			&
			\vdk \in L_\kk
	\label{eqn:cataphyllPremillennialist}
	\\
	{\omega}\left(  \vxNk; \mfk, \linFeasApproxk, \normTwo{{\bullet}} \right)
		&=
			-\max_{\ell} \langle \dmf_{\ell}(\vxNk), \vd_2 \rangle,
			\;
			&
				\normTwo{\vd_2}
				&
					\le 1,
			\;
			&
				\vd_2 \in L_\kk
	\label{eqn:materielsMashgihim}
	\end{alignat}
	Then
	\begin{equation}
	\frac{1}{\CCnorm}
	\le
	\normTwo{\vdk}
	\le
	\CCnorm
	\quad\text{and}\quad
	\frac{1}{\CCnorm}
	\le
	\normTwo{\vd_2}
	\le
	\CCnorm
  .
	\label{eqn:fractiousnessBowdleriser}
	\end{equation}
	The upper bounds are trivial because of $\CCnorm \ge 1$.
	Suppose the first lower bound is violated, i.e.,
	$\normTwo{\vdk} < \nicefrac{1}{\CCnorm}.$
	According to observation \ref{item:denominatingConcocted},
	it then follows that $\normIt{\vdk} \le 1$.
	The vector $\vdk \in L_\kk$ is then also feasible for \eqref{eqn:materielsMashgihim} and
	the optimality of $\vd_2$ implies
	$$
	\underbrace{%
	- \max_{\ell}
		\langle \dmf_{\ell}(\vxNk), \vd_2 \rangle
	}_{= {\omega}\left(  \vxNk; \mfk, \linFeasApproxk, \normTwo{{\bullet}} \right)}
	\ge
	\underbrace{%
	- \max_{\ell}
		\langle  \dmf_{\ell}(\vxNk), \vdk \rangle
	}_{%
	= {\omega}\left(  \vxNk; \mfk, \linFeasApproxk, \normIt{{\bullet}} \right)
	},
	$$
	in contradiction to \eqref{eqn:antiorganizationTalkinesses}.\\
	Suppose the second lower bound is violated, i.e.,
	$
	\normTwo{\vd_2} < \nicefrac{1}{\CCnorm}.
	$
	According to our observations from above, for $\vd_2$ we also have $\normIt{\vd_2} \le 1$,
	and it is therefore feasible for the problem in \eqref{eqn:cataphyllPremillennialist}.
	With \eqref{eqn:antiorganizationTalkinesses} we even see that it is strictly optimal, contradicting
	the optimality of $\vdk$ in \eqref{eqn:cataphyllPremillennialist}.
	Thus, the second lower bound in \eqref{eqn:fractiousnessBowdleriser} must hold, too.\\
	Equation \eqref{eqn:fractiousnessBowdleriser} shows that both vectors $\vdk$ and $\vd_2$ are feasible for the
	problems defining ${\omega}_\text{max}$ and ${\omega}_\text{min}$.
	We again apply the definition \eqref{eqn:itr_tangential_theta} to see that
	$$
	{\omega}_\text{min} \le {\omega}\left(  \vxNk; \mfk, \linFeasApproxk, \normTwo{{\bullet}} \right) \le {\omega}_\text{max}
	\quad\text{and}\quad
	{\omega}_\text{min} \le {\omega}\left(  \vxNk; \mfk, \linFeasApproxk, \normIt{{\bullet}} \right) \le {\omega}_\text{max}.
	$$
	Plugging in \eqref{eqn:antiorganizationTalkinesses} and \eqref{eqn:heldentenorsAlineation} results in
	$$
	\begin{aligned}
	{\omega}\left(  \vxNk; \mfk, \linFeasApproxk, \normTwo{{\bullet}} \right)
	<
	{\omega}\left(  \vxNk; \mfk, \linFeasApproxk, \normIt{{\bullet}} \right)
	<
	{\omega}_\text{max}
	\le
	\CCnorm^2 {\omega}_{\text{min}}
	\le
	\CCnorm^2
		{\omega}\left(  \vxNk; \mfk, \linFeasApproxk, \normTwo{{\bullet}} \right),
	\end{aligned}
	$$
	implying
	\begin{equation}
	\frac{1}{\CCnorm^2}
		{\omega}\left(  \vxNk; \mfk, \linFeasApproxk, \normIt{{\bullet}} \right)
	\le
		{\omega}\left(  \vxNk; \mfk, \linFeasApproxk, \normTwo{{\bullet}} \right)
	\le
	\CCnorm^2
		{\omega}\left(  \vxNk; \mfk, \linFeasApproxk, \normIt{{\bullet}} \right)
		.
	\label{eqn:parklySubmicrons}
	\end{equation}

	The case
	$$
	{\omega}\left(  \vxNk; \mfk, \linFeasApproxk, \normTwo{{\bullet}} \right)
	<
	{\omega}\left(  \vxNk; \mfk, \linFeasApproxk, \normIt{{\bullet}} \right)
	$$
	is treated analogously and also leads to \eqref{eqn:parklySubmicrons}.
	Thus, \eqref{eqn:overclaimedOrchestrate} holds for any constant
	$\CCcrit \ge \CCnorm^2 \ge 1$ and
	\Cref{thm:critical_values_equivalent} is valid.
\end{proof}

\subsection*{Sufficient Decrease via Backtracking}
This section of the appendix is concerned with justifying the sufficient decrease
bound~\eqref{eqn:sufficient_decrease} in \Cref{ass:sufficient_decrease}.
In the single-objective case, there are many possibilities to achieve the bound,
some of which, e.g., exact and inexact line-search, we have previously shown to work in the
multi-objective as well if the global feasible set is compact and convex~\cite{ours}.
Now at least the approximated linearized feasible sets are convex, which allows us to
cite and utilize the following result:

\begin{mylemma}[Backtracking Decrease, {\cite[Lemma 3]{ours}}]%
  \label{thm:backtracking}
  Let $\linFeasLetter$ be a convex set and $\vd$ be a descent direction for $\mf$ at $\vx \in \linFeasLetter$ and let $\stepsize\ge 0$ be a step-size such that
  $\vx + \nicefrac{\stepsize}{\norm{\vd}} \cdot \vd \in \linFeasLetter$
  and let $\norm{{\bullet}}$ be any vector norm that fulfills~\eqref{eqn:equivalent_norms}.
  Then, for any fixed constants $\armijoConstant, \backtrackConstant \in (0,1)$ and for
  $\scalShort = \scal{\mf}$
  %or $\scalShort = \mfL_{\ell}, {\ell}\in \{0,\ldots , \fDim\}$,
  there is an integer $j\in \mathbb{N}_0$ such that
  \begin{equation}
      \scalShort \left( \ve x \right) -
      \scalShort\left( \vx + \frac{\backtrackConstant^j \stepsize}{\norm{\vd}} \vd \right)
      \ge
      \armijoConstant \frac{\stepsize \backtrackConstant^j}{\norm{\vd}}
      {\omega},
      \label{eqn:backtracking_eqn}
  \end{equation}
  where ${\omega} = -\max_{\ell} \langle \dmf_{\ell}(\vx), \vd \rangle$.
  Moreover,
  there is a constant $\CCsd \in (0,1)$ such that
  if $j$ is the smallest $j\in \mathbb{N}_0$ that satisfies~\eqref{eqn:backtracking_eqn},
  then
  \begin{equation*}
      \scalShort \left( \ve x \right)
      -
      \scalShort\left( \vx + \frac{\backtrackConstant^j \stepsize}{\norm{\vd}} \vd \right)
      \ge
      \CCsd {\omega}
      \min \left\{
      \frac{{\omega}}{ \norm{\vd}^2 \CCnorm^2 H},
      \frac{\stepsize}{\norm{\vd}}
      \right\},
  \end{equation*}
  where
  \begin{equation}
      H = \max_{\vd \in \linFeasLetter - \vx}
      \max_{{\ell}=1,\ldots ,\fDim}
          \normTwo{
              \hess \mfL_{\ell}( \ve x + \vd )
          }
      \text{ s.t. } \norm{\vd}\le \backtrackConstant^j\stepsize.
      \label{eqn:hess_norm_max}
  \end{equation}
\end{mylemma}

For~\cref{thm:backtracking} to be applicable in our setting we require the surrogate function
$\mf$ to be twice continuously differentiable and defined on $\covering$,
which is guaranteed by
\cref{%
	ass:covering,%
	ass:lipschitz_gradients_true,%
	ass:fully_linear_models%
}
Furthermore,~\cref{ass:hess_norm_bounded} may be used instead of~\eqref{eqn:hess_norm_max} to
bound the model Hessians, which appear due to a Taylor approximation of
the components of $\mf$ in the proof of~\cref{thm:backtracking}.
\\
We now have to deal with the possibly different norms $\normTrIt{{\bullet}}$ 
(defining the trust-region) and $\normTr{{\bullet}}$ (used in \eqref{eqn:itr_tangential})
and how to choose the initial backtracking step-length $\stepsize$ in~\cref{thm:backtracking}.
We want to choose $\stepsize \ge 0$ as large as possible and so that for a 
descent direction $\vdk$ it holds that $\normTrIt{\vnk + \stepsize\vdk} \le \Deltak$ 
and $\vxNk + \stepsize \vdk \in \linFeasApproxk$.
Luckily, there is the following bound on the optimal $\stepsize$:

\begin{mylemma}
For $\vxk$ let $\linFeasApproxk$ be the linearized feasible set and $\trk$ be a trust-region of radius
$\Deltak$ w.r.t.  $\normTrIt{{\bullet}}$.
Let $\vdk$ be a minimizer of \eqref{eqn:itr_tangential} with $\normIt{\vdk}\le 1$.
Then there is an initial step-length $\bar{\stepsize} \ge 0$ with
$\vxNk + \frac{\bar{\stepsize}}{\normTrIt{\vdk}} \vdk \in \linFeasApproxk \cap \trk$ 
and
\begin{equation}
    \bar{\stepsize} \ge
    \min
        \left\{
            \Deltak - \normTrIt{\vnk},
            \normTrIt{\vdk}
        \right\}.
    \label{eqn:step-size_init}
\end{equation}
\end{mylemma}

\begin{proof}
  There are, of course, better ways to determine $\bar{\stepsize}$, 
  but if $\normTrIt{\vnk + \vdk} \leq \Deltak$ we can always choose 
  $\bar{\stepsize} = \normTrIt{\vdk}$.
  If $\normTrIt{\vnk + \vdk} > \Deltak$ we can equate either side of the 
  triangle inequality
  $$
	\normTrIt{\vnk} + \bar{\stepsize}
	  \ge
	  \normTrIt{\vnk + \frac{\bar{\stepsize}}{\normTrIt{\vdk}} \vdk}
	$$
  with $\Deltak$ and solve for $\bar{\stepsize}$.

\end{proof}

%%%%%%%%%%%%%%%%%%%%%%%%%%%%%%%%%%%%%%%%%%%%%%%%%%%%%%%%%%%%%%%%m

Finally, we are able to derive the sufficient decrease bound \eqref{eqn:sufficient_decrease} when backtracking is used.
Of course, in case that $\vxk$ is critical, the bound \eqref{eqn:sufficient_decrease} is automatically fulfilled.
Else, we use~\cref{thm:backtracking} with $\bar{\stepsize}$ satisfying \eqref{eqn:step-size_init}.
\Cref{ass:equivalent_norms} and $\normIt{\vdk} \le 1$ together imply
\begin{equation}
	\normTrIt{\vdk} \le \CCnorm^2 \normIt{\vdk} \le \CCnorm^2.
	\label{eqn:evolvesBanoffees}
\end{equation}
\Cref{thm:backtracking} and the fact that we assume $\vnk$ to be compatible lead to
\begin{align*}
    &
    \scalk \left( \vxNk \right)
    -
    \scalk \left( \vxNk + \stepsizek \vdk \right)
    \\
    &\qquad\ge
    \CCsdT \omegak
    \min \left\{
    \frac{\omegak}{ \normTrIt{\vdk}^2 \CCnorm^2 \CChessF},
    \frac{\bar{\stepsize}}{\normTrIt{\vdk}}
    \right\}
    \\
    &\qquad\stackrel{\mathclap{\eqref{eqn:step-size_init}}}\ge
    \CCsdT \omegak
    \min \left\{
    \frac{\omegak}{ \normTrIt{\vdk}^2 \CCnorm^2 \CChessF},
    \frac{ \Deltak - \normTrIt{\vnk}}{\normTrIt{\vdk}},
    \frac{\normTrIt{\vdk}}{\normTrIt{\vdk}}
    \right\}
    \\
    &\qquad\stackrel{\mathclap{\eqref{eqn:evolvesBanoffees},\eqref{eqn:compatible_normal_step}}}%
	\ge
	\quad
    \CCsdT\omegak
    \min \left\{
    \frac{\omegak}{\CCnorm^4 \CCnorm^2 \CChessF},
    \frac{1-\CCDelta}{\CCnorm^2}\Deltak,
    \frac{1}{\CCnorm^2}
    \right\}
    \\
    &\qquad\stackrel{\mathclap{\eqref{eqn:critical_values_equivalent}}}=
    \frac{(1 - \CCDelta)\CCsdT}{\CCcrit\CCnorm^2}\omegak_2
    \min \left\{
    \frac{\omegak_2}{(1 - \CCDelta) \CCcrit\CCnorm^4 \CChessF},
    \Deltak,
    \frac{1}{(1 - \CCDelta)}
    \right\}
    ,
\end{align*}
and~\eqref{eqn:sufficient_decrease} follows with
$\CCsd := \dfrac{\CCsdT (1 - \CCDelta)}{\CCcrit\CCnorm^2} \in (0,1)$
and from the fact that $\nicefrac{1}{(1 - \CCDelta)} \ge 1$
and with $\CCcritDenom := \CCcrit\CCnorm^4 \CChessF \ge 1$, which we may assume \ac{wlog}

\bibliographystyle{plain+eid}

%  This inserts the bib file
\bibliography{references.bib}

\end{document}